\documentclass[11pt,reqno]{amsart}

\usepackage[utf8]{inputenc} 
\usepackage[T1]{fontenc}
\usepackage[english]{babel}
\usepackage{csquotes}
\usepackage{xfrac}

\usepackage{tabularx}
\usepackage{array}
\usepackage{diagbox}
\usepackage[table]{xcolor}

\usepackage{aliascnt}

\usepackage{amsfonts, amsmath, amsthm, amssymb, pdfsync, extarrows, bigints}
\usepackage{mathalpha, mathrsfs, MnSymbol}

\usepackage[title]{appendix}
\usepackage{esint}
\usepackage{a4wide}
\usepackage{graphicx}
\usepackage{float}
\usepackage{color}
\usepackage{microtype}
\usepackage{indentfirst}
\usepackage{lettrine}
\usepackage{setspace}

\usepackage{imakeidx}
\usepackage{enumitem}
\usepackage{cancel}

\usepackage[colorlinks=true, allcolors=blue]{hyperref}
\usepackage{cleveref}

\usepackage[backend=biber,style=alphabetic]{biblatex}
\numberwithin{equation}{section}
\usepackage{xfrac}

\newtheorem{theo}{Theorem}[section]

\newaliascnt{lemma}{theo}
\newtheorem{lemma}[lemma]{Lemma}
\aliascntresetthe{lemma}

\newaliascnt{prop}{theo}
\newtheorem{prop}[prop]{Proposition}
\aliascntresetthe{prop}

\newaliascnt{cor}{theo}

\aliascntresetthe{cor}

\newaliascnt{dhef}{theo}
\newtheorem{dhef}[dhef]{Definition}
\aliascntresetthe{dhef}

\newaliascnt{ohss}{theo}
\newtheorem{ohss}[ohss]{Remark}
\aliascntresetthe{ohss}

\newaliascnt{example}{theo}

\aliascntresetthe{example}

\newaliascnt{ex}{theo}

\aliascntresetthe{ex}

\newaliascnt{eser}{theo}

\aliascntresetthe{eser}

\newaliascnt{cohmm}{theo}

\aliascntresetthe{cohmm}

\crefname{theo}{Theorem}{Theorems}
\crefname{lemma}{Lemma}{Lemmas}
\crefname{prop}{Proposition}{Propositions}
\crefname{cor}{Corollary}{Corollaries}
\crefname{dhef}{Definition}{Definitions}
\crefname{ohss}{Remark}{Remarks}
\crefname{example}{Example}{Examples}
\crefname{ex}{Example}{Examples}
\crefname{eser}{Exercise}{Exercises}
\crefname{cohmm}{Important Comment}{Important Comments}

\def\N{\mathbb{N}}
\def\Z{\mathbb{Z}}

\def\R{\mathbb{R}}
\definecolor{headercolor}{RGB}{220,235,245}
\definecolor{columncolor}{RGB}{240,240,240}

\definecolor{magenta_custom}{rgb}{1,0,1} %

\newcommand{\mres}{\mathbin{\vrule height 1.6ex depth 0pt width 0.13ex\vrule height 0.13ex depth 0pt width 1.3ex}}

\def\car#1{\raise1pt\hbox{$\chi$}_{#1}}
 \def\norm#1{ \|#1 \|}

\newcommand{\union}{\Z \cup S}

\usepackage{mathtools}
\mathtoolsset{showonlyrefs}

\title{Non-local phase transition models
for partial singularities}
\date{}

\author[A. Garroni]{Adriana Garroni}
	\address {Adriana Garroni \newline \indent
		Dipartimento di Matematica, Università La Sapienza di Roma, 
        \newline \indent Piazzale Aldo Moro 5, 00185 Roma, ITALY}
    \email{\href{mailto:garroni@mat.uniroma1.it}{garroni@mat.uniroma1.it}}
	\author[M. A. Palladino]{Maria Antonietta Palladino}
	\address {Maria Antonietta Palladino \newline \indent
		Classe di Scienze, Scuola Normale Superiore,\newline \indent
		Piazza dei Cavalieri 7, I-56126 Pisa, ITALY}
	\email{\href{mailto:mariaantonietta.palladino@sns.it}{mariaantonietta.palladino@sns.it}}

\begin{document}
\begin{abstract}
    \noindent We study  $\Gamma$-convergence, under appropriate scalings, of the energy functionals
\[
F_\varepsilon(u) := [u]_{H^s(Q)}^2 + \int_{Q} W_\varepsilon(u)\,dx,
\]
where $s \in (0,1)$, $[\cdot]_{H^s(Q)}$ denotes the $H^s$-seminorm in $Q\subset \R^{n}$, and $W_\varepsilon$ is a multi-well potential featuring infinitely many equally spaced wells together with additional intermediate wells whose height vanishes at a prescribed rate. Depending on this rate, different limiting behaviors emerge, ranging from classical models to new intermediate regimes. 
\end{abstract}
\maketitle
\tableofcontents
\section{Introduction}

In this paper, we study the asymptotics via $\Gamma$-convergence of a general class of non-local diffuse interface models for phase transitions with infinitely many stable states. Our prototype example is the functional
\begin{equation}\label{garroni}
\frac{1}{|\log \varepsilon|}
\iint_{Q \times Q} 
\frac{|u(x) - u(y)|^2}{|x-y|^{n+1}} \, dx \, dy
+
\frac{1}{\varepsilon |\log \varepsilon|}
\int_Q W_\varepsilon(u(x)) \, dx,
\end{equation}
where $n\in \N\setminus \{0\}$ and $Q$ is a cube in $\mathbb{R}^n$. The first term is the square of the $H^{1/2}$ Gagliardo seminorm, and the second one is a potential that vanishes exactly on the integers, i.e.\\ $\{t\,:\,W_\varepsilon(t)=0\}=\mathbb{Z}$.

This energy shows the typical structure of variational models for phase transitions. It consists of a singular perturbation of a multi-well potential (the prototype being the Cahn--Hilliard functional \cite{Cahn-Hilliard}) whose global minima represent the pure phases, which are the energetically preferred configurations of the system. After rescaling the Cahn--Hilliard models one obtains the Modica--Mortola functional (see \cite{Modica-Mortola, Modica}) and in the limit recovers a sharp interface model. In the limit regime, only pure phases are allowed, and the effective energy concentrates on the jump set: the interface that separates different phases. 

This variational scheme applies to a wide range of problems in physics and materials science. An extensive body of literature has explored variants of the classical Modica--Mortola functional to account for anisotropy (see, e.g., \cite{Owen-Sternberg, Bouchitte}), multiphase vectorial systems \cite{Baldo}, and homogenization in perforated or heterogeneous media (see \cite{Anzellotti-Baldo-Percivale, Caffarelli-de-la-Llave}, \cite{Fonseca-Tartar} and references therein).
Parallel to these gradient-based developments, a significant line of research has focused on phase transition models with non-local regularizations, where the gradient term is replaced by convolution-type \cite{A3} or fractional interaction kernels (\cite{Alberti-Bouchitte-Seppecher}). This framework has been further extended from an analytical perspective to phase transitions governed by arbitrary fractional operators, such as Gagliardo-type seminorms (\cite{Savin,Palatucci}).

In this setting one can also formulate the classical Nabarro--Peierls energy for dislocations (see \cite{Koslowski}), leading to the development of phase-field models for dislocations driven by non-local \(H^{1/2}\)-type regularizations, as in \eqref{garroni}, with $W_\varepsilon(t)=\operatorname{dist}^2(t, \mathbb{Z})$. Here the periodic potential reflects the underlying crystalline structure: the order parameter $u:Q \to \mathbb{R}$ describes the plastic slips of the form  $u\,\mathbf{e}_1$, while the $H^{1/2}$ seminorm represents the elastic energy induced by the slip $u$ along a single slip plane. 
In the limit this model produces a line-tension model for crystal plasticity (see, e.g.,  \cite{Garroni-Muller}, see also \cite{Kurzke, Kurzke2} for the case $n=1$ with applications to boundary vortices). Generalization of this model to vector-valued slips and to the case of multiple slip planes can be found in \cite{Conti-Garroni-Muller}.

In this paper, motivated by the study of partial dislocations and the emergence of stacking faults in crystals, we introduce a novel mathematical ingredient into this framework. Specifically, we analyze the asymptotics of a singular, non-local perturbation of a multi-well potential - in the spirit of the aforementioned phase-field models for dislocations \cite{Garroni-Muller}. A key feature of our model is that the potential exhibits stable local minima. These minima represent ``intermediate'' metastable states that are energetically permissible, though less favorable than the global ones. In the mechanics literature, the energy profile mapping these configurations is known as the stacking fault energy or \(\gamma \)-surface \cite{Vitek}. 

For concreteness, we consider a potential $W_\varepsilon:\R \rightarrow [0,+\infty)$ of the form
\begin{equation}\label{intro_We}
W_\varepsilon(u) := \min \left\{ \operatorname{dist}^2\!(u, \mathbb{Z}), \, \operatorname{dist}^2\!\left(u, \mathbb{Z}+\tfrac{1}{2}\right) + \delta_\varepsilon \right\},
\end{equation}
where $\{\delta_\varepsilon\}\subset \R^+ $ represents the energy gap (or ``cost'') associated with the secondary wells. The effect of intermediate wells depends on how fast the gap $\delta_\varepsilon$
vanishes as $\varepsilon\to 0$. Intuitively, the asymptotic height of the intermediate well determines whether the limiting energy allows partial slips and, if so, at what energetic cost. The goal of this work is to identify the relevant asymptotic regimes corresponding to different limiting behaviors. In particular, in the case $\delta_\varepsilon=\varepsilon|\log\varepsilon|$, we prove that the limit energy of \eqref{garroni}, choosing $W_\varepsilon$ as in \eqref{intro_We}, is 
    \begin{equation}\label{limit-intro}
  \omega_{n-1} \int_{J_u} |u^+- u^-| \,d\mathcal{H}^{n-1} +  \left|\left\{ u\in \mathbb{Z}+\tfrac12\right\}\right| \qquad \text{for } u \in BV\left(Q; \tfrac{1}{2}\mathbb{Z}\right),
    \end{equation}
   where $\omega_{n-1}$ is the volume of the unit ball in $\mathbb{R}^{n-1}$ ($\omega_0=1$).
Therefore, in this regime, the permanence in intermediate states has a cost which is captured by the limit energy and shows the emergence of the stacking fault term, which, in the context of dislocations, represents an energy concentrated on the regions of the slip plane separating partial dislocation lines. This is the main result of this paper.

It is worth noting that variational models involving partial defects have attracted significant interest in recent years. On one hand models similar in spirit to the one proposed here have been considered in the discrete setting for two-dimensional systems (see \cite{Badal-Cicalese-DeLuca-Ponsiglione} for a generalized XY model and \cite{Bach-Cicalese-Garroni-Orlando} for a next-to-nearest neighbor lattice model for partial edge dislocations) where the asymptotic energy concentrates on fractional singularities (vortices or dislocations) connected by line-like string defects. On the other hand, several continuum models share a similar analytical structure, where the coupling between fractional point singularities and line discontinuities is topologically induced, leading to a structure where integer-degree defects split into fractional components connected by penalized interfaces (see \cite{Goldman-Merlet-Millot}).

\bigskip

For completeness we study the asymptotics in terms of $\Gamma$-convergence of a more general problem in which we consider all possible regimes for the parameter $\delta_\varepsilon$, we assume that the local minima may have a more general structure (not necessarily at half integers) and we also include a more general non-local regularization which covers all fractional seminorms $H^s$.

\section{Formulation of the problem and main result}

Let $Q \subset \mathbb{R}^n$ be an $n$-dimensional cube with $n\ge 1$. For $s \in \left[\tfrac12, 1\right)$, we investigate the asymptotic behavior, as $\varepsilon$ goes to $0$, of the family of functionals $F_\varepsilon: L^1(Q) \to [0, +\infty]$ defined as follows:
\begin{equation}\label{eq:main_functional}
F_\varepsilon(u) := 
\begin{cases} 
\displaystyle \frac{1}{\theta_\varepsilon} \left( \varepsilon^{2s} \int_{Q} \int_{Q} \frac{|u(x) - u(y)|^2}{|x - y|^{n + 2s}} \, dx \, dy + \int_{Q} W_\varepsilon(u(x)) \, dx \right), & \text{if } u \in H^s(Q), \\[1em]
+\infty, & \text{otherwise,}
\end{cases}
\end{equation}
where the scaling factor $\theta_\varepsilon$ is chosen as in \cite{SavinValdinoci2010}:
\begin{equation}\label{eq:theta_scaling}
\theta_\varepsilon := 
\begin{cases} 
\varepsilon |\log \varepsilon|, & \text{if } s = \tfrac12, \\[1em] 
\varepsilon, & \text{if } s \in \left(\tfrac12, 1\right). 
\end{cases}
\end{equation}

The non-local term represents the square of the Gagliardo $H^s$-seminorm of $u$ in $Q$. 
For the sake of brevity, we shall denote it simply by
\[
[u]_{H^s(Q)}^2:=\int_{Q} \int_{Q} 
\frac{|u(x) - u(y)|^2}{|x - y|^{n + 2s}} \, dx \, dy.
\]

The potential $W_\varepsilon: \mathbb{R} \to [0, +\infty)$ is designed to model a system with multiple hierarchical states. It vanishes on $\mathbb{Z}$ and has local minima at $S := \{s_z\}_{z \in \mathbb{Z}}$ representing partial slips or intermediate stable configurations. To ensure that the wells are well-separated, we assume that there exists a constant $d\in \left( 0,\tfrac12\right]$ such that:
\begin{equation}\label{d}
z + d \le s_z \le z + 1 - d \quad \text{for every } z \in \mathbb{Z}.
\end{equation} The energetic cost at these intermediate configurations is $\delta_\varepsilon\geq 0$ with
\begin{equation}\label{c}
\lim_{\varepsilon\to 0}\delta_\varepsilon=0.
\end{equation}
For concreteness, the potential $W_\varepsilon$ is defined by
\begin{equation}\label{eq:potential_def}
W_\varepsilon(u) := \min \left\{ \operatorname{dist}^2\!(u, \mathbb{Z}), \, \operatorname{dist}^2\!(u, S) + \delta_\varepsilon \right\}.
\end{equation}
 The rate at which $\delta_\varepsilon$ 
vanishes plays a crucial role in our analysis giving rise to different 
variational regimes.
\begin{figure}[htbp]
    \centering
    \includegraphics[trim=0cm 0cm 0cm 0.9cm, clip, width=0.8\textwidth]{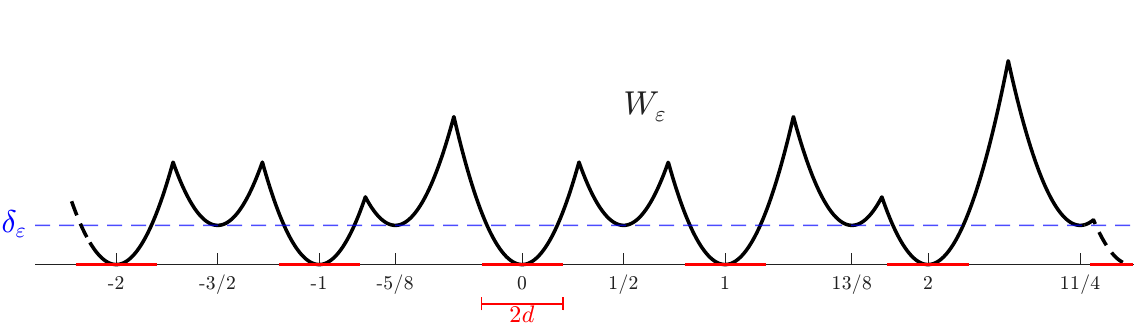}
    \caption{Example of potential $W_\varepsilon$}
    \label{fig:potenziale_W}
\end{figure}
We study how the competition between the fractional interaction term and 
the penalization of the intermediate wells through $\delta_\varepsilon$ 
determines both the set of admissible limit values (the state space) 
and the form of the $\Gamma$-limit.

We distinguish two different regimes. We assume that the relevant limits exist, and classify the regimes as follows:
\begin{itemize}
    \item[(R1)]\hypertarget{R1}{ \textbf{Supercritical/Critical regime}: We are in the \textit{supercritical} or \textit{critical regime} if
\begin{equation}\label{eq-super-regime}
    K := \lim_{\varepsilon \to 0} \frac{\delta_\varepsilon}{\theta_\varepsilon} \in [0, +\infty).
\end{equation}
We refer to the case $K = 0$ as the \emph{supercritical} regime, and to the case $K \in (0, +\infty)$ as the \emph{critical} one.}

    \item[(R2)] \hypertarget{R2}{\textbf{Subcritical regime}: We are in the \textit{subcritical regime} if
    \begin{equation}\label{eq-sub-regime}
        \lim_{\varepsilon \to 0} \frac{\delta_\varepsilon}{\theta_\varepsilon} = +\infty.
    \end{equation}}
    In the case $s=\tfrac12$ and regime \hyperlink{R2}{(R2)}, we additionally assume
\begin{equation}\label{ell}
\ell :=
\lim_{\varepsilon \to 0}
\frac{|\log \delta_\varepsilon|}{|\log \varepsilon|}
\in [0,1].
\end{equation}
\end{itemize}
As we shall see, these regimes lead to different target spaces for the limit functions. In regime \hyperlink{R1}{(R1)}, the intermediate wells are ``cheap'' enough to be visited in the limit. In this case, the state space of the limit function is $BV(Q; \mathbb{Z} \cup S)$.  Conversely, in \hyperlink{R2}{(R2)}
the intermediate wells become energetically too expensive. Consequently, the limit function $u$ is forced to take values in $\mathbb{Z}$.

\medskip

\subsection{The limiting jump costs}
In the sharp interface limit only pure phases are present: the limit $u$ takes values in $\union$ in regime \hyperlink{R1}{(R1)} and in $\mathbb{Z} $ in \hyperlink{R2}{(R2)}. 
Moreover, the limit energy is concentrated on the jump set $J_u$ with density $\Psi(u^+,u^-)$, where $u^+, u^-\in \union$ denote the values of $u$ around a jump point.
We refer to $\Psi(\cdot,\cdot)$ as the \textit{jump cost}. In all cases, it is the product of three factors
\[
\Psi(\cdot,\cdot)
=
\underbrace{\lambda_{s,n}}_{\substack{\text{Dimensional}\\ \text{const.}}}
\cdot
\overbrace{\kappa_s}^{\mathclap{\substack{\text{Transition profile}\\ \text{const.}}}}
\cdot
\underbrace{\Phi(\cdot,\cdot)}_{\substack{\text{Position}\\ \text{of the wells}}}
\]
that we define as follows.
\begin{dhef}[Position of the wells]\label{one-dim_jump}
     Let $\mathcal{V}$ be a discrete set in $\R$.
        For  $a< b $ in $\mathcal{V} $, let $a = p_0 < p_1 < \dots < p_N = b$ be the points of $\mathcal{V}\cap [a, b]$. We define
\begin{equation}\label{dhef:critical_jump}
\Phi_\mathcal{V}(a, b) := \sum_{i=1}^N (p_i - p_{i-1})^2,
\end{equation}
extended by $\Phi_\mathcal{V}(a, b)=\Phi_\mathcal{V}(b, a)$ and $\Phi_\mathcal{V}(a, a)=0.$ 
\end{dhef}
\begin{dhef} [Transition profile constant]
We set $\kappa_{1/2}=2$ and for $s \in\left(\tfrac12,1\right),$
\begin{equation}\label{eq:kappa_s}
\kappa_s:=\inf\left\{
[w]^2_{H^s(\R)}+\int_\R \operatorname{dist}^2(w,\{0,1\})\,dx\,:\,w \in H^s_{\rm loc}(\R), \,\lim_{x \to -\infty }w(x)=0,\,\lim_{x\to +\infty}w(x)=1
\right\},
\end{equation}
where
$$
[w]^2_{H^s(\R)}=\int_\R\int_\R\frac{|w(x)-w(y)|^2}{|x-y|^{1+2s}}\,dx\,dy.
$$
\end{dhef}
\noindent The constant $\kappa_s$ is the one-dimensional cost of a transition from $0$ to $1$ for the standard double-well potential $\operatorname{dist}^2(\cdot,\{0,1\})$: if $s>\tfrac12$, it is given by the optimal profile problem (see, e.g., \cite{Palatucci,SavinValdinoci2010}); if $s=\tfrac12$, with the logarithmic scale, the cost does not depend on the profile, and it is equal to $2$ (see, e.g., \cite{ABS1,Garroni-Muller,Kurzke}).
\begin{dhef}[Dimensional constant]
     For $n \ge 1$ and $s\in \left[\tfrac12,1\right)$, we set 
\begin{equation}\label{dhef:C_ns}
C_{s,n} := 
\pi^{\frac{n-1}{2}}
\frac{\Gamma\!\left(s+\frac12\right)}
{\Gamma\!\left(\frac{n+2s}{2}\right)}
=\begin{cases}
    \displaystyle \int_{\mathbb{R}^{n-1}} (1+|t|^2)^{-\frac{n+2s}{2}} \, dt, & \text{if } n \ge 2, \\[10pt]
    1, & \text{if } n=1
\end{cases}\qquad \text{and}\qquad \lambda_{s,n}:=(C_{s,n})^{\frac{1}{2s}}.
\end{equation}
In particular, with this definition we have
$$
\lambda_{s,1}=1,\,\quad \lambda_{1/2,n}=C_{1/2,n}=\omega_{n-1},
$$
where $\omega_{n-1}$ is the volume of the unit ball of $\R^{n-1}$ (with $\omega_0=1$).
\end{dhef}
Based on the above definitions, we can define the jump costs.

\medskip

\begin{dhef}[Jump costs]\label{Jump_costs}
    For $n\ge 1$ and $s\in \left[\tfrac12,1\right)$,
we define
    $$
\Psi_{\mathrm{sup}}:(\mathbb Z\cup S)\times(\mathbb Z\cup S)\rightarrow [0,+\infty),\qquad \Psi_{\mathrm{sub}}:\mathbb Z\times\mathbb Z\to [0,+\infty),
$$
defined by 
\begin{equation}\label{eq:psi-sup-sub}
\Psi_{\mathrm{sup}} :=\lambda_{s,n}\,\kappa_s\, \Phi_{\Z\cup S},
\qquad 
\Psi_{\mathrm{sub}} :=\lambda_{s,n}\,\kappa_s\cdot 
\begin{cases}
(1-\ell)\,\Phi_\Z+\ell \,\Phi_{\Z\cup S}, & \text{if }s=\tfrac12\\
\Phi_{\Z \cup S}, & \text{if }s>\tfrac12,
\end{cases}
\end{equation}
where $\ell \in [0,1]$ is the parameter in \eqref{ell}.
\end{dhef}
Observe that, for $s>\tfrac12$, the cost $\Psi_{\rm sub} $ is the restriction of $\Psi_{\rm sup} $ to integer phases; if instead $s=\tfrac12$, it is the interpolation of the cost with the wells in $\Z$ and in $\Z \cup S$ (see \Cref{remark_interpolazione}).
Moreover, both jump costs are additive with respect to the order of the wells, namely
\[
\Psi_{\mathrm{sup}}(\alpha,\gamma)=\Psi_{\mathrm{sup}}(\alpha,\beta)+\Psi_{\mathrm{sup}}(\beta,\gamma)
\qquad \text{for } \alpha<\beta<\gamma \hbox{ in } \Z \cup S,
\]
\[
\Psi_{\mathrm{sub}}(\alpha,\gamma)=\Psi_{\mathrm{sub}}(\alpha,\beta)+\Psi_{\mathrm{sub}}(\beta,\gamma)
\qquad \text{for } \alpha<\beta<\gamma \hbox{ in } \Z. 
\] 
This means that a jump between non-consecutive phases costs exactly as much as a \textit{staircase} of elementary (admissible) jumps visiting all intermediate phases. Therefore, it is sufficient to characterize the jump costs between consecutive wells in $\mathbb{Z}\cup S$ and $\mathbb{Z}$, respectively.  Let $a<b$ be two consecutive elements of $\mathbb{Z}\cup S$, and let $k\in\mathbb{Z}$. By \Cref{Jump_costs}, we have
\begin{equation}\label{eq_sup:psi-consecutive}
\Psi_{\mathrm{sup}}(a,b)=\begin{cases}
\omega_{n-1}\,2\,(b-a)^2, & \text{if }s=\frac12,\\[1em]
\lambda_{s,n}\,\kappa_s\,(b-a)^2, & \text{if }s>\frac12,
\end{cases}
\end{equation}
\begin{equation}\label{eq_sub:psi-consecutive}
\Psi_{\mathrm{sub}}(k,k+1)=
\begin{cases}
\omega_{n-1}\,2\,\big(1-\ell\,2w_k(1-w_k)\big), & \text{if }s=\frac12,\\[1em]
\lambda_{s,n}\,\kappa_s\,\big(1-2w_k(1-w_k)\big), & \text{if }s>\frac12
\end{cases}\qquad \text{with}\qquad w_k:=s_k-k,
\end{equation}
where for $s=\tfrac12$ we used $\lambda_{1/2,n}=\omega_{n-1}$. 

\subsection{The \texorpdfstring{$\Gamma$}{Gamma}-convergence theorem}
We can now state the main results of the paper, which provide a complete characterization of the $\Gamma$-limit in the different regimes. 

Throughout, for $u\in BV(Q;\mathcal I)$ with $\mathcal I\subset\mathbb R$ 
discrete, $J_u$ denotes the (approximate) jump set, $u^\pm$ the one-sided 
traces on $J_u$ with the orientation chosen so that $u^+>u^-$, $[u]:=u^+-u^-$ the jump 
amplitude, and $\mathcal H^{n-1}$ the $(n-1)$-dimensional Hausdorff measure (see \cite{Ambrosio}). 

We write $|A|$ for the Lebesgue 
measure of a measurable set $A\subset\mathbb R^n$.

\begin{theo}[Compactness]\label{thm:compactness}
Let $s \in \left[\tfrac12,1\right)$ and let $\{u_\varepsilon\} \subset L^1(Q)$ satisfy
$$\sup_{\varepsilon} F_\varepsilon(u_\varepsilon) < +\infty.$$
\\\noindent
Then there exists a sequence of integers $\{a_\varepsilon\} \subset \mathbb{Z}$ such that the translated sequence
$\{u_\varepsilon - a_\varepsilon\}$ is relatively compact in $L^1(Q)$.

Moreover, for small $\varepsilon$, the sequence
$\{u_\varepsilon - a_\varepsilon\}$ is uniformly bounded in 
$L^{\frac{n}{n-1}}(Q)$ if $n \ge 2$, and  in 
$L^{q}(Q)$ for every 
$ q \in[1,+\infty) $ if $n=1$; if $n=1$ and $s > \tfrac12$, it is uniformly bounded in $L^\infty(Q)$.
\\\noindent
Finally, if $u_\varepsilon \to u$ in $L^1(Q)$, then $u \in BV(Q;\Z\cup S)$. In regime \hyperlink{R2}{(R2)}, we have  $u \in BV(Q;\Z)$.
\end{theo}

\begin{ohss}\label{rem:compactness}
The limit points of $\{u_\varepsilon-a_\varepsilon\}$ need not take values in $\Z\cup S$ if $S$ is not periodic, because $F_\varepsilon$ is not invariant under integer translations. For instance, take $u_\varepsilon\equiv s_{k_\varepsilon}$ with $k_\varepsilon\to\infty$. In general these limit points take values in $\Z\cup S^*$, where $S^*:=\{
\lim_\varepsilon s_{z+a_\varepsilon}-
a_\varepsilon:\ z\in\Z\}$ (see the proof in \Cref{sottosezione_compattezza}). 
\end{ohss}
\begin{theo}[$\Gamma$-convergence]\label{thm:main}
Let $s \in \left[\frac12, 1\right)$ and $Q\subset \R^n$ be an open cube. As $\varepsilon\to 0$, the sequence of functionals $F_\varepsilon$ in \eqref{eq:main_functional} $\Gamma$-converges, with respect to the strong $L^1(Q)$ topology, to the following limits.
\begin{itemize}
    \item \textbf{Regime \hyperlink{R1}{(R1)}:} If $\delta_\varepsilon$ satisfies \eqref{eq-super-regime} with $K \in [0,+\infty)$, the $\Gamma$-limit is
    \begin{equation}
    F_{\mathrm{sup}}(u) := \begin{cases}
    \displaystyle \int_{J_u} \Psi_{\mathrm{sup}}(u^+, u^-) \, d\mathcal{H}^{n-1} + K |\{x \in Q : u(x) \in S\}|, & \text{if } u \in BV(Q; \union), \\[1em]
    +\infty, & \text{otherwise,}
    \end{cases}
    \end{equation}
   where $\Psi_{\mathrm{sup}}$ is the jump cost in \Cref{Jump_costs}. 
        \item \textbf{Regime \hyperlink{R2}{(R2)}:} If $\delta_\varepsilon$ satisfies \eqref{eq-sub-regime}, the $\Gamma$-limit is 
    \begin{equation}
    F_{\mathrm{sub}}(u) := \begin{cases}
    \displaystyle \int_{J_u} \Psi_{\mathrm{sub}}(u^+, u^- )\, d\mathcal{H}^{n-1}, & \text{if } u \in BV(Q; \mathbb{Z}), \\[1em]
    +\infty, & \text{otherwise,}
    \end{cases}
    \end{equation}
     where $\Psi_{\mathrm{sub}}$ is the jump cost in  \Cref{Jump_costs}.
\end{itemize}

\end{theo}
It should be noted that the two regimes differ in the set of admissible phases and in the role of the secondary wells $S.$
\begin{itemize}
    \item In regime \hyperlink{R1}{(R1)}, the phases in $S$ are admissible and each has a volume cost $K$;

    \item  In regime \hyperlink{R2}{(R2)}, the phases in $S$ are \textit{not} admissible; however, the set $S$ affects the jump cost depending on $s$:
    \begin{itemize}
    \item{ If $s>\tfrac12,$ the cost of the jump $k \to k+1$ is equal to the cost of a double jump $k \to s_k\to k+1,$ passing through the intermediate well $s_k$;}
        
        \item If $s=\tfrac12,$ also the parameter $\ell$ plays a role: it represents the interpolation between a jump that ignores $s_k$ and the one that passes through it (see \Cref{remark_interpolazione}).
    \end{itemize}

\end{itemize}

For future reference, given a discrete set $\mathcal{V}\subset \R$ and an open domain $A\subseteq\R^n$, we define 
\begin{equation}\label{eq:periodic_general} F^{\mathcal{V}}_\varepsilon(u;A) := \frac{1}{\theta_\varepsilon} \left( \varepsilon^{2s} \int_{A} \int_{A} \frac{|u(x) - u(y)|^2}{|x - y|^{n + 2s}} \, dx \, dy + \int_{A} \operatorname{dist}^2(u(x), \mathcal{V}) \, dx \right),
\end{equation}
and we omit $A$ when  $A=Q$.
\begin{ohss}[Periodic potential]\label{no_buca}
\Cref{thm:main} also covers the $\Gamma$-limit of the functional 
$F_\varepsilon^\Z$, obtained from \eqref{eq:periodic_general} with $\mathcal{V}=\Z$.
This functional was studied for $s=\tfrac12$ in \cite{Kurzke} $(n=1)$ and \cite{Garroni-Muller} $(n=2)$. Its $\Gamma$-limit is
\begin{equation}\label{eq:remark_periodic}
\lambda_{s,n}\,\kappa_s\int_{J_u}[u]\,d\mathcal{H}^{n-1}, \quad u \in BV(Q;\Z),
\end{equation}
which, for $s=\tfrac12$, equals 
$2\omega_{n-1}\int_{J_u}[u]\, d\mathcal{H}^{n-1}$, 
that is precisely the result already obtained in \cite{Kurzke} and \cite{Garroni-Muller}.
The functional \eqref{eq:remark_periodic} can be obtained by a simple rescaling; indeed we have $$F^{\mathbb{Z}}_\varepsilon(u) = 4 \,F^{\frac12\mathbb{Z}}_\varepsilon\left(\frac{u}{2}\right),$$ where 
$F^{\frac12\mathbb{Z}}_\varepsilon$
is the functional in \eqref{eq:main_functional} with  $S = \mathbb{Z}+\frac12$ and $\delta_\varepsilon = 0$, i.e. the regime \hyperlink{R1}{(R1)} with $K=0$. 
Hence, the $\Gamma$-limit of $F^{\frac12\mathbb{Z}}_\varepsilon(\cdot/2)$ is, by \Cref{thm:main}, the following 
$$ \lambda_{s,n}\,\kappa_s\int_{J_u}\Phi_{\frac12\Z}(u^+/2,u^-/2) \, d\mathcal{H}^{n-1}=\lambda_{s,n}\,\kappa_s\int_{J_u}\frac{[u]}{4}\,d\mathcal{H}^{n-1}, \quad \quad
u \in BV(Q; \mathbb{Z}).$$
\end{ohss}

\begin{ohss}[Role of $\ell$ for $s=\tfrac12$]\label{remark_interpolazione} In the subcritical regime \hyperlink{R2}{(R2)} for $s=\tfrac12,$ the limit energy \( F_{\mathrm{sub}} \)
 can be interpreted as an interpolation, through the parameter $\ell$, between the two extreme cases: when the potential wells are only at integers, and when wells occur in $\union$. 
Indeed, by definition, we have 
\begin{equation}\label{eq:interpolazione}
F_{\mathrm{sub}}(u)=(1-\ell)\left[2\,\omega_{n-1}\int_{J_u}[u]\,d\mathcal H^{n-1}\right]+\ell\,\left[2\,\omega_{n-1}\int_{J_u}\Phi_{\Z\cup S}(u^-,u^+)\,d\mathcal H^{n-1}\right],\text{ for }u\in BV(Q;\Z),
\end{equation}
where the first term is the $\Gamma$-limit of $F^\Z_\varepsilon$ (as we observe in \Cref{no_buca}) and the second is the $\Gamma$-limit of $F^{\Z\cup S}_\varepsilon$ (as follows from \Cref{thm:main} choosing $\delta_\varepsilon=0).$
\end{ohss}
Let us notice that the jump cost appearing in the limit functional has a different
nature according to $s$.
\begin{itemize}
    \item For $s>1/2$ the energy of a transition between two consecutive wells is 
the minimal energy of a one-dimensional profile connecting them: the cost $\kappa_s$ in \eqref{eq:kappa_s} is defined by an optimal profile problem on $\R$;
\item  For $s=1/2$, there is no optimal profile; indeed, because of the logarithmic scaling, the cost depends only on the end states: the cost is constant, $\kappa_s=2$.
\end{itemize}
 In both cases, the
$n$-dimensional problem reduces to a one-dimensional one across the
interface, and the dimension only enters through a multiplicative
constant $\lambda_{s,n}$. The next remark makes this reduction explicit and explains
how the constant affects the cost in each case.
\begin{ohss}[The dimensional constant $\lambda_{s,n}$]\label{rmk:slicing_constant}
The constant $C_{s,n}$ represents the energy of a planar transition per unit
area of interface; indeed, fix a one-dimensional function $v :\R \to \R$ and $\nu \in \mathbb{S}^{n-1}$. Then, for every $E \subset \nu^\perp $, we have
$$
\int_{\mathcal{C}_E}\int_{\R^n}\frac{|v(x\cdot \nu)-v(y\cdot \nu)|^2}{|x-y|^{n+2s}}\,dx\,dy=C_{s,n}[v]^2_{H^s(\R)}\, \mathcal{H}^{n-1}(E),
$$
with $\mathcal{C}_E=\{z+t\nu:z\in E,\,t \in \R\}$ the cylinder with base $E$ in the direction $\nu$. Thus, on profiles depending only on $x\cdot\nu$, the $n$-dimensional energy is the one-dimensional one with the weight $C_{s,n}$ in front of the seminorm.
How this weight enters the cost depends on $s.$ 
\begin{itemize}
    \item{ For $s>\tfrac12, $ the two terms of the one-dimensional energy scale differently: for $L>0 $ and $v_L(t):=v(t/L)$, we have
    $$
     \quad \quad C_{s,n}[v_L]^2_{H^s(\R)}+\int_\R \operatorname{dist}^2(v_L,\Z\cup S)= C_{s,n}L^{1-2s}[v]_{H^s(\R)}^2+L\int_\R \operatorname{dist}^2(v,\Z\cup S).
    $$
    Choosing $L=\lambda_{s,n}$ the two coefficients coincide, since $C_{s,n}\lambda_{s,n}^{1-2s}=\lambda_{s,n}$. 
    
    Thanks to this observation, given two consecutive points $a<b$ of $\union$, we can compute the optimal transition between them, i.e.
     \begin{align*}
     &\qquad\,\,\,\inf\Big\{C_{s,n}[v]_{H^s(\mathbb{R})}^2+\int_{\mathbb{R}}\operatorname{dist}^2(v,\mathbb{Z}\cup S)\,dx : v\in H^s_{\rm loc}(\mathbb{R}),\ v(-\infty)=a,\ v(+\infty)=b\Big\}\\
     &\qquad=\inf\Big\{C_{s,n}[v]_{H^s(\mathbb{R})}^2+\int_{\mathbb{R}}\operatorname{dist}^2(v,\{a,b\})\,dx : v\in H^s_{\rm loc}(\mathbb{R}),\ v(-\infty)=a,\ v(+\infty)=b\Big\}.
     \end{align*}
Let $w\in H^s_{\rm loc}(\mathbb{R})$ be an admissible profile from $0$ to $1$ and set
$v(x):=a+(b-a)\,w\!\left(\frac{x}{\lambda_{s,n}}\right).$
The map $w\mapsto v$ is a bijection between admissible profiles from $0$ to $1$ and those from $a$ to $b$. Since
\[
\big[w(\cdot/\lambda_{s,n})\big]_{H^s(\mathbb{R})}^2=\lambda_{s,n}^{1-2s}[w]_{H^s(\mathbb{R})}^2,
\qquad C_{s,n}\lambda_{s,n}^{1-2s}=\lambda_{s,n}
\]
and $\operatorname{dist}(v,\{a,b\})=(b-a)\operatorname{dist}(w(\cdot/\lambda_{s,n}),\{0,1\})$, by the definition of $\kappa_s$ in \eqref{eq:kappa_s} the infimum above equals
\begin{align*}\qquad
&\lambda_{s,n}(b-a)^2\inf\Big\{[w]_{H^s(\mathbb{R})}^2+\int_{\mathbb{R}}\operatorname{dist}^2(w,\{0,1\})\,dx : \\
&\hspace{6em} w\in H^s_{\rm loc}(\mathbb{R}),\ w(-\infty)=0,\ w(+\infty)=1\Big\}
=\lambda_{s,n}(b-a)^2\kappa_s =\Psi_{\rm sup}(a,b).
\end{align*}
   In other words, the optimal profile in dimension $n$ is the one-dimensional one stretched by the factor $\lambda_{s,n}$, and its cost is the one-dimensional one multiplied by $\lambda_{s,n}$.
   The general computation for $a,b\in \union$ (not necessarily consecutive) is in \Cref{Section:optimal} where we show the additivity of the cost of the optimal transition with respect to multiple jumps. }
    \item For $s=\tfrac12,$ we have $L^{1-2s}=1$, so stretching the profile does not change the seminorm and no rescaling can balance the two terms. However, at logarithmic scale, the energy concentrates in the non-local term, hence the cost is the one-dimensional one multiplied by the constant itself
    $$
    C_{1/2,n}=\omega_{n-1},
    $$
    the volume of the unit ball in $\R^{n-1}$.
\end{itemize}
\end{ohss}

\subsection{Crucial comparison estimate and strategy of the proof}

Besides the presence of the secondary wells and their effect in the different regimes, the present model shows two major difficulties: lack of coercivity of the potential and non-locality of the singular perturbation.

The lack of coercivity is clear: every integer constant has zero energy. Hence compactness can
only be expected up to integer translations, as in the periodic case studied in
\cite{Garroni-Muller,Kurzke} ($n=1,2$, $s=1/2$). However, in our setting the potential is not periodic 
and then the energy $F_\varepsilon$ is not invariant
under integer translations, so the limit of a translated sequence need not take values in
$\Z\cup S$ (see \Cref{rem:compactness}). 

A key ingredient in the analysis is the comparison  
between our original functional and a simpler reference functional $F_\varepsilon^{\Z\cup S}$, obtained from \eqref{eq:periodic_general} with $\mathcal{V}=\Z\cup S$, i.e.
\begin{equation}\label{F0}
F_\varepsilon^{\mathbb{Z} \cup S}(u) := \frac{1}{\theta_\varepsilon} \left( \varepsilon^{2s} \int_Q \int_Q \frac{|u(x) - u(y)|^2}{|x - y|^{n+2s}} \, dx \, dy + \int_Q \operatorname{dist}^2(u(x), \mathbb{Z} \cup S) \, dx \right).
\end{equation}
 Clearly 
$
F_\varepsilon^{\mathbb{Z} \cup S}(u) \leq F_\varepsilon(u)
.$ 
A more careful comparison allows us to isolate and 
quantify the specific contribution of the energy gap $\delta_\varepsilon$ in 
the asymptotic behavior of $F_\varepsilon$.
 
To see this, for any $\delta\in(0,1]$ we define the set
\begin{equation}
    \label{eq-set-comparison}
    A_\delta:=\{t\in\R:\, \operatorname{dist}(t,S)<\delta\operatorname{dist}(t,\Z)\},
\end{equation}
and note that if $\delta \in (0,1)$ and $\varepsilon$ is small enough (depending on $\delta$ and $d$), it holds
$$
A_\delta\subseteq \{t\in \R:\ \operatorname{dist}^2(t,S)+\delta_\varepsilon<\operatorname{dist}^2(t,\Z)\}\subseteq A_1.
$$
Moreover, by the very definition of the potential $W_\varepsilon$ we have
that $W_\varepsilon(t)=\operatorname{dist}^2(t,\Z)=\operatorname{dist}^2(t,\Z\cup S)$ if $t\notin A_1$ and $W_\varepsilon(t)\leq \operatorname{dist}^2(t,S)+\delta_\varepsilon=\operatorname{dist}^2(t,\Z\cup S)+\delta_\varepsilon$ if $t\in A_1$, with equality if $t\in A_\delta$ (for $\varepsilon$ small). 
This, for every given $\delta\in(0,1)$ and $\varepsilon>0$ small enough, leads to the following two-sided estimate
\begin{equation}\label{eq:two_sided_estimate}
F_\varepsilon^{\mathbb{Z} \cup S}(u) + \frac{\delta_\varepsilon}{\theta_\varepsilon} |A_{\delta}(u)| \leq F_\varepsilon(u) \leq F_\varepsilon^{\mathbb{Z} \cup S}(u) + \frac{\delta_\varepsilon}{\theta_\varepsilon} |A_{1}(u)|,
\end{equation}
where
\begin{equation}
      \label{eq-set-u-comparison}
      A_\delta(u):=\{x\in Q:\, u(x)\in A_\delta\}.
\end{equation}

The estimate \eqref{eq:two_sided_estimate} is a fundamental technical tool in our analysis. It effectively decouples the total energy into two distinct contributions:
\begin{enumerate}
    \item A \emph{fractional transition energy} ($F_\varepsilon^{\mathbb{Z} \cup S}$), which accounts for the sharp interface energy;
    \item A \emph{penalization term} proportional to the measure of the set where $u$ is close to the secondary wells in $S$.
\end{enumerate}
This decomposition clarifies the emergence of the two regimes: in the subcritical regime \hyperlink{R2}{(R2)}, the penalization factor $\delta_\varepsilon/\theta_\varepsilon$ diverges, acting as a constraint that forces the limit function to take values only in $\mathbb{Z}$. Conversely, in the critical/supercritical regime  \hyperlink{R1}{(R1)}, this factor remains bounded or vanishes, allowing the secondary wells to survive in the limit as a bulk contribution.

Another elementary tool, used both in the proof of compactness and in the lower bound, is a truncation inequality.
 Let $\mathcal V=\{v_i\}_{i\in \Z}\subset\R$ be a discrete set with $v_i<v_{i+1}$ for all $i \in \Z$, and let $T_it:=(t\vee v_i)\wedge v_{i+1}$ be the truncation between two consecutive values of $\mathcal{V}$. For all $t,t'\in\R$, we have $$\sum_{i\in \Z}|T_it-T_it'|^2\le|t-t'|^2\quad \text{and} \quad \sum_{i\in \Z}W(T_it)=W(t)\text{ for every $W\ge0$ vanishing on $\mathcal V$.}$$ Hence, we obtain
\begin{equation}\label{eq:truncation}
\sum_{i\in \Z}F^{\mathcal V}_\varepsilon(T_iu)\le F^{\mathcal V}_\varepsilon(u),
\end{equation}
where $F^{\mathcal V}_\varepsilon$ is defined in \eqref{eq:periodic_general}. Thus the energy of $u$ controls the sum of the energies of its truncations, and the analysis can often be reduced to a single transition between two consecutive points of $\mathcal V$.
\subsubsection*{Strategy of the proof}
The compactness result (\Cref{thm:compactness}) is proved in \Cref{sec:compact}. By the lower estimate
in \eqref{eq:two_sided_estimate} it is enough to establish the compactness result for a sequence uniformly bounded for the energy $F^{\Z\cup S}_\varepsilon$. Here, a bi-Lipschitz
change of variables mapping $\Z\cup S$ onto $\Z$ reduces the problem to the periodic energy
$F^\Z_\varepsilon$. For the latter, we prove uniform $L^q$ bounds up to integer translations, extending
\cite{Garroni-Muller,Kurzke}, and we deduce compactness by truncating between consecutive integers and
using the classical two-well results. Since $F_\varepsilon$ is not invariant under integer translations,
the limit takes values in $\Z\cup S$ (in $\Z$ in regime \hyperlink{R2}{(R2)}) only for sequences that
converge without translations.

In \Cref{sec:lower-bound} we establish the lower bound for the functional $F_\varepsilon$. 
By the truncation inequality \eqref{eq:truncation}, it is enough to consider a single transition between two consecutive wells of the limit. In regime \hyperlink{R1}{(R1)}, the lower estimate \eqref{eq:two_sided_estimate} allows us to 
distinguish the contribution of the fixed wells in $\mathbb Z$ from that of the local minima in $S$. In particular we can study $F_\varepsilon^{\Z\cup S},$ whose potential does not depend on $\varepsilon$; the same works in regime \hyperlink{R2}{(R2)} if $s>1/2.$ The result in dimension one (\Cref{sub:lower_1d}) reduces to the classical two-well setting, for which the lower bound is known (see, e.g., \cite{ABS1,Palatucci,SavinValdinoci2010}). 
The genuinely new case is regime \hyperlink{R2}{(R2)} for $s=1/2$ (see \Cref{lemma:lb_1d_sub}) where an almost optimal transition between $k$ and $k+1$ rests at $s_k$ on an interval of length $\varepsilon/\delta_\varepsilon$: this lowers the cost by an amount depending on $\ell.$
 Then the $n$-dimensional case is obtained in Section \ref{sub:lower_nd} by slicing in all directions splitting the potential among them with a probability density on the sphere; optimizing this density we recover the dimensional constant $\lambda_{s,n}$. 

In \Cref{sec:upper_bound} we establish the upper bound. By density, it is enough to construct a recovery sequence for each function with a polyhedral jump set that jumps between consecutive wells in the limit. 
The recovery sequence depends on the regime and on $s$. In particular, for $s=1/2$, in regime \hyperlink{R1}{(R1)}, it is a transition between consecutive wells that occurs at scale $\varepsilon$;
in regime \hyperlink{R2}{(R2)}, the sequence additionally rests at the intermediate well on a layer of width of order $\varepsilon/\delta_\varepsilon$.
 For $s\in (1/2,1)$, each transition is described by an optimal profile problem defined in \Cref{Section:optimal}. We first define the sequence in one dimension and then extend it to higher dimensions by considering the one-dimensional profile along the normal direction to the jump set.

Finally, in \Cref{s<1/2}, we consider the case $s<1/2$ with the scaling $\theta_\varepsilon:=\varepsilon^{2s}$: the $\Gamma$-limit is the non-local energy $[u]^2_{H^s(Q)}$, plus the bulk term in regime \hyperlink{R1}{(R1)}, for functions which take values in $\union$ (in $\Z$ for \hyperlink{R2}{(R2)}).

\section{Compactness}\label{sec:compact}
In this section we prove \Cref{thm:compactness}. The proof has three steps.

\smallskip
\textit{1) Uniform bounds for the periodic energy.} We show that a sequence with
$\sup_\varepsilon F^\Z_\varepsilon(u_\varepsilon)<+\infty$ is bounded in $L^q(Q)$, for some $q>1$, after an
integer translation. 
\begin{itemize}
\item \textbf{$n\ge2$, $s\in[\frac12,1)$} (\Cref{lem:uniform_bound_nD}): bound in $L^{\frac n{n-1}}(Q)$. We
extend the argument of \cite{Garroni-Muller} ($n=2$, $s=\frac12$). Its key ingredient is a lower bound
for the non-local interaction between disjoint sets (\Cref{lem:kernel,lem:stima_bassa}), which gives an
isoperimetric-type inequality for the super-level sets of $u_\varepsilon$.
\item \textbf{$n=1$, $s=\frac12$} (\Cref{lem:uniform_bound_1D}): bound in $L^q(Q)$ for every $q<\infty$,
proved in \cite[Proposition~2.7]{Kurzke}.
\item \textbf{$n=1$, $s>\frac12$} (\Cref{lem:uniform_bound_1D}): bound in $L^\infty(Q)$. Every transition
between consecutive integers costs at least a fixed amount, so the number of transitions is bounded.
\end{itemize}

\textit{2) Compactness for the periodic energy} (Step~1 in \Cref{sottosezione_compattezza}). By
\eqref{eq:truncation}, the truncations of the translated sequence between consecutive integers have
bounded energy. Each of them is a two-well problem, hence compact in $L^1$ by \cite{ABS1} ($n=1$,
$s=\frac12$), \cite{Palatucci} ($n=1$, $s>\frac12$) and \cite{SavinValdinoci2010} ($n\ge2$). The $L^q$
bound of \textit{1)} controls the tails, and the limit belongs to $BV(Q;\Z)$.

\textit{3) General potential} (Steps~2--5 in \Cref{sottosezione_compattezza}). The piecewise affine homeomorphism $\phi^S$, defined in \textit{Step 2}, maps $\Z\cup S$
onto $\frac12\Z$ and satisfies $F_\varepsilon(u)\ge d^2F^\Z_\varepsilon(2\phi^S\circ u)$. Applying compactness to
$2\phi^S\circ u_\varepsilon$ gives compactness of $u_\varepsilon$ up to integer translations. If $u_\varepsilon$
converges, the translations are eventually constant up to a subsequence, as shown in \textit{Step 4}, and the limit belongs to $BV(Q;\Z\cup S)$. In
regime \hyperlink{R2}{(R2)}, the lower estimate in \eqref{eq:two_sided_estimate} forces $|\{u\in S\}|=0$.

\subsection{\texorpdfstring{Uniform  $L^q$ bounds in dimension $1$}{Uniform Lq bounds in dimension 1}}
The uniform bounds in the case $s=
\tfrac12$ were established by Kurzke in 
\cite[Proposition~2.7]{Kurzke}. For the reader’s convenience we recall this result 
in the following proposition. For
$s>\tfrac12$ we give a direct proof which yields an $L^\infty$ 
estimate.

\begin{prop}[Uniform $L^q$-estimates with $n=1$]\label{lem:uniform_bound_1D}
Let $n=1$ and let $I \subset \mathbb{R}$ be a bounded and open interval. Suppose that $\{u_\varepsilon\} \subset L^1(I)$ is a sequence with uniformly bounded energy, i.e.,
\[
\sup_{\varepsilon } F_\varepsilon^\Z(u_\varepsilon;I) < +\infty.
\]
Then, there exists a sequence of integers $\{a_\varepsilon\} \subset \mathbb{Z}$ such that the translated sequence $\{u_\varepsilon - a_\varepsilon\}$ is uniformly bounded in $L^q(I)$ for every $1 \le q < \infty$. 
\\\noindent
Furthermore, if $s \in \left(\tfrac12, 1\right)$, the translated sequence is also uniformly bounded in $L^\infty(I)$.
\end{prop}
\begin{proof}
We only prove the result for $s>\tfrac12$, since the case $s=\tfrac12$ is already proved in \cite{Kurzke}.
For each $\varepsilon > 0$, we define a translation $a_\varepsilon \in \mathbb{Z}$ as:
\[
a_\varepsilon := \max \left\{ a \in \mathbb{Z} : |\{ u_\varepsilon > a-1 \}| \geq \frac{|I|}{4} \right\}.
\]
This choice ensures that 
\(
|\{ u_\varepsilon \le a_\varepsilon \}| \geq \frac{|I|}{4}  \text{ and }  |\{ u_\varepsilon > a_\varepsilon - 1 \}| \geq \frac{|I|}{4}.
\)
Consider the truncated functions
\[
v_\varepsilon := (u_\varepsilon - a_\varepsilon)_+ \quad \text{and} \quad w_\varepsilon := (a_\varepsilon - 1 - u_\varepsilon)_+;
\]
then $|u_\varepsilon - a_\varepsilon| \leq 1 + v_\varepsilon + w_\varepsilon$. 
Since the energy is invariant under integer translations and non-increasing under this truncation operation, we have $F^\Z_\varepsilon(v_\varepsilon;I) \leq F^\Z_\varepsilon(u_\varepsilon;I)$ and $F^\Z_\varepsilon(w_\varepsilon;I) \leq F^\Z_\varepsilon(u_\varepsilon;I)$. Thus, $\{v_\varepsilon\}$ and $\{w_\varepsilon\}$ are sequences with uniformly bounded energy.
Furthermore, we may reduce the problem to the case of monotone functions via rearrangement. Let $v_\varepsilon^*$ denote the monotone non-decreasing rearrangement of $v_\varepsilon$. The potential term $\int_I \operatorname{dist}^2(v_\varepsilon,\Z)$ is invariant under rearrangement due to equimeasurability, while the fractional Gagliardo seminorm decreases, $[v_\varepsilon^*]_{H^s} \leq [v_\varepsilon]_{H^s}$ (see \cite{Garsia}). Moreover $\|v^*_\varepsilon\|_\infty=\|v_\varepsilon\|_\infty$ and $|\{v^*_\varepsilon=0\}|=|\{u_\varepsilon\le a_\varepsilon\}|\ge|I|/4$. An analogous argument applies to $w_\varepsilon$, with $|\{w^*_\varepsilon=0\}|\ge|\{u_\varepsilon>a_\varepsilon-1\}|\ge|I|/4$.
Therefore, it is sufficient to prove that $\|u_\varepsilon\|_\infty\le C$ for a sequence $\{u_\varepsilon\}$ with $F^\Z_\varepsilon(u_\varepsilon;I)\le M$ that satisfies
\begin{equation}\label{ipotesi}
u_\varepsilon \geq 0, \quad u_\varepsilon \text{ is non-decreasing}, \quad \text{and} \quad |\{u_\varepsilon = 0\}| \geq \frac{|I|}{4}.
\end{equation}
By the Sobolev embedding, \(u_\varepsilon \in C^{0,\alpha}(\bar I)\) with \(\alpha = s-\tfrac12\) (see \cite[Theorem 8.2]{Frazionari}). For each $\varepsilon$ we can define 
$$N_\varepsilon:= \left\lfloor\sup_I u_\varepsilon\right\rfloor.$$
We want to prove that if $
F^\Z_\varepsilon(u_\varepsilon)\leq M,
$
then $N_\varepsilon\leq C$. To this end, we show that each transition of $u_\varepsilon$ between two consecutive integers carries a strictly positive cost; hence, if the total energy is finite, such transitions can occur only a finite number of times.
\\\noindent
Let $u \in L^\infty(I) $ be a continuous function that satisfies \eqref{ipotesi}, and let $N:=\left\lfloor\sup_I u\right\rfloor$.
\\\noindent
For each integer \(0\leq j\leq N-1\), define
\(
I^j= [x^j, y^j]
\)
to be an interval such that
\[
u(x^j) = j + \frac{1}{4},
\qquad
u(y^j) = (j+1) - \frac{1}{4}.
\]
By monotonicity, we have
\[
j + \frac{1}{4} \le u(x) \le (j+1) - \frac{1}{4}, \quad \forall \,x \in I^j,
\]
and the intervals $I^j$ are pairwise disjoint.
Finally, set
$h^j = y^j - x^j,$
the length of the interval \(I^j\). We have that
\[
[u]_{H^s(I)}^2\geq \sum_{j=0}^{N-1}\int_{I^j}\int_{I^j} 
\frac{|u(x)-u(y)|^2}{|x-y|^{1+2s}} \, dx\,dy .
\]
We estimate it from below considering the rescaled function
\[
\hat u_j(\hat{x}) = u(x^j + h^j \hat{x})\quad \text{with }\hat x\in[0,1].
\]
Observe that
\[
|\hat u_j(1)-\hat u_j(0)| = \frac{1}{2}.
\]
By a change of variables
we obtain
\[
\int_{I^j}\int_{I^j} 
\frac{|u(x)-u(y)|^2}{|x-y|^{1+2s}} \, dx\,dy
=
(h^j)^{1-2s}
[\hat u_j]_{H^s(0,1)}^2.
\]
Now, the function \(\hat u_j\) spans a vertical distance 
\(1/2\) on \([0,1]\); therefore its fractional seminorm is bounded from below by a positive constant \(C\) depending only on \(s\), namely
\[
[\hat u_j]_{H^s(0,1)}^2
\ge C .
\]
Indeed, by the Sobolev embedding (\cite[Theorem 8.2]{Frazionari}) applied to $\hat u_j-m_j,$ with $m_j=\int_0^1 \hat u_j$, which has the same increments as $\hat u_j$, there exists a constant $C_H=C_H(s)>0 $ such that 
$$
|\hat u_j(x)-\hat u_j(y)|\leq C_H \|\hat u_j-m_j\|_{H^s(0,1)}|x-y|^{s-1/2}
$$
By the Poincaré inequality for fractional Sobolev spaces (see \cite[Theorem 1.17]{Leoni}), we have
$$
\|\hat u_j-m_j\|_{L^2(0,1)}\leq C_P [\hat u_j]_{H^s(0,1)},
$$
hence, taking $x=1, y=0,$
\[
[\hat u_j]_{H^s(0,1)}^2
\ge
\frac{1}{C_H^2(1+C_P^2)} (1/2)^2.
\]
We conclude that
\[
[u]_{H^s(I^j)}^2
\ge
c^*  \, (h^j)^{1-2s},\quad \text{with }c^*=\frac{1}{4C_H^2(1+C_P^2)} .
\]
Moreover, since $\operatorname{dist}(u,\Z)\ge \frac14$ on $I^j$, for the potential term we have
\[
\frac{1}{\varepsilon} \int_{I} \operatorname{dist}^2(u(x),\Z) \, dx
\geq
\frac{1}{16\,\varepsilon}  \sum_{j=0}^{N-1} h^j.
\]
In conclusion, we obtain
$$
F_\varepsilon^\Z(u)\geq \varepsilon^{2s-1} c^*
\sum_{j=0}^{N-1} (h^j)^{1-2s}
+
\frac{1}{16\,\varepsilon} \sum_{j=0}^{N-1} h^j\ge C(s)
\sum_{j=0}^{N-1}\left(\varepsilon^{2s-1}  (h^j)^{1-2s}+\frac{1}{\varepsilon}h^j\right).
$$
By Young's inequality in the form $x^\frac{1}{2s}y^\frac{2s-1}{2s} \le \frac{1}{2s}x + \frac{2s-1}{2s}y\le x+y$, we estimate
\[
F_\varepsilon^{\mathbb{Z}}(u) \ge C(s) \sum_{j=0}^{N-1} 
\left(\varepsilon^{2s-1}(h^j)^{1-2s}\right)^{\frac{1}{2s}}
\left(\frac{h^j}{\varepsilon}\right)^{\frac{2s-1}{2s}}
= C(s) N.
\]
Therefore, if $F_\varepsilon^{\mathbb{Z}}(u_\varepsilon) \le M$, then $N_\varepsilon \le C(M,s)$ and consequently $\|u_\varepsilon\|_\infty \le C(M,s)$.
\end{proof}

\subsection{Uniform \texorpdfstring{$L^q$}{Lq} bounds in dimension \texorpdfstring{$n\ge 2$}{n ≥ 2}}
\label{sec:Lq_bounds_n_ge_2}

The proof is based on a summability estimate for the super-level sets of the sequence.
The heuristic idea is that the energy $F^\Z_\varepsilon(u_\varepsilon)$ behaves asymptotically like the $BV$-seminorm. By comparing consecutive super-level sets, one expects that the perimeter of the set $\{u_\varepsilon > k\}$ is controlled by the energy of the truncation $T_{k-1}u_\varepsilon := (u_\varepsilon \vee (k-1) )\wedge k$. In combination with the isoperimetric inequality, this would formally lead to an estimate of the form 
\begin{equation}\label{isoperimetrica}
| \{u_\varepsilon > k\} |^{\frac{n-1}{n}} \le C F^\Z_\varepsilon(T_{k-1}u_\varepsilon),
\end{equation}
from which the $L^{\frac{n}{n-1}}$ bound would follow.
The key tool is a lower bound for the non-local interaction between disjoint sets (\Cref{lem:kernel,lem:stima_bassa}), in the spirit of \cite{Savin}.

In what follows $B_r(x)$ will denote the ball centered at $x$ with radius $r$ (we will omit the center when $x=0$ and simply write $B_r$). Similarly $Q_r(x)$ will denote the cube of side $r$ centered at $x$. We denote simply by $Q_r$ the cube of side $r$ centered at $0$.

\begin{prop}[Proposition 14, 
\cite{Garroni-Muller}]\label{nucleo_regolare}
Let $n\ge 2,$ and let 
$\varphi \in C_c^\infty(\R^n)$ satisfy $\varphi\ge 0,$ $\varphi>0$ in $B_{1/4}$ (and $\varphi$ even),
$\int \varphi = 1$, $\operatorname{supp}\varphi \subseteq B_{1/2}$ and set $\varphi_\delta(x):=\delta^{-n}\varphi(x/\delta)$. 

For every $\mu\in \left(0,\frac12\right]$ there exists $c(\mu) > 0$ (depending also on $n$ and $\varphi$) such that for every measurable set 
$A \subset Q_1$ with $\mu\le |A| \le 1-\mu$, 
\begin{equation}\label{regular_kernel}
\frac{1}{\delta}
\int_A \int_{Q_1\setminus A}
\varphi_\delta(x-y)\,dy\,dx
\ge c(\mu),
\qquad \forall \delta \in (0,1].
\end{equation}
\end{prop}
The proof in \cite{Garroni-Muller} is written for $n=2$ but holds verbatim in any dimension.

\begin{lemma}\label{lem:kernel}
Let $n\ge 2,$ and let $A, B \subseteq Q_1$ be two disjoint subsets. Define $E := Q_1 \setminus (A \cup B)$. For every $\mu \in \left(0, \frac12\right]$, let $c(\mu)$ be the constant of \Cref{nucleo_regolare}. There exists $C(\mu)>0$ (depending also on $\varphi$, the dimension $n$ and $s$) such that, if $\mu \leq |A| \leq 1 - \mu$ and $|E|<c(\mu)/2$, then 
\begin{equation}\label{stima_bassa_kernel}
\int_A \int_B \frac{1}{|x - y|^{n+2s}} \, dy \, dx \ge
\begin{cases} 
\displaystyle
C(\mu) \log \frac{c(\mu)}{2|E|}, & \text{if } s = \frac12, \\ \\
\displaystyle C(\mu) \left(\frac{c(\mu)}{|E|}\right)^{2s-1}, & \text{if } s \in \left(\frac12, 1\right). 
\end{cases}
\end{equation}
\end{lemma}

\begin{proof} 
If $|E|=0 $ the right-hand side is $+\infty$, and the argument below with $m=\infty$
 shows that so is the left-hand side. Assume $|E|>0.$
To obtain the estimate \eqref{stima_bassa_kernel}, we approximate the singular kernel with a dyadic sum of regular kernels.
Let $\varphi_\delta$ be as in \Cref{nucleo_regolare} and assume 
$\frac{2|E|}{c(\mu)} \le \delta \le 1$. Then 
\begin{equation}\label{44}
\frac{1}{\delta} \int_A \int_B \varphi_\delta(x-y)\,dy\,dx
=
\frac{1}{\delta} \int_A \int_{Q_1\setminus A} \varphi_\delta(x-y)\,dy\,dx
-
\frac{1}{\delta} \int_A \int_E \varphi_\delta(x-y)\,dy\,dx
\ge
c(\mu) - \frac{|E|}{\delta}
\ge
\frac{c(\mu)}{2}.
\end{equation}
 Let $m \in \mathbb{N}$ be defined by
\begin{equation}\label{4.1}
m:=\left\lfloor\log_2\frac{c(\mu)}{2|E|}\right\rfloor,\quad \text{that is}\quad 
2^{-(m+1)} \leq \frac{2|E|}{c(\mu)} \leq 2^{-m}, 
\end{equation}
which is well defined since $|E|< c(\mu)/2.$
Choose $\delta_i = 2^{-i}$ for $i = 0, \dots, m$, so that $\frac{2|E|}{c(\mu)} \leq\delta_m\le \delta_i\le \delta_0=1.$
We claim that there exists a constant $C$ (depending on $n$, $s$ and $\varphi$) such that 
\begin{equation}\label{kernel_approximation}
\frac{1}{|x - y|^{n+2s}} \ge C \sum_{i=0}^m \frac{1}{\delta_i^{2s}} \varphi_{\delta_i}(x - y),\quad \forall\,x,y\in Q_1, x\neq y.
\end{equation}
Recall that $\varphi_\delta(x - y) = 0$ whenever $|x - y| \ge \delta/2$; hence, if $|x-y|\ge\delta_0=1$ then every term of the sum vanishes.\\\noindent
If $|x - y| < \delta_m$, then
\[\displaystyle
\sum_{i=0}^m \frac{1}{\delta_i^{2s}} \varphi_{\delta_i}(x - y) \le \sup \varphi \sum_{i=0}^m \frac{1}{\delta_i^{n+2s}} 
\leq \sup \varphi \, \frac{1}{\delta_m^{n+2s}} \sum_{i=0}^\infty \left(\frac{1}{2^{n+2s}}\right)^i 
\le C(n) \, \sup \varphi \, \frac{1}{\delta_m^{n+2s}} 
\le \frac{C(n,\varphi)}{|x - y|^{n+2s}}.
\]
On the other hand, if $\delta_{i_0+1} \le |x - y| < \delta_{i_0}$ for some $i_0 \in \{0, \dots, m-1\}$, then $\varphi_{\delta_i}(x-y)=0$ for all $i\ge i_{0}+2,$ hence, since $|x-y|<\delta_{i_0}=2\delta_{i_0+1},$
\[\displaystyle
\sum_{i=0}^m \frac{1}{\delta_i^{2s}} \varphi_{\delta_i}(x - y) \le \sup \varphi \sum_{i=0}^{i_0+1} \frac{1}{\delta_i^{n+2s}} 
\le  \sum_{i=0}^\infty \left(\frac{1}{2^{n+2s}}\right)^i \frac{ \sup \varphi } {\delta_{i_0+1}^{n+2s}} 
\le \frac{C(n,\varphi)}{|x - y|^{n+2s}}.
\]
This completes the proof of \eqref{kernel_approximation}. Integrating \eqref{kernel_approximation} over $A \times B$ and using \eqref{44}, we obtain
\begin{align*}
\int_A \int_B \frac{1}{|x - y|^{n+2s}} \, dy \, dx 
    &\ge C \sum_{i=0}^m \frac{1}{\delta_i^{2s-1}} \int_A \int_B \frac{1}{\delta_i} \varphi_{\delta_i}(x - y) \, dy \, dx \ge C \frac{c(\mu)}{2} \sum_{i=0}^m (2^{2s-1})^i.
\end{align*}
If $s=\frac12$ the sum equals $m+1\ge \log_2\frac{c(\mu)}{2|E|}$ which gives the first line in \eqref{stima_bassa_kernel}. If $s\in \left(\frac12,1\right)$ then by \eqref{4.1},
\[
\sum_{i=0}^m(2^{2s-1})^i\ge (2^m)^{2s-1}\ge \left(\frac{c(\mu)}{4|E|}\right)^{2s-1},
\]
which gives the second line of \eqref{stima_bassa_kernel}, up to renaming $C(\mu).$
\end{proof}

\begin{ohss}[Necessity of the assumption $|E|<c(\mu)/2$]\label{rem:necessita1} For $s\in\left(\frac12,1\right)$, without the smallness assumption on $|E|$, the statement is false. Indeed, let $A$ be the left half of $Q_1$, let $B\subset Q_1$ be a ball of measure $\eta$ contained in the right half and at distance at least $1/4$ from $A$, and set $E=Q_1\setminus(A\cup B)$. Then $|E|=1/2-\eta$, while \[ \int_A\int_B\frac{dy\,dx}{|x-y|^{n+2s}} \le |A||B|4^{n+2s}\xrightarrow{\eta\to0}0. \] On the other hand, the right-hand side of \eqref{stima_bassa_kernel} converges to the positive constant $C(\mu)(2c(\mu))^{2s-1}$. Thus, some smallness assumption on $|E|$ (equivalently, a lower bound on $|B|$) is necessary. 
Instead, if $s=\frac12$ the assumption is unnecessary: if $|E|\ge c(\mu)/2$, the right-hand side is non-positive and the inequality is trivial. See also the corresponding Lemma 15 in \cite{Garroni-Muller}. \end{ohss}

\begin{lemma}\label{lem:stima_bassa} Let $n \ge 2$ and $s\in \left[\frac12,1\right).$
Let $Q$ be an open cube in $\R^n$. For every $\mu\in (0,1)$ there exist $\gamma_0>0$ and two constants $C_1,C_2>0$ (depending on $\mu$, $n$ and $s$) such that for every $A,B \subseteq Q$ with $A \cap B = \emptyset$, $|A| \le (1-\mu)|Q|$ and 
\begin{equation}\label{ipotesi_su_E}
   |E|\le \gamma_0|A|,\quad \text{where } E:=Q\setminus (A \cup B),
\end{equation}
the following inequality holds
\begin{equation}\label{stima_bassa}
\int_A \int_B \frac{1}{|x - y|^{n+2s}} \, dy \, dx \ge
\begin{cases} 
C_1|A|^{\frac{n-1}{n}} \log(C_2|A|/|E|),\quad  & \text{if } s = \frac12, \\ 
C_1 |A|^{\frac{n-2s}{n}} (|A|/|E|)^{2s-1} & \text{if } s \in \left(\frac12, 1\right). 
\end{cases}
\end{equation}
\end{lemma}

\begin{proof} Assume that $|E|,|A|>0,$ otherwise it is trivial.
Without loss of generality we may assume $Q=(-1/2,1/2)^n$, since both sides of 
\eqref{stima_bassa} scale as $r^{n-2s}$ under $x\mapsto r x$ (the ratio 
$|A|/|E|$, the volume constraint $|A|\le(1-\mu)|Q|$ and hypothesis \eqref{ipotesi_su_E} being scale 
invariant); the general case follows by rescaling.\\\noindent
Let $Q'=Q/3$ denote the concentric subcube. We first prove the interior estimate
\begin{equation}\label{stima_interna}
\int_A \int_B \frac{1}{|x - y|^{n+2s}} \, dy \, dx \ge
\begin{cases} 
C_1|A\cap Q'|^{\frac{n-1}{n}} \log(C_2|A\cap Q'|/|E|),\quad & \text{if } s = \frac12, \\ 
C_1 |A\cap Q'|^{\frac{n-2s}{n}} (|A\cap Q'|/|E|)^{2s-1}, \quad& \text{if } s \in \left(\frac12, 1\right),
\end{cases}
\end{equation}
under the hypothesis $|E|\le \gamma_1|A\cap Q'|$ and $|A|\le 3^{-n-1}.$
Then we derive the global estimate \eqref{stima_bassa} by a reflection and extension argument.

Set $\mu':=\min\{\mu,3^{-n-1}\}$ and suppose first that $|A|>\mu'$. Then
$\mu'\le|A|\le1-\mu\le1-\mu'$ and, by \eqref{ipotesi_su_E},
$|E|\le\gamma_0|A|\le\gamma_0<c(\mu')/2$, provided
\begin{equation}\label{eps0_a}
\gamma_0<\frac{c(\mu')}{2}.
\end{equation}
Hence \Cref{lem:kernel} applies and gives \eqref{stima_bassa} at once, since $|A|\le 1,$ with $C_1,C_2$
depending only on $\mu$ ($n$ and $s$). 

We may thus assume $|A| \le 3^{-n-1}$ and we prove \eqref{stima_interna} by a covering argument.
\\\noindent
By scaling the estimate \eqref{stima_bassa_kernel}, with $\mu =\frac14$, we obtain the following: let $Q_r \subset Q$ be a cube of side $r>0$, and let $A \subseteq Q$ be such that 
\[
\frac14 \le \frac{|A \cap Q_r|}{|Q_r|} \le \frac34 \quad \text{and}\quad |E\cap Q_r|< \frac{c(1/4)}{2}r^n;
\]
then
\begin{equation}\label{stima_cubo}
\int_{A \cap Q_r} \int_{B \cap Q_r} \frac{1}{|x - y|^{n+2s}} \, dy \, dx \ge 
\begin{cases} 
\displaystyle C \, r^{n-1} \log \frac{c \, r^n}{2 |E \cap Q_r|}, & \text{if } s = \frac12, \\[0.8em] 
 \displaystyle C \, r^{n-2s} \left(\frac{c\,r^n}{2|E \cap Q_r|}\right)^{2s-1}, & \text{if } s \in \left(\frac12, 1\right),
\end{cases}
\end{equation}
where $c=c(1/4)$ and $C>0$ is a constant depending only on $s$ and the dimension $n$.

We claim that for almost every $x \in A \cap Q'$ there exists a cube $Q_{r_x}(x)$, centered at $x$ with side length $r_x\le 1/3$ (hence $Q_{r_x}(x)\subseteq Q$), such that
\begin{equation}\label{densita_proof}
\frac{1}{4} \leq \frac{|A \cap Q_{r_x}(x)|}{|Q_{r_x}(x)|} \le \frac{3}{4}.
\end{equation}
Indeed, for almost every $x \in A \cap Q'$, the density of $A$ at $x$ is $1$, while for the cube of side $1/3$ centered at $x$, i.e., $Q_{1/3}(x) \subset Q$, we have
\[
\frac{|A \cap Q_{1/3}(x)|}{|Q_{1/3}(x)|} \le 3^n |A| \le \frac{1}{3}.
\]
By the continuity of $\frac{|A \cap Q_r(x)|}{|Q_r(x)|}$ with respect to $r$, this guarantees the existence of an intermediate scale $r_x$ satisfying \eqref{densita_proof}.

By the Besicovitch covering theorem we obtain a family of disjoint cubes $Q_{r_i}(x_i)$, $i \in I$, satisfying \eqref{densita_proof} such that
\[ \sum_{i \in I} r_i^n = \sum_{i \in I} |Q_{r_i}(x_i)| \geq \tilde{c} |A \cap Q'|, \]
where $\tilde{c}$ depends only on $n$. Consider a subfamily $J$ corresponding to cubes where $E$ has sufficiently low density,
\begin{equation}\label{def_J}
J := \left\{ i \in I : \frac{|E \cap Q_{r_i}(x_i)|}{r_i^n} < \frac{2}{\tilde{c}} \frac{|E|}{|A \cap Q'|} \right\}.
\end{equation}
We have
\begin{equation}\label{stima_quadrati}
\sum_{i \in J} r_i^n = \sum_{i \in I} r_i^n - \sum_{i \notin J} r_i^n \geq \tilde{c} |A \cap Q'| - \frac{\tilde{c}}{2} \sum_{i \notin J} \frac{|E \cap Q_{r_i}(x_i)|}{|E|} |A \cap Q'| \geq \frac{\tilde{c}}{2} |A \cap Q'|.
\end{equation}
We now check that \eqref{stima_cubo} is applicable to each cube $Q_{r_i}(x_i)$ with $i \in J$. The density condition is satisfied thanks to  \eqref{densita_proof}; the estimate on $E$ is satisfied by \eqref{def_J} and the hypothesis $|E|\le \gamma_1|A\cap Q'|$ of \eqref{stima_interna}; indeed
$$
\frac{|E\cap Q_{r_i}|}{r_i^n}<\frac{2}{\tilde{c}}\frac{|E|}{|A \cap Q'|}\le \frac{2\gamma_1}{\tilde c}\le \frac{c}{2}\quad \text{provided}\quad \gamma_1\le \frac{\tilde c\,c}{4}.
$$
The same choice guarantees $\frac{c\tilde c|A\cap Q'|}{4|E|}\ge 1,$ so that the logarithm below is non-negative.
\\\noindent
Then we can fix $\gamma_1:=3^{-n}\min\{c(\mu')/4,\ \tilde c\,c/4\}$.

By \eqref{stima_cubo} and \eqref{stima_quadrati} we get
\begin{align*}
\iint_{A \times B} \frac{dx \, dy}{|x-y|^{n+2s}} 
&\geq \sum_{i \in J} \iint_{(A \cap Q_{r_i}) \times (B \cap Q_{r_i})} \frac{dx \, dy}{|x-y|^{n+2s}} \\
&\geq
\begin{cases} 
\sum_{i \in J} C r_i^{n-1} \log \left( \frac{c r_i^n}{2|E \cap Q_{r_i}|} \right), & \text{if } s = \frac12 \\[10pt] 
\sum_{i \in J} C r_i^{n-2s} \left( \frac{c r_i^n}{2|E \cap Q_{r_i}|} \right)^{2s-1}, & \text{if } s \in \left(\frac12, 1\right) 
\end{cases} \\
&{\geq} 
\begin{cases} 
\left( \sum_{i \in J} r_i^{n-1} \right) C \log \left( \frac{c \tilde{c} |A \cap Q'|}{4|E|} \right), \\[10pt] 
\left( \sum_{i \in J} r_i^{n-2s} \right) C \left( \frac{c \tilde{c} |A \cap Q'|}{4|E|} \right)^{2s-1} 
\end{cases} \\
&{\geq} 
\begin{cases} 
C |A \cap Q'|^{\frac{n-1}{n}} \log \left( \frac{c \tilde{c} |A \cap Q'|}{4|E|} \right), & \text{if } s = \frac12 ,\\[10pt] 
C |A \cap Q'|^{\frac{n-2s}{n}} \left( \frac{c \tilde{c} |A \cap Q'|}{4|E|} \right)^{2s-1}, & \text{if } s \in \left(\frac12, 1\right), 
\end{cases}
\end{align*}
where we used that
$$
\sum_{i\in J}r_i^k\geq \left(\sum_{i\in J}r_i^n\right)^\frac{k}{n}\quad \text{if }r_i\geq 0, 0<k\le n.
$$
 Now, to establish the full estimate \eqref{stima_bassa}, we use a reflection argument as in \cite{Garroni-Muller}. 
We extend $A, B$, and $E$ to the larger cube $\tilde{Q} = 3Q$ by reflecting them across the hyperplanes $x_k = \pm 1/2$ for $k=1, \dots, n$. Specifically, let $\mathcal{I} = \{-1, 0, 1\}^n$ be the set of multi-indices. For each $\alpha \in \mathcal{I}$, we define the reflection $\pi_\alpha(x) = (\pi_{\alpha_1}(x_1), \dots, \pi_{\alpha_n}(x_n))$, where $\pi_0(t)=t$, $\pi_{-1}(t)=-1-t$, and $\pi_1(t)=1-t$. We then set $\tilde{A} = \bigcup_{\alpha \in \mathcal{I}} \pi_\alpha(A)$, and similarly for $\tilde{B}$ and $\tilde{E}$. 
\\\noindent
Since $|\pi_\alpha(x) - \pi_\beta(y)| \geq |x - y|$ for all $\alpha, \beta \in \mathcal{I}$ and $x,y\in Q$, the integral over the extended sets satisfies:
\begin{equation}
\int_{\tilde{A}} \int_{\tilde{B}} \frac{1}{|x-y|^{n+2s}} \, dx \, dy \leq 3^{2n} \int_A \int_B \frac{1}{|x-y|^{n+2s}} \, dx \, dy. 
\end{equation}
Moreover, denoting by $\tilde Q'=\tilde Q/3=Q$ the subcube concentric to $\tilde Q$,
$$
|\tilde Q'\cap \tilde A|=|A|,\quad |\tilde E|=3^n|E|,\quad |\tilde A|=3^n|A|\leq 3^n(1-\mu)=(1-\mu)|\tilde Q|,
$$
 so the volume constraints are preserved ($|\tilde A|/|\tilde Q|=|A|\le 3^{-n-1}$), and the hypothesis $|E|\le \gamma_1|A\cap Q'|$, required in \eqref{stima_interna}, corresponds to $3^n|E|\le\gamma_1|A| $, which is satisfied thanks to \eqref{ipotesi_su_E} since $\gamma_0=3^{-n}\gamma_1.$
Applying the interior estimate \eqref{stima_interna} to $\tilde{A}, \tilde{B} \subset \tilde{Q}$ and rescaling by a factor of 3, we obtain

\begin{equation}
\int_{\tilde A}\int_{\tilde B}\frac{dx\,dy}{|x-y|^{n+2s}}\ge
\begin{cases}
\tilde C_1\,|\tilde A\cap Q|^{\frac{n-1}{n}}\log\dfrac{\tilde C_2|\tilde A\cap Q|}{|\tilde E|} & s=\frac12,\\[8pt]
\tilde C_1\,|\tilde A\cap Q|^{\frac{n-2s}{n}}\Big(\dfrac{|\tilde A\cap Q|}{|\tilde E|}\Big)^{2s-1} & s\in\left(\frac12,1\right),
\end{cases}
\end{equation}
which, together with the previous inequality, gives \eqref{stima_bassa} with constants $C_1, C_2$ depending only on $n,s$ and $\mu$.
\end{proof}
\begin{ohss}
    \label{rem:necessita2} By \Cref{rem:necessita1}, in \Cref{lem:stima_bassa}, \eqref{ipotesi_su_E} is not necessary for $s=\frac12:$
   it suffices to choose $C_2\le\gamma_0$.
\end{ohss}
We now establish a uniform $L^{\frac{n}{n-1}}$ bound for the sequence, up to integer translations. 
 The proof in the case $s=\frac12$ follows the same lines as the two-dimensional argument presented in \cite[Theorem 13]{Garroni-Muller}.

\begin{prop}\label{lem:uniform_bound_nD}
Let $n\ge 2$ and $s \in \left[\frac12,1\right)$, and let $Q\subset \R^n$ be an open cube. Given $M > 0$, there exist $\varepsilon_0 \in (0,1)$ and a constant $C(M,n,Q,s)$ such that the following holds. If $0 < \varepsilon < \varepsilon_0$ and
\begin{equation}\label{eq:energy_bound}
F^\Z_\varepsilon(u_\varepsilon) \leq M,
\end{equation}
then there exists $a_\varepsilon \in\mathbb{Z}$ such that
\begin{equation}\label{eq:L2_bound}
\|u_\varepsilon-a_\varepsilon\|_{L^\frac{n}{n-1}(Q)} \, \leq C(M,n,Q,s).
\end{equation}
\end{prop}

\begin{proof}
Throughout the proof, $\varepsilon_0$ is chosen small enough, depending only on $M$, $n$, $s$ and $\sigma$ below.
    We define 
    $a_\varepsilon = \max \{ a \in \mathbb{Z} : |\{ u_\varepsilon - a > 0 \}| \geq \frac{1}{2} |Q| \}$. With this choice we have that the translated sequence $v_\varepsilon:=u_\varepsilon-a_\varepsilon$ satisfies 
    \begin{equation}\label{normalizzazione}
        |\{v_\varepsilon \geq 0\}| \geq \frac{1}{2}|Q| \quad \text{and}\quad  |\{v_\varepsilon \leq 1\}| \geq \frac{1}{2}|Q|. 
        \end{equation}
Since $F^\Z_\varepsilon$ is invariant under integer translations,
$F^\Z_\varepsilon(v_\varepsilon)=F^\Z_\varepsilon(u_\varepsilon)\le M$.
    Fix $\sigma \in (0, \frac{1}{4})$ and $k \in \mathbb{Z}$. Fix the following level sets
    \begin{equation*}
        A_k^\varepsilon = \{x \in Q : v_\varepsilon > k + 1 - \sigma\}, \quad  B_k^\varepsilon = \{x \in Q : v_\varepsilon < k + \sigma\} \quad\text{and}
    \end{equation*}
    \begin{equation*}
        E_k^\varepsilon := \{x \in Q : k + \sigma \le v_\varepsilon \le k + 1 - \sigma\}=Q\setminus (A_k^\varepsilon\cup B_k^\varepsilon).
    \end{equation*}
    \textit{Step 1.}  We prove that there exists $C = C(M, n, Q, \sigma,s)$ such that 
    \begin{equation}\label{Lq_unif_final}
        \sum_{k=0}^{\infty} |A_k^\varepsilon|^{\frac{n-1}{n}} \le C \text{.}
    \end{equation}
    We treat the positive part, i.e., the case $k\geq 1$; the estimate for the negative part follows by symmetry (applying \Cref{lem:stima_bassa} with the roles of $A$ and $B$ exchanged, which is
legitimate since the kernel is symmetric, although the right-hand side of \eqref{stima_bassa} is
not).

    Throughout, $T_kv_\varepsilon:=(v_\varepsilon\vee k)\wedge(k+1)$ denotes the truncation at levels
$k$ and $k+1$, so that $|T_kv_\varepsilon(x)-T_kv_\varepsilon(y)|\ge1-2\sigma$ for
$x\in A_k^\varepsilon$, $y\in B_k^\varepsilon$.
\\
Let us observe that for $k\ge1$ one has
$A_k^\varepsilon\subseteq\{v_\varepsilon>1\}$, so \eqref{normalizzazione} gives
$|A_k^\varepsilon|\le\frac12|Q|$: \Cref{lem:stima_bassa} is therefore applicable with
$\mu=1/2$. (For $k\le-1$ one argues symmetrically, using
$B_k^\varepsilon\subseteq\{v_\varepsilon<0\}$.)
    \begin{itemize}[leftmargin=0.5pt]
        \item \emph{Case $s = \frac12$.}
    By the energy bound, we have $|E_k^\varepsilon| \leq \frac{M}{\sigma^2} \varepsilon |\log \varepsilon|$. 
    Fix $\alpha \in (0, 1/3)$; then $|E_k^\varepsilon| \leq \varepsilon^{3\alpha}$ for $\varepsilon<\varepsilon_0.$

\begin{itemize}
    \item \emph{If  $|A_k^\varepsilon|^{\frac12}>\varepsilon^{\alpha}$,} then
 \Cref{lem:stima_bassa} applies with
$\mu=1/2$ (by \Cref{rem:necessita2}, we do not need to check \eqref{ipotesi_su_E}):
\[
\int_{A_k^\varepsilon}\int_{B_k^\varepsilon}\frac{dy\,dx}{|x-y|^{n+1}}
\ \ge\ C_1|A_k^\varepsilon|^{\frac{n-1}{n}}\log\frac{C_2|A_k^\varepsilon|}{\varepsilon^{3\alpha}} .
\]
Since $|A_k^\varepsilon|>\varepsilon^{2\alpha}$,
\begin{equation}\label{stima_log}
\log\frac{C_2|A^\varepsilon_k|}{\varepsilon^{3\alpha}}\ \ge\
\log\frac{C_2\varepsilon^{2\alpha}}{\varepsilon^{3\alpha}}
=\alpha|\log\varepsilon|+\log C_2\ \ge\ \frac{\alpha}{2}|\log\varepsilon|
\end{equation}
for $\varepsilon$ small, whence, we have
\[
|A_k^\varepsilon|^{\frac{n-1}{n}}\ \le\ C\,F^{\Z}_\varepsilon(T_kv_\varepsilon).
\]
    \item\emph{If  $|A_k^\varepsilon|^{1/2}\le\varepsilon^{\alpha}$.} 
 Let $\gamma \in (1, 2)$. We claim that there exists $C=C(n,Q,\alpha,\gamma,M,\varepsilon_0)$
 such that
\begin{equation}\label{51}
|A_k^\varepsilon|^\frac{n-1}{n} \leq \frac{C}{(k - \sigma)^\gamma}.
\end{equation}
Indeed by the Sobolev embedding $H^{1/2}(Q) \hookrightarrow L^{\frac{2n}{n-1}}(Q)$ (see \cite[Theorem 6.7]{Frazionari}) and the fact that $|\{v_\varepsilon \le 1\}| \geq \frac{1}{2}|Q|$ (that allows us to estimate the $L^2$ norm of $(v_\varepsilon-1)_+$ with the seminorm), we have
\[
\|(v_\varepsilon - 1)_+ \|_{L^\frac{2n}{n-1}(Q)} \leq C(n,Q) [(v_\varepsilon - 1)_+]_{H^{1/2}(Q)} \leq C (|\log \varepsilon| F^\Z_\varepsilon(v_\varepsilon))^{1/2}.
\]

By Hölder's inequality we get
\begin{align*}
|A_k^\varepsilon|(k - \sigma) &\leq \int_{A_k^\varepsilon} |(v_\varepsilon - 1)_+| \, dx \leq |A_k^\varepsilon|^\frac{n+1}{2n} \|(v_\varepsilon - 1)_+\|_{L^\frac{2n}{n-1}(Q)}.
\end{align*}
Hence
\[
|A_k^\varepsilon|^\frac{n-1}{2n} \leq \frac{C}{(k - \sigma)} (|\log \varepsilon| F^\Z_\varepsilon(v_\varepsilon))^{1/2}.
\]
Since  $|A_k^\varepsilon|^{1/2} \leq \varepsilon^\alpha$ and
$$
\beta:=\frac{n-1}{\gamma n}-\frac{n-1}{2n}>0\quad \text{for}\quad \gamma\in (1,2),
$$
we have that $\varepsilon^{2\alpha\beta} |\log \varepsilon|^{1/2}$ is bounded for $\varepsilon \leq \varepsilon_0$.
This implies that
\[
|A_k^\varepsilon|^{\frac{n-1}{\gamma n}} \leq C \frac{1}{(k - \sigma)} (F^\Z_\varepsilon(v_\varepsilon))^{1/2},
\]
and \eqref{51} follows by raising the last inequality to the power $\gamma$.
In conclusion, fixing $\alpha$ and $\gamma$ as before, we have \eqref{Lq_unif_final} by \eqref{eq:truncation}; indeed
\[
\sum_{k=1}^\infty |A_k^\varepsilon|^\frac{n-1}{n} \leq C \sum_{k=1}^\infty \left( \frac{1}{(k - \sigma)^\gamma} + F^\Z_\varepsilon(T_k v_\varepsilon) \right) \leq C(1 + F^\Z_\varepsilon(v_\varepsilon))\le C(1 + M).
\]
\end{itemize}

    \item\noindent \emph{Case $s > \frac12$.} By the energy bound, we have $|E_k^\varepsilon| \leq \varepsilon\frac{M}{\sigma^2}$. 

Fix $\tau\in(0,1)$ and $\gamma>1$ such that
\begin{equation}\label{scelta_theta_gamma}
\tau\,\beta_s\ \ge\ \frac{2s-1}{2},\qquad\text{where}\qquad
\beta_s:=\frac{n-1}{\gamma n}-\frac{n-2s}{2n}>0 .
\end{equation}
Such a choice is possible: as $\gamma\to1^+$ and $\tau\to1^-$ one has
$\beta_s\to\frac{n-2+2s}{2n}$, and
\[
\frac{n-2+2s}{2n}>\frac{2s-1}{2}
\iff (n-2+2s)-n(2s-1)>0
\iff 2(n-1)(1-s)>0,
\]
which holds for every $n\ge2$ and $s\in\left(\frac12,1\right)$.
 \begin{itemize}
     \item
\emph{If  $|A_k^\varepsilon|>\varepsilon^{\tau}$,} then
\begin{equation}
\frac{|E_k^\varepsilon|}{|A_k^\varepsilon|}\ \le\ \frac{M}{\sigma^2}\,
\varepsilon^{1-\tau}\ \le\ \gamma_0\qquad\text{for }\varepsilon\text{ small,}
\end{equation}
so \eqref{ipotesi_su_E} holds and \Cref{lem:stima_bassa} applies with $\mu=\frac12$. Therefore
\begin{align*}
F^{\Z}_\varepsilon(T_kv_\varepsilon)
&\ \ge(1-2\sigma)^2\,\varepsilon^{2s-1}\int_{A^\varepsilon_k}\int_{B^\varepsilon_k}
\frac{dy\,dx}{|x-y|^{n+2s}}
+\frac{1}{\varepsilon}\int_{E^\varepsilon_k}\operatorname{dist}^2(v_\varepsilon,\Z)\,dx\\
&\ \ge\ C\,\varepsilon^{2s-1}|A^\varepsilon_k|^{\frac{n-2s}{n}}
\Big(\frac{|A^\varepsilon_k|}{|E^\varepsilon_k|}\Big)^{2s-1}
+\frac{\sigma^2}{\varepsilon}|E^\varepsilon_k|\\
&\ =\ C\,\varepsilon^{2s-1}|A^\varepsilon_k|^{\frac{2s(n-1)}{n}}|E^\varepsilon_k|^{1-2s}
+\frac{\sigma^2}{\varepsilon}|E^\varepsilon_k| .
\end{align*}
Using Young's inequality with exponents $p=2s$ and $q=\frac{2s}{2s-1}$,
\[
F^{\Z}_\varepsilon(T_kv_\varepsilon)\ \ge\ C(C_1,\sigma,s)
\Big(\varepsilon^{2s-1}|A^\varepsilon_k|^{\frac{2s(n-1)}{n}}|E^\varepsilon_k|^{1-2s}\Big)^{\frac{1}{2s}}
\Big(\frac{|E^\varepsilon_k|}{\varepsilon}\Big)^{\frac{2s-1}{2s}}
=C(C_1,\sigma,s)\,|A_k^\varepsilon|^{\frac{n-1}{n}}.
\]
\item
\emph{If  $|A_k^\varepsilon|\le\varepsilon^{\tau}$,} we argue as in the case
$s=\frac12$. By the Sobolev embedding $H^{s}(Q)\hookrightarrow L^{\frac{2n}{n-2s}}(Q)$ (admissible since $n \ge 2>2s$), \eqref{normalizzazione}, the Poincaré inequality  for functions vanishing on a set with positive measure,  
$[(v_\varepsilon-1)_+]^2_{H^s}\le\varepsilon^{1-2s}F^\Z_\varepsilon(v_\varepsilon)$ and Hölder's
inequality, we have
\[
|A_k^\varepsilon|^{\frac{n-2s}{2n}}\ \le\ \frac{C}{k-\sigma}\,
\varepsilon^{-\frac{2s-1}{2}}\big(F^{\Z}_\varepsilon(v_\varepsilon)\big)^{1/2} .
\]
Multiplying by $|A_k^\varepsilon|^{\beta_s}\le\varepsilon^{\tau\beta_s}$ and using
\eqref{scelta_theta_gamma},
\[
|A_k^\varepsilon|^{\frac{n-1}{\gamma n}}\ \le\ \frac{C}{k-\sigma}\,
\varepsilon^{\tau\beta_s-\frac{2s-1}{2}}\big(F^{\Z}_\varepsilon(v_\varepsilon)\big)^{1/2}
\ \le\ \frac{C\sqrt M}{k-\sigma},
\]
and raising to the power $\gamma$ we obtain again \eqref{51}.
 \end{itemize}
Combining the two cases and using \eqref{eq:truncation}, we have
\begin{equation}\label{Lq_unif}
\sum_{k=1}^{\infty} |A^{\varepsilon}_{k}|^{\frac{n-1}{n}}
\leq C \sum_{k=1}^\infty \left( \frac{1}{(k - \sigma)^\gamma} + F^\Z_\varepsilon(T_k v_\varepsilon) \right) \leq C(1 + F^\Z_\varepsilon(v_\varepsilon))\le C(1 + M).
\end{equation}
\end{itemize}
\textit{Step 2.}
For any non-increasing sequence $\{b_k\}$ of non-negative numbers we have
\begin{equation}\label{A0}
\sum_{k=1}^{\infty} k^\frac{1}{n-1} b_k \leq \left( \sum_{k=1}^{\infty} b_k^\frac{n-1}{n} \right)^\frac{n}{n-1}. 
\end{equation}
Indeed
\[
\sum_{k=1}^{\infty} k^{\frac{1}{n-1}} b_k
=
\sum_{k=1}^{\infty}
\left(k\, b_k^{\frac{n-1}{n}}\right)^{\frac{1}{n-1}}
b_k^{\frac{n-1}{n}}\quad \text{and}\quad
k\, b_k^{\frac{n-1}{n}}
\le
\sum_{j=1}^{\infty} b_j^{\frac{n-1}{n}}.
\]
Thus, we have
\[
\sum_{k=1}^{\infty} k^{\frac{1}{n-1}} b_k
\le
\left(
\sum_{j=1}^{\infty} b_j^{\frac{n-1}{n}}
\right)^{\frac{1}{n-1}}
\sum_{k=1}^{\infty} b_k^{\frac{n-1}{n}}
=
\left(
\sum_{k=1}^{\infty} b_k^{\frac{n-1}{n}}
\right)^{\frac{n}{n-1}}.
\]We apply inequality \eqref{A0} with $b_k = |A^\varepsilon_{k-2}|$ in combination with \eqref{Lq_unif_final}. This yields
\begin{equation*}
\int_Q |(v_\varepsilon)_+|^\frac{n}{n-1} \, dx = \frac{n}{n-1}\int_0^\infty t^\frac{1}{n-1}|\{v_\varepsilon > t\}| \, dt \leq2 \sum_{k=1}^\infty k^\frac{1}{n-1} |\{v_\varepsilon > k-1\}| \leq 2\sum_{k=1}^\infty k^\frac{1}{n-1} |A^\varepsilon_{k-2}| \leq C.
\end{equation*}
The claim follows by arguing in a similar way for the negative part of $v_\varepsilon$.

\end{proof}
\subsection{Proof of the compactness result}\label{sottosezione_compattezza}
The proof follows the scheme described at the beginning of this section.

\begin{proof}[Proof of \Cref{thm:compactness}] Throughout, $p=\frac{n}{n-1}$ if $n\ge 2$ and $p=2$ if $n=1.$

\textit{Step 1. Strong compactness for the periodic case.} First we prove that if  $ \{u_\varepsilon\}\subset L^1(Q)$ satisfies $F^\Z_\varepsilon(u_\varepsilon)\leq M$, then there exists $b_\varepsilon\in \Z$ such that $\{u_\varepsilon-b_\varepsilon\}$ is precompact in $L^1(Q)$ and all its limit points are in $BV(Q;\Z)$.

Let \(b_{\varepsilon} \in \mathbb{Z}\) be chosen as in \Cref{lem:uniform_bound_1D} if $n=1$ and \Cref{lem:uniform_bound_nD} if $n\ge 2,$ and consider the translated sequence $v_\varepsilon:=u_\varepsilon-b_\varepsilon$. 
Since $F^\Z_\varepsilon$ is invariant under integer translations, $F_\varepsilon^\Z(v_\varepsilon)=F^\Z_\varepsilon(u_\varepsilon)\le M. $
This sequence is uniformly bounded in $L^{p}$; hence, up 
to extraction of a subsequence, $v_\varepsilon \rightharpoonup v$ weakly in $L^{p}(Q)$. To obtain strong convergence in $L^1$ we use a truncation argument. 
For all $k \in \mathbb{Z}$, we consider the truncation operator $T_k u = (u \vee k) \wedge (k + 1)$. We have
\begin{equation}
\frac{1}{\theta_\varepsilon} \left(\varepsilon^{2s} [T_k v_\varepsilon]^2_{H^s(Q)} +  \int_Q \operatorname{dist}^2(T_k v_\varepsilon, \{k, k + 1\}) \, dx \right) \leq  F^\Z_\varepsilon(v_\varepsilon) \leq M.
\end{equation}
By the compactness and lower bound results for double-well potentials (established in \cite{ABS1} for $s=\frac12$ and $n=1$, in \cite[Theorem 1.2]{Palatucci} for $s>\frac12$ and $n=1$, and in \cite[Theorem 1.3]{SavinValdinoci2010} for $n \ge 2$), for every $k \in \mathbb{Z}$ there exists a  function $\bar u_k \in BV(Q; \{k, k+1\})$ such that, up to the extraction of a subsequence, the sequence $T_k v_\varepsilon$ converges to $\bar u_k$ in $L^1(Q)$ and 
\[
c_0\operatorname{Per}(\{\bar u_k=k+1\};Q)\le\liminf_{\varepsilon\to0}F^\Z_\varepsilon(T_k v_\varepsilon),
\]
with $c_0>0$ depending only on $n$ and $s$. By a diagonal argument, the same subsequence works for every $k \in \Z.$
Now consider $N \in \mathbb{N}$ and the truncation $T^Nv = (v \vee (-N)) \wedge N$.
We can express this new truncation operator in terms of the previous operators \( T_k \) as
$$
T^Nv=T_{-N}v+\sum_{i=-N+1}^{N-1}(T_i v-i).
$$
Hence, there exists $v^N\in BV(Q;\mathbb{Z})$, namely $v^N=\bar u_{-N}+\sum_{i=-N+1}^{N-1}(\bar u_i-i)$, such that, up to a subsequence, 
$$
T^N\,v_\varepsilon\rightarrow v^N\quad\text{in }L^1(Q).
$$
Moreover, by Fatou's lemma and the truncation inequality in \eqref{eq:truncation},
\[
|Dv^N|(Q)\le\sum_{k=-N}^{N-1}\operatorname{Per}(\{\bar u_k=k+1\};Q)\le\frac1{c_0}\liminf_{\varepsilon\to0}\sum_{k\in\Z}F^\Z_\varepsilon(T_k v_\varepsilon)\le\frac{M}{c_0}.
\]
For each fixed $N\in \mathbb{N}$, we have 
$$
\begin{aligned}
\|v^N-v\|_{L^1(Q)}&\leq \liminf_{\varepsilon\rightarrow 0}\|T^N v_\varepsilon-v_\varepsilon\|_{L^1(Q)}=
\liminf_{\varepsilon\rightarrow 0}\int_{\{|v_\varepsilon|>N\}} (|v_\varepsilon|-N)\,dx\\&
\leq
\liminf_{\varepsilon\rightarrow 0} \int_{\{|v_\varepsilon|>N\}} |v_\varepsilon|\,dx\leq \frac{C}{N^{p-1}}.
\end{aligned}
$$
In the first inequality we used the weak lower semicontinuity of the $L^1$ norm (since $T^Nv_\varepsilon-v_\varepsilon\rightharpoonup v^N-v$) and in the last one we used that on the set $\{|v_\varepsilon| > N\}$ we have $|v_\varepsilon| \leq |v_\varepsilon|^{p}/N^{p-1}$ and the $L^{p}(Q)$  bound.
In conclusion, we have
$$
\begin{aligned}\|v_\varepsilon-v\|_{L^1(Q)}&\le\|v_\varepsilon-T^Nv_\varepsilon\|_{L^1(Q)}+\|T^Nv_\varepsilon-v^N\|_{L^1(Q)}+\|v^N-v\|_{L^1(Q)}\\&\le\frac{2C}{N^{p-1}}+\|T^Nv_\varepsilon-v^N\|_{L^1(Q)}.\end{aligned}
$$
Letting first $\varepsilon\to0$ and then $N\to\infty$, we get $v_\varepsilon$ convergent to $v$ in $L^1(Q)$.
Thus, $v^N=T^Nv\in BV(Q;\mathbb{Z})$ with $|DT^Nv|(Q)\le M/c_0$ for every $N,$ and then $v\in BV(Q;\mathbb{Z})$ by lower semicontinuity of the total variation as $N \to \infty$.

\textit{Step 2. Bi-Lipschitz rectification and extension to the general potential.} We now reduce
the general potential to the periodic one.
 Let $w_z := s_z - z$. 
We define the function $\phi^S: \mathbb{R} \to \mathbb{R}$ as follows:
\begin{equation}
\phi^S(x) := \begin{cases} 
z + \displaystyle\frac{x - z}{2w_z} & \text{if } x \in [z, s_z) \\[10pt]
z + \displaystyle\frac{1}{2} + \frac{x - s_z}{2(1 - w_z)} & \text{if } x \in [s_z, z+1)
\end{cases}
\end{equation}
where $z = \lfloor x \rfloor$. Note that $\phi^S(\mathbb{Z} \cup S) = \frac{1}{2}\mathbb{Z}$ and $\phi^S$ is a strictly increasing continuous function.
The function $\phi^S$ is bi-Lipschitz on $\mathbb{R}$. Specifically by \eqref{d} for all $x, y \in \mathbb{R}$, $$\frac{1}{2(1-d) }|x - y|\leq |\phi^S(x) - \phi^S(y)| \le \frac{1}{2d} |x - y|.$$
Since $2\phi^S$ is $\frac1d$-Lipschitz and maps $\Z \cup S $ onto $\Z$, for every $t \in \R$
\begin{equation}
 \operatorname{dist}^2(t,\Z\cup S) \geq d^2 \operatorname{dist}^2(2\phi^S(t),\Z)  ,
\end{equation}
and therefore, by the lower estimate in \eqref{eq:two_sided_estimate},
$$
F_\varepsilon(u_\varepsilon)\geq F^{\Z\cup S}_\varepsilon(u_\varepsilon)\geq d^2 F^\Z_\varepsilon(2 \phi^S\circ u_\varepsilon).
$$
By \textit{Step 1} applied to $2\phi^S\circ u_\varepsilon$, there exists a sequence  $\{b_\varepsilon\} \subset \mathbb{Z}$ such that, up to a subsequence,
\begin{equation}
    \psi_\varepsilon := 2\phi^S(u_\varepsilon) - b_\varepsilon \longrightarrow \psi \in BV(Q;\Z)\quad \text{strongly in } L^1(Q)\text{ and a.e.}
\end{equation}

\textit{Step 3. Compactness up to translations.} 
Write $b_\varepsilon=2a_\varepsilon+r_\varepsilon$ with $a_\varepsilon\in \Z$ and $r_\varepsilon\in \{0,1\}$. We prove that $\{u_\varepsilon-a_\varepsilon\}$ is relatively compact.
Up to a subsequence, $r_\varepsilon = r \in \{0, 1\}$ is constant.
For each $k \in \Z $ set $g_k(y) := (\phi^S)^{-1}(y + k) - k$. Then
$$
u_\varepsilon-a_\varepsilon=g_{a_\varepsilon}(y_\varepsilon), \qquad y_\varepsilon:=\frac{\psi_\varepsilon+r}{2}\to y :=\frac{\psi+r}{2}\in \tfrac12\Z \quad \text{a.e.}
$$
Note that for every $k \in \mathbb{Z}$,
$g_k(0) = (\phi^S)^{-1}(k) - k = k - k = 0$ and each $g_k$ is  2-Lipschitz; hence $|g_k(y)|\le 2|y|.$ In particular $|u_\varepsilon-a_\varepsilon|\le 1+|\psi_\varepsilon|$, so $\{u_\varepsilon-a_\varepsilon\}$ satisfies the bounds in the statement.
For each $m \in \Z$ the sequence $\{g_{a_\varepsilon}(m/2)\}$ is bounded, so up to a subsequence (independent of $m$, by a diagonal argument)
$$
g_{a_\varepsilon}(m/2)\to\Gamma(m)\in \R, \quad \text{for each }m\in \Z.
$$
Thus
$$
u_\varepsilon-a_\varepsilon\to \Gamma(\psi+r)\quad \text{strongly in }L^1(Q)\text{ and a.e.}
$$
Finally, $|\Gamma(m+1)-\Gamma(m)|\le 1$ because $g_k$ is $2$-Lipschitz. By the coarea formula, since $\psi\in BV(Q;\Z),$ we get
$$
\displaystyle|D[\Gamma(\psi+r)]|(Q)\le\sum_{j\in\Z}|\Gamma(j+1+r)-\Gamma(j+r)|\operatorname{Per}(\{\psi>j\};Q)\le|D\psi|(Q)<\infty .
$$
Note that $\Gamma(2z)=z$ and $\Gamma(2z+1)=\lim_\varepsilon(s_{z+a_\varepsilon}-a_\varepsilon)=z+\lim_\varepsilon w_{z+a_\varepsilon}$, which gives the set $S^*$ of \Cref{rem:compactness}.

\textit{Step 4. Limit in $\Z \cup S$.} We want to show that if $u_\varepsilon$ converges to $u$ in $L^1$ with  $F_\varepsilon(u_\varepsilon)\leq M$ then $u \in BV(Q;\Z\cup S )$.
We have
\[
|a_\varepsilon|\,|Q|\le\|u_\varepsilon\|_{L^1(Q)}+\|u_\varepsilon-a_\varepsilon\|_{L^1(Q)},
\]
and the right-hand side is bounded by the convergence in $L^1$ of $u_\varepsilon$ and \textit{Step 3}. Hence, up to a subsequence, $a_\varepsilon=a$ for all $\varepsilon$. Then $\Gamma(m)=g_a(\frac m2)$ and
\[
u=\lim_{\varepsilon\to0}u_\varepsilon=a+g_a(y)=(\phi^S)^{-1}(y+a).
\]
Since $y+a\in\frac12\Z$ a.e. and $(\phi^S)^{-1}(\frac12\Z)=\Z\cup S$, and $u\in BV(Q)$ by \textit{Step 3}, we get $u\in BV(Q;\Z\cup S)$.

\textit{Step 5. Regime \hyperlink{R2}{(R2)}.} We want to prove that if $u_\varepsilon\to u$ in $L^1(Q)$
and $\delta_\varepsilon / \theta_\varepsilon \to \infty$, then $u \in BV(Q;\Z).$
Fix $\delta \in (0,1)$. By
\eqref{eq:two_sided_estimate} and $\delta_\varepsilon / \theta_\varepsilon \to \infty$ we have  
\[
 |A_{\delta}(u_\varepsilon)| \le \frac{M\theta_\varepsilon}{\delta_\varepsilon }\to 0,
\]
where $A_\delta(u_\varepsilon)$ is defined in \eqref{eq-set-u-comparison}.
On the other hand, by \eqref{d}, 
$$
\{x \in Q\,:\,u(x)\in S \text{ and }|u_\varepsilon(x)-u(x)|<\delta d/2\}\subset A_\delta(u_\varepsilon).
$$
Therefore
\[
|\{u\in S\}|\le|A_\delta(u_\varepsilon)|+|\{|u_\varepsilon-u|\ge\delta d/2\}|\longrightarrow0,
\]
where for the second term we use convergence in measure. Hence $|\{u\in S\}|=0$, and by \textit{Step 4}, $u\in BV(Q;\Z)$.
\end{proof}
 \section{The optimal profile problem for \texorpdfstring{$s > \frac12$}{s > 1/2}}\label{Section:optimal}
In this section, we fix $s\in\left(\frac12,1\right)$, $n \ge 1$ and we set $W^S=\operatorname{dist}^2(\cdot,\Z \cup S)$.
We show that the jump cost $\lambda_{s,n}\kappa_s\Phi_{\Z\cup S}$ of \Cref{Jump_costs} is the cost of the optimal transition for the potential $W^S, $ and we introduce the truncated minimizer $v^T$ used for the recovery sequence in \Cref{sec:upper_bound}.

For $v \in H_{\rm loc}^s(\R)$, let 
\begin{equation}\label{E_s}
    \mathcal{E}_{s,n}(v):=C_{s,n}[v]^2_{H^s(\R)}+\int_\R W^S(v)\,dx,
\end{equation}
where $C_{s,n}$ is the constant defined in \eqref{dhef:C_ns} and $[v]_{H^s(\R)}$ is the Gagliardo seminorm on $\R$ (possibly $+\infty$). 
For $a,b \in \Z\cup S$ and $T>0$, we set 
\begin{equation}\label{theta_s,n}
    \Theta_{s,n}(a,b):=\inf\left\{
     \mathcal{E}_{s,n}(v)\,:\,v \in H_{\rm loc}^s(\R),\,\lim_{x\to -\infty}v(x)=a,\,\lim_{x\to +\infty}v(x)=b
    \right\},
\end{equation}
and
\begin{equation}\label{theta_s_T}
    \Theta_{s,n}^T(a,b):=\inf\left\{
     \mathcal{E}_{s,n}(v)\,:\,v \in H_{\rm loc}^s(\R),\,v=a\, \text{ on }(-\infty,-T],\,v=b \text{ on }[T,+\infty)
    \right\}.
\end{equation}

We prove the following result, which is the general version of the observation presented in \Cref{rmk:slicing_constant} for a jump between any two values of $\union$.
\begin{prop}\label{costo_theta}
    Let $a,b  \in\Z \cup S$ with $a<b.$ Then 
    \begin{equation}
        \Theta_{s,n}(b,a)=\Theta_{s,n}(a,b)=\lambda_{s,n}\kappa_s\Phi_{\Z\cup S}(a,b),
    \end{equation}
    where $\lambda_{s,n}$ is defined in \eqref{dhef:C_ns}, $\kappa_s$ in \eqref{eq:kappa_s} and $\Phi_{\Z\cup S}$ in \Cref{one-dim_jump}. 
    Moreover, if $v^T$ is a minimizer of \eqref{theta_s_T} (see \Cref{lemma:optimal_profile}), then
    \begin{equation}
        \lim_{T\to \infty}\mathcal{E}_{s,n}(v^T)=\lambda_{s,n}\kappa_s\Phi_{\Z\cup S}(a,b).
    \end{equation}
\end{prop}
The proposition follows from the following properties of $\Theta_{s,n}.$
\begin{lemma}\label{lemma:optimal_profile}
Let $s \in \left(\frac12, 1\right)$ and $a, b \in \Z \cup S$ with $a < b$. The following properties hold.
\begin{enumerate}
    \item  For any $T > 0$, the infimum $\Theta^T_{s,n}(a,b)$ is attained by a non-decreasing minimizer $v^T \in H^s_{\mathrm{loc}} (\mathbb{R})$, with $a\leq v^T\leq b$.
    \item  The truncated energy recovers the full-line infimum as $T \to \infty$, that is, 
    \[\lim_{T\to +\infty} \Theta^T_{s,n}(a,b) =\Theta_{s,n}(a,b) . \]
    \item If there exists $c \in (a, b)\cap (\Z \cup S)$, then the infimum is additive:
    \[ \Theta_{s,n}(a,b) = \Theta_{s,n}(a,c) + \Theta_{s,n}(c,b). \]
      \item  If there exists an intermediate well $c \in (a, b)\cap (\Z \cup S)$, then the infimum $\Theta_{s,n}(a,b)$ is not attained, otherwise it is attained.
\end{enumerate}
\end{lemma}
\begin{proof}
\textit{(1)} Since the truncation $(v\vee a)\wedge b$ does not increase the seminorm and $W^S(a)=W^S(b)=0$, we can take a minimizing sequence $\{v_k\}$, for \eqref{theta_s_T}, with $a\le v_k\le b$. Then $\{v_k\}$ is bounded in $H^s((-2T,2T))$ and, since $v_k$ is fixed outside $[-T,T]$, up to a subsequence it converges in $L^2_{\mathrm{loc}}(\R)$ and a.e.\ in $\R$. By Fatou's lemma the limit is a minimizer. Finally, the non-decreasing rearrangement of a minimizer on $[-T,T]$ is still admissible, has the same potential energy and does not increase the seminorm (see \cite{Garsia}); hence we can take $v^T$ non-decreasing.

\textit{(2)} 
Since $\Theta^T_{s,n}(a,b)$ is non-increasing in $T$ and is bounded from below by $\Theta_{s,n}(a,b),$ the limit exists and is greater than or equal to $\Theta_{s,n}(a,b).$ The opposite inequality 
follows by a standard argument for fractional transition layers; we refer to \cite[Proposition 3.1]{Palatucci} for the detailed construction, which only uses that $W^S$ vanishes at $a$ and $b$ and is quadratic near the wells.

\textit{(3)} Let $c\in(a,b)\cap (\Z \cup S)$. 
We first prove that
\begin{equation}
    \label{eq-lemma-subadd}
    \Theta_{s,n}(a,b) \leq \Theta_{s,n}(a,c)+\Theta_{s,n}(c,b).
\end{equation}
Let $\eta>0$. By definition of the infimum and by the observation in the proof of \textit{(1)-(2)}, there exist $T>0$ and two non-decreasing functions $v_1,v_2$, with $a\le v_1\le c\le v_2\le b$, such that
\[
\mathcal{E}_{s,n}(v_1)\le \Theta_{s,n}(a,c)+\eta,
\qquad
\mathcal{E}_{s,n}(v_2)\le \Theta_{s,n}(c,b)+\eta,
\]
with 
\[
v_1(x)=a \,\,\text{for } x\le -T,
\,\,\,
v_1(x)=c \,\,\text{for } x\ge 0
\quad
\hbox{and}
\quad
v_2(x)=c \,\, \text{for } x\le 0,
\,\,\,
v_2(x)=b \,\, \text{for } x\ge T .
\]
Given $L>0$ we construct a competitor by placing the two profiles far apart:
\[
g_L(x)=
\begin{cases}
v_1(x), & x\le 0,\\[2mm]
v_2(x-L), & x\ge 0.
\end{cases}
\]
The energy decomposes as
\[
\mathcal{E}_{s,n}(g_L)
=
\mathcal{E}_{s,n}(v_1)
+
\mathcal{E}_{s,n}(v_2)
+
2 I_L,
\]
where the interaction term between the two profiles is
\[
I_L
=
 C_{s,n}\int_{\mathbb{R}}\int_{\mathbb{R}}
\frac{(v_1(x)-v_1(y))(v_2(x-L)-v_2(y-L))}
{|x-y|^{1+2s}}
\,dx\,dy.
\]
An easy computation shows that $I_L$ tends to $0$ as $L\to\infty$, since $s>\frac12$.
Since $g_L$ is admissible for $\Theta_{s,n}(a,b)$, we get
\[
\Theta_{s,n}(a,b)\le\mathcal E_{s,n}(g_L)\le\Theta_{s,n}(a,c)+\Theta_{s,n}(c,b)+2\eta+2I_L.
\]
Letting $L\to +\infty$ and $\eta \to 0$ we have \eqref{eq-lemma-subadd}.

Now we prove the opposite inequality,
\begin{equation}
    \label{eq-lemma-supadd}
    \Theta_{s,n}(a,b) \geq \Theta_{s,n}(a,c)+\Theta_{s,n}(c,b).
\end{equation}
Let $v$ be an admissible function in the definition of $\Theta_{s,n}(a,b)$ with finite energy,
define
\[
w_1:=\min\{v,c\},
\qquad
w_2:=\max\{v,c\},
\]
with $w_1$ admissible for $\Theta_{s,n}(a,c)$ and $w_2$ admissible for $\Theta_{s,n}(c,b).$ Moreover, since $W^S(c)=0,$
\[
\int_{\mathbb{R}} W^S(v)\,dx
=
\int_{\mathbb{R}} W^S(w_1)\,dx
+
\int_{\mathbb{R}} W^S(w_2)\,dx .
\]
Moreover, for each $x,y \in \R$,
\begin{equation}\label{eq:pointwise_identity}
|v(x) - v(y)|^2 = |w_1(x) - w_1(y)|^2 + |w_2(x) - w_2(y)|^2 + 2(c - w_1(x))(w_2(y) - c)+ 2(c - w_1(y))(w_2(x) - c).
\end{equation}
Hence, we obtain
\[
[v]_{H^s}^2
=
[w_1]_{H^s}^2
+
[w_2]_{H^s}^2
+
4\int_{\mathbb{R}}\int_{\mathbb{R}}
\frac{\,(c-w_1(x))(w_2(y)-c)}{|x-y|^{1+2s}}
\,dx\,dy
>
[w_1]_{H^s}^2 + [w_2]_{H^s}^2 ,
\]
where the strict inequality follows by the fact that $
c-w_1>0$ and $w_2-c>0$ on sets with positive measure.
Hence
\begin{equation}\label{optimal_4}
\mathcal{E}_{s,n}(v)>\mathcal{E}_{s,n}(w_1)+\mathcal{E}_{s,n}(w_2)\geq \Theta_{s,n}(a,c)+\Theta_{s,n}(c,b),
\end{equation}
and taking the infimum over $v$ we get \eqref{eq-lemma-supadd}.

\textit{(4)}
The strict inequality in \eqref{optimal_4} shows that the minimum is not achieved if an intermediate well is present. Moreover,
if no intermediate wells are present then the infimum $\Theta_{s,n}(a,b)$ is attained by a strictly increasing profile (see \cite[Proposition~3.3]{Palatucci}).
\end{proof}
\begin{proof}[Proof of \Cref{costo_theta}]
Let $a=p_0<p_1<\dots<p_N=b$ be the points in $[a,b]\cap (\Z \cup S)$;
by \Cref{lemma:optimal_profile} \textit{(3)}, we have
$$
\Theta_{s,n}(a,b)=\sum_{i=0}^{N-1} \Theta_{s,n}(p_{i},p_{i+1}).
$$
Since truncating at $p_{i}$ and $p_{i+1}$ does not increase $\mathcal{E}_{s,n}$, we can rewrite $ \Theta_{s,n}(p_i,p_{i+1})$ as
    $$
    \Theta_{s,n}(p_i,p_{i+1})=\inf\left\{C_{s,n}[v]^2_{H^s(\R)}+\int_\R \operatorname{dist}^2(v,\{p_i,p_{i+1}\})\,dx\,:\,v(-\infty)=p_i,\,v(+\infty)=p_{i+1}\right\}.
    $$
    Let $w\in H^s_{\rm loc}(\R)$ be an admissible profile for the above problem, we consider the function  $v(x)=p_i+(p_{i+1}-p_i)\,w(x/\lambda_{s,n})$.
Since $[w(\cdot/\lambda)]^2_{H^s(\R)}=\lambda^{1-2s}[w]^2_{H^s(\R)}$ and $C_{s,n}\lambda_{s,n}^{1-2s}=\lambda_{s,n},$ we get
$$
\mathcal{E}_{s,n}(v)=\lambda_{s,n}(p_{i+1}-p_i)^2\left(
[w]^2_{H^s(\R)}+\int_\R\operatorname{dist}^2(w,\{0,1\})\,dx
\right).
$$
The map $w \mapsto v$ is a bijection between the admissible profiles from $0$ to $1$ and those from $p_i$ to $p_{i+1};$ hence, by definition of $\kappa_s$ in \eqref{eq:kappa_s}, $\Theta_{s,n}(p_i,p_{i+1})=\lambda_{s,n}\kappa_s(p_i-p_{i+1})^2.$
Finally, $\mathcal E_{s,n}(v^T)=\Theta^T_{s,n}(a,b)\to\Theta_{s,n}(a,b)$ by (2) in \Cref{lemma:optimal_profile}.
\end{proof}
\section{Lower bound}\label{sec:lower-bound}
In this section, we establish the lower bound for the functional $F_\varepsilon$.
In what follows, we set $W^S:=\operatorname{dist}^2(\cdot,\union)$ and $\mathcal{V}$ the set of wells of the limit ($\mathcal V=\Z\cup S$ in regime \hyperlink{R1}{(R1)}, $\mathcal V=\Z$ in regime \hyperlink{R2}{(R2)}). 
The proof is divided into three steps.

\textit{One-dimensional lower bound}. For each regime, we prove the lower bound for a single transition between two consecutive points $a<b$ in $\mathcal{V}$. This is the step where the jump cost is computed:
\begin{itemize}
    \item  \emph{Regime \hyperlink{R1}{(R1)}} (\Cref{lemma:lb_1d_critico}). By the lowest estimate \eqref{eq:two_sided_estimate}, it suffices to find the lower bound of the reference energy $F^{\Z\cup S}_\varepsilon$, whose potential $W^S$ does not depend on $\varepsilon$; this reduction is carried out in \Cref{sub:proof_lower}. Between two consecutive points of $\Z\cup S$, the potential is the classic double-well potential, whose one-dimensional $\Gamma$-limit was computed in \cite{ABS1} for $s=\frac12$ and in \cite{SavinValdinoci2010,Palatucci} for $s>\frac12$.

\item  \emph{Regime \hyperlink{R2}{(R2)}, $s>\frac12$} (\Cref{lemma:lb_1d_sub}). By \eqref{eq:two_sided_estimate}, the analysis reduces to the reference energy $F_\varepsilon^{\mathbb Z\cup S}$. A transition from $k$ to $k+1$ costs as much as a double transition $k \to s_k\to k+1.$ In particular, the height $\delta_\varepsilon$ of the well at $s_k$ plays no role.

\item \emph{Regime \hyperlink{R2}{(R2)}, $s=\frac12$} (\Cref{lemma:lb_1d_sub}).
This is the genuinely new case. The comparison with $F^{\Z\cup S}_\varepsilon$ is still valid, but no longer sharp. Heuristically, an almost optimal transition from $k$ to $k+1$ rests at $s_k$ on an interval of length about $\varepsilon/\delta_\varepsilon$. Consequently, the transition cost between the integer minima $k$ and $k+1$ must be explicitly calculated, leading to a limiting energy that depends on the parameter $\ell$ in \eqref{ell}.
\end{itemize}

\textit{Dimension $n$} (\Cref{lb_mu}).  We slice in all directions and apply the one-dimensional lemmas on each line. On the slices the constants in front of the seminorm and of the potential are no longer $1$, so we use a version of the one-dimensional lemmas with arbitrary constants, which follows from them by scaling  (\Cref{lemma:slicing}).

\textit{Conclusion} (\Cref{sub:proof_lower}). Truncations between consecutive points of $\mathcal V$ reduce the general case to a single transition. In regime \hyperlink{R1}{(R1)}, the lower estimate in \eqref{eq:two_sided_estimate} also produces the bulk term $K|\{u\in S\}|$.

\subsection{Lower bounds for the one-dimensional case}\label{sub:lower_1d}
We first recall the classical lower-bound results for two-well potentials, applied to $F^{\Z\cup S}_\varepsilon.$
\begin{lemma}[Lower bound for $F^{\Z\cup S}_\varepsilon$]\label{lemma:lb_1d_critico} Let $n=1$. Let $v_i$ and $v_{i+1}$ be two consecutive points of $\Z\cup S.$
Let $I\subset \R$ be a bounded open interval and  \( u \in BV(I; \{v_i,v_{i+1}\}) \). Let \( \{u_\varepsilon \} \subseteq L^1(I; [v_i,v_{i+1}]) \) be a sequence such that \( u_\varepsilon \) converges to \( u \) in \( L^1(I) \). Then
\begin{equation}
\liminf_{\varepsilon \to 0} F^{\Z\cup S}_\varepsilon(u_\varepsilon;I) \geq \,\Psi_{\mathrm{sup}}(v_i,v_{i+1})\,\mathcal{H}^{0}(J_u)=\kappa_s(v_{i+1}-v_i)^2\,\mathcal{H}^{0}(J_u).
\end{equation}
\end{lemma}
\begin{proof}
 Since the sequence $\{u_\varepsilon\}$ takes values in the interval $[v_i, v_{i+1}]$, with $v_i,v_{i+1}$ two consecutive points of $\Z\cup S$, we have $\operatorname{dist}^2(u_\varepsilon,\Z\cup S)=\operatorname{dist}^2(u_\varepsilon,\{v_i,v_{i+1}\})$. Then, the claim is the classical
one-dimensional lower bound for double-well potentials.

\emph{Case $s>\frac12$.} This is the fractional counterpart of the Modica--Mortola theorem,
proved in \cite{SavinValdinoci2010,Palatucci}; the lower bound is given by the optimal transition cost $\Theta_{s,1}(v_i, v_{i+1})$ defined in \eqref{theta_s,n}, as shown in \cite[Proposition 5.1]{Palatucci}. Hence the thesis follows by \Cref{costo_theta}.

\emph{Case $s=\frac12$.} The lower bound $2(v_{i+1}-v_i)^2\,\mathcal H^0(J_u)$ is proved in \cite{ABS1},
and $\kappa_{1/2}=2$. Here the energy does not exhibit equipartition, and the limit is determined solely by the $H^{1/2}$ interactions between the two wells of the potential.
\end{proof}
The following lemma addresses the specific features of regime \hyperlink{R2}{(R2)}. For $s = \frac12$, the presence of the additional well $s_k$ modifies the limiting energy cost.
\begin{lemma}[Regime \hyperlink{R2}{(R2)} -- Wells in $\Z$]\label{lemma:lb_1d_sub} 
Let $n=1$ and assume we are in regime  \hyperlink{R2}{(R2)}. Let $I\subset \R$ be a bounded open interval, \( u \in BV(I; \{k,k+1\}) \) with $k \in \Z$, and \( \{u_\varepsilon \} \subseteq L^1(I; [k,k+1]) \) be a sequence such that \( u_\varepsilon \) converges to \( u \) in \( L^1(I) \). Then
\begin{equation}\displaystyle
\liminf_{\varepsilon \to 0} F_\varepsilon(u_\varepsilon;I)
\geq\mathcal{H}^{0}(J_u)\,\Psi_{\mathrm{sub}}(k,k+1)=
\mathcal{H}^{0}(J_u)
\begin{cases}
\kappa_s(1-2(s_k-k)(k+1-s_k)), & \text{if } s>\frac12,\\[8pt]
2-4\ell(s_k-k)(k+1-s_k), & \text{if } s=\frac12.
\end{cases}
\end{equation}
\end{lemma}
\begin{proof} 
\textit{Case $s > \frac12$.}\quad For $s > \frac12$, the thesis follows directly from the estimate $F_\varepsilon(u_\varepsilon) \ge F^{\Z\cup S}_\varepsilon(u_\varepsilon)$. Truncating the function in $[k,s_k]$ and $[s_k,k+1]$ we can apply  \Cref{lemma:lb_1d_critico} and conclude that 
$$
\liminf_{\varepsilon \to 0} F^{\Z\cup S}_\varepsilon(u_\varepsilon) \geq \,\Psi_{\mathrm{sup}}(k,k+1)\,\mathcal{H}^{0}(J_u)=\kappa_s((s_k-k)^2+(k+1-s_k)^2)\,\mathcal{H}^{0}(J_u).
$$
\textit{Case $s=\frac12.$
} It is sufficient to prove the result for a jump set consisting of a single point, $J_u=\{x\}$. Indeed, for a general jump set, we can consider a finite collection of disjoint intervals $\{I_j\}$, each containing exactly one jump point of $u$, and then exploit the superadditivity of the energy:
\begin{equation}
    F_\varepsilon(u_\varepsilon; I) \ge \sum_{j} F_\varepsilon(u_\varepsilon; I_j).
\end{equation}
We can reduce the problem to the case where $u_\varepsilon$ is a sequence of monotone non-decreasing functions via rearrangement, since the potential term $\int_I W_\varepsilon(u_\varepsilon) \, dx$ is invariant under rearrangement, while the fractional Gagliardo seminorm decreases by rearrangement (see, e.g., \cite{Garsia}).
Let $w_k:=s_k-k.$
    Fix $\delta \in (0,d/4)$ with $d$ as in \eqref{d}, let us consider the following disjoint sets
    $$
    A_\varepsilon:=\{u_\varepsilon\in [k,k+\delta)\},
   \qquad
    B_\varepsilon:=\{u_\varepsilon\in (k+1-\delta,k+1]\},
    $$
    $$
    C_\varepsilon:=\{u_\varepsilon\in (s_k-\delta,s_k+\delta)\},
    \qquad
    D_\varepsilon:=I\setminus (A_\varepsilon\cup B_\varepsilon\cup C_\varepsilon).
    $$
    We can easily estimate the potential term as follows
$$
      \int_I W_\varepsilon(u_\varepsilon(x))\,dx\ge \int_{D_\varepsilon} W_\varepsilon(u_\varepsilon(x))\,dx+\int_{C_\varepsilon} W_\varepsilon(u_\varepsilon(x))\,dx\ge |D_\varepsilon|\inf_{D_\varepsilon}W_\varepsilon \circ u_\varepsilon+|C_\varepsilon|\inf_{C_\varepsilon}W_\varepsilon \circ u_\varepsilon.
$$
  Now, for small $\varepsilon$ we have
  $$
  \inf_{D_\varepsilon}W_\varepsilon \circ u_\varepsilon\ge\min\{W_\varepsilon(x)\,:\,x \in [k+\delta,s_k-\delta]\cup [s_k+\delta,k+1-\delta]\}=\delta^2
  $$
  and
  $$
  \inf_{C_\varepsilon}W_\varepsilon \circ u_\varepsilon\ge\min\{W_\varepsilon(x)\,:\,x \in [s_k-\delta, s_k+\delta]\}=\delta_\varepsilon.
  $$
  Hence,
  \begin{equation}\label{potenziale}
       \int_I W_\varepsilon(u_\varepsilon(x))\,dx\ge|D_\varepsilon|\delta^2+|C_\varepsilon|\delta_\varepsilon.
  \end{equation}
For the Gagliardo seminorm, we have
\begin{equation}
\begin{split}
&\frac12\int_I\int_I \frac{|u_\varepsilon(x)-u_\varepsilon(y)|^2}{|x-y|^2} dx\, dy  \\
&\ge  
\int_{A_\varepsilon}\int_{B_\varepsilon} \frac{|u_\varepsilon(x)-u_\varepsilon(y)|^2}{|x-y|^2}dx\, dy+\int_{C_\varepsilon}\int_{B_\varepsilon} \frac{|u_\varepsilon(x)-u_\varepsilon(y)|^2}{|x-y|^2}dx\, dy+\int_{A_\varepsilon}\int_{C_\varepsilon} \frac{|u_\varepsilon(x)-u_\varepsilon(y)|^2}{|x-y|^2}dx\, dy
\\
&\ge 
(1-2\delta)^2 \int_{A_\varepsilon}\int_{B_\varepsilon}\frac{1}{|x-y|^2}dx\, dy+(1-w_k-2\delta)^2\int_{C_\varepsilon}\int_{B_\varepsilon}\frac{1}{|x-y|^2} dx\, dy\\
&\qquad + (w_k-2\delta)^2\int_{A_\varepsilon}\int_{C_\varepsilon}\frac{1}{|x-y|^2}dx\, dy.
\end{split}  
\end{equation}
In order to estimate each integral we use the explicit one-dimensional computation:
if \(\alpha_1 < \beta_1 < \alpha_2 < \beta_2\), then
\begin{equation}\label{conto_espl}
\begin{aligned}
\int_{\alpha_1}^{\beta_1} \int_{\alpha_2}^{\beta_2} \frac{1}{|x - y|^2} \, dx \, dy
= \log \left( \frac{(\beta_2 - \beta_1)(\alpha_2 - \alpha_1)}{(\beta_2 - \alpha_1)(\alpha_2 - \beta_1)} \right).
\end{aligned}
\end{equation}
   Therefore we get
$$
\int_{A_\varepsilon}\int_{B_\varepsilon}\frac{1}{|x-y|^2}=\log\left(\frac{(|A_\varepsilon|+|C_\varepsilon|+|D_\varepsilon|)(|B_\varepsilon|+|C_\varepsilon|+|D_\varepsilon|)}{(|C_\varepsilon|+|D_\varepsilon|)|I|}\right).
$$
  To estimate the last two terms, we rely on a rearrangement argument. Since $u_\varepsilon$ transitions from $k$ to $k+1$, the corresponding sets are spatially ordered in $I$ as $A_\varepsilon$, followed by a portion of $D_\varepsilon$, then $C_\varepsilon$, the remaining part of $D_\varepsilon$, and finally $B_\varepsilon$. 
Noting that the interaction integral decreases as the distance between the sets increases, a lower bound for the double integral over $A_\varepsilon \times C_\varepsilon$ (and similarly over $B_\varepsilon \times C_\varepsilon$) is obtained by placing the entire measure of the transition set $D_\varepsilon$ between them. A rigorous justification of this argument, based on the properties of monotone rearrangements, is provided in \cite[Lemma 2]{ABS1}. This leads to the following estimates
$$
\int_{B_\varepsilon}\int_{C_\varepsilon}\frac{1}{|x-y|^2}\ge\log\left(\frac{(|C_\varepsilon|+|D_\varepsilon|)(|B_\varepsilon|+|D_\varepsilon|)}{(|C_\varepsilon|+|D_\varepsilon|+|B_\varepsilon|)|D_\varepsilon|}\right),
$$
  and 
$$
\int_{A_\varepsilon}\int_{C_\varepsilon}\frac{1}{|x-y|^2}\ge\log\left(\frac{(|C_\varepsilon|+|D_\varepsilon|)(|A_\varepsilon|+|D_\varepsilon|)}{(|C_\varepsilon|+|D_\varepsilon|+|A_\varepsilon|)|D_\varepsilon|}\right).
$$
Now, by the convergence in measure of the sequence $u_\varepsilon$ we have
$$
\liminf_{\varepsilon\to 0 }|A_\varepsilon|\ge \liminf_{\varepsilon\to 0 }|\{u=k\}\cap \{|u_\varepsilon-u|<\delta\}|\ge |\{u=k\}|-\lim_{\varepsilon\to 0 } |\{|u_\varepsilon-u|\ge \delta\}|=|\{u=k\}|.
$$
Similarly
\begin{gather*}
\begin{aligned}
\liminf_{\varepsilon\to 0}|B_\varepsilon|
&\ge \liminf_{\varepsilon\to 0}|\{u=k+1\}\cap\{|u_\varepsilon-u|<\delta\}| \\
&\ge |\{u=k+1\}| - \lim_{\varepsilon\to 0}|\{|u_\varepsilon-u|\ge\delta\}|
= |\{u=k+1\}|,
\end{aligned} \\[6pt]
\lim_{\varepsilon\to 0}|C_\varepsilon|
\le \lim_{\varepsilon\to 0}|\{|u_\varepsilon-u|\ge\delta\}| = 0
\quad\text{and}\quad
\lim_{\varepsilon\to 0}|D_\varepsilon|
\le \lim_{\varepsilon\to 0}|\{|u_\varepsilon-u|\ge\delta\}| = 0.
\end{gather*}
Using the estimate from below of the integrals and these limits, we can estimate the seminorm as follows
\begin{equation}\label{seminorma}
\begin{aligned}
\int_I\int_I & \frac{|u_\varepsilon(x)-u_\varepsilon(y)|^2}{|x-y|^2} \, dx \, dy \\[10pt]
&\ge 2 \left[ (1-2\delta)^2 \log \left( \frac{1}{|C_\varepsilon|+|D_\varepsilon|} \right) + \left( (1-w_k-2\delta)^2 + (w_k-2\delta)^2 \right) \log \left( \frac{|C_\varepsilon|+|D_\varepsilon|}{|D_\varepsilon|} \right) \right] + O(1) \\[10pt]
&= -\left( 4w_k(1-w_k) - 8\delta^2 \right) \log(|C_\varepsilon| + |D_\varepsilon|) - 2 \left[ (1-w_k-2\delta)^2 + (w_k-2\delta)^2 \right] \log |D_\varepsilon| + O(1) \\[10pt]
&= -\left( 4w_k(1-w_k) - 8\delta^2 \right) \log(|C_\varepsilon| + |D_\varepsilon|) - (2 -4w_k(1-w_k)- 8\delta + 16\delta^2) \log |D_\varepsilon| + O(1).
\end{aligned}
\end{equation}

Using estimates \eqref{seminorma} and \eqref{potenziale}, we obtain
\begin{equation}
\begin{aligned}
F_\varepsilon(u_\varepsilon) \ge \, & \frac{1}{|\log\varepsilon|} \Biggl[ -\left( 4w_k(1-w_k) - 8\delta^2 \right) \log(|C_\varepsilon| + |D_\varepsilon|) - (2 -4w_k(1-w_k)- 8\delta + 16\delta^2) \log |D_\varepsilon| \\
& \qquad \quad + \frac{|D_\varepsilon|\delta^2}{\varepsilon} + \frac{|C_\varepsilon|\delta_\varepsilon}{\varepsilon} \Biggr] + o(1)= \\[10pt]
=& \frac{4w_k(1-w_k) - 8\delta^2}{|\log\varepsilon|} \left( -\log(|C_\varepsilon| + |D_\varepsilon|) + (|C_\varepsilon| + |D_\varepsilon|) \frac{\delta_\varepsilon}{\varepsilon (4w_k(1-w_k) - 8\delta^2)} \right)+ \\[10pt]
 +& \frac{(2 -4w_k(1-w_k)- 8\delta + 16\delta^2)}{|\log\varepsilon|} \left( -\log |D_\varepsilon| + |D_\varepsilon| \frac{\delta^2 - \delta_\varepsilon}{\varepsilon (2 -4w_k(1-w_k)- 8\delta + 16\delta^2)} \right) + o(1).
\end{aligned}
\end{equation}\\
Now, using that if $x,M>0 $ then $\log x \leq Mx -\log M$, we have (using that $\delta_\varepsilon\leq \delta^2$ eventually)

\begin{equation}
\begin{aligned}
F_\varepsilon(u_\varepsilon) \ge \, & \frac{4w_k(1-w_k) - 8\delta^2}{|\log\varepsilon|} \log \left( \frac{\delta_\varepsilon}{\varepsilon (4w_k(1-w_k) - 8\delta^2)} \right)+ \\[10pt]
& + \frac{(2 -4w_k(1-w_k)- 8\delta + 16\delta^2)}{|\log\varepsilon|} \log \left( \frac{\delta^2 - \delta_\varepsilon}{\varepsilon (2 -4w_k(1-w_k)- 8\delta + 16\delta^2)} \right) + o(1)=\\[10pt]
= \, & (4w_k(1-w_k) - 8\delta^2) \left( 1 - \frac{|\log \delta_\varepsilon|}{|\log \varepsilon|} + \frac{O(1)}{|\log \varepsilon|} \right) \\
& + (2 -4w_k(1-w_k)- 8\delta + 16\delta^2) \left( 1 + \frac{O(1)}{|\log \varepsilon|} \right) + o(1) \\[10pt]
\xrightarrow{\varepsilon \to 0} \, & (4w_k(1-w_k) - 8\delta^2)(1 - \ell) + (2 -4w_k(1-w_k)- 8\delta + 16\delta^2).
\end{aligned}
\end{equation}
Taking the limit as $\delta \to 0 $ we have that
$\liminf_{\varepsilon\to 0}F_\varepsilon(u_\varepsilon)\ge 2-4\ell w_k(1-w_k) $ and therefore the thesis follows.
\end{proof}

\subsection{\texorpdfstring{Lower bounds for the $n$-dimensional case}{Lower bounds for the n-dimensional case}}\label{sub:lower_nd}
In this section, we extend the one-dimensional lower bound to $n \ge 2$ by slicing. 
Differently from the local energy where we slice with respect to the normal of the jump set $J_u$, the non-local energy is a superposition of one-dimensional $H^s$ energies along \emph{all} directions $\xi\in \mathbb{S}^{n-1}$ (see \eqref{1}). 

We prove the following result.
\begin{prop}\label{lb_mu}  Let $n\ge2$ and $s\in\left[\frac12,1\right)$.\\
{\bf Regime (R1).} Assume \hyperlink{R1}{(R1)} holds and let $v_i$ and $v_{i+1}$ be two consecutive points in $\Z\cup S$. Let \( u \in BV(Q; \{v_i,v_{i+1}\}) \). Let \( \{u_\varepsilon \} \subseteq L^1(Q; [v_i,v_{i+1}]) \) be a sequence such that \( u_\varepsilon \) converges to \( u \) in \( L^1(Q) \). Then
\begin{equation}\displaystyle
    \liminf_{\varepsilon\to 0} F_\varepsilon^{\Z\cup S}(u_\varepsilon)\ge \Psi_{\mathrm{sup}}(v_i,v_{i+1})\mathcal{H}^{n-1}(J_u).
\end{equation}

\noindent{\bf Regime (R2).} Assume \hyperlink{R2}{(R2)} holds. Let \( u \in BV(Q; \{k,k+1\}) \) with $k \in \Z$. Let \( \{u_\varepsilon \} \subseteq L^1(Q; [k, k+1]) \) be a sequence such that \( u_\varepsilon \) converges to \( u \) in \( L^1(Q) \).
 Then
\begin{equation}
    \liminf_{\varepsilon\to 0} F_\varepsilon(u_\varepsilon)\ge \Psi_{\mathrm{sub}}(k,k+1)\mathcal{H}^{n-1}( J_u).
\end{equation}
\end{prop}
\subsubsection*{Notation}  Fix $n\ge2$ and $s\in\left[\frac12,1\right)$; to treat both regimes at once, we use the following notation. 
\begin{itemize}
    \item In regime \hyperlink{R1}{(R1)}: $a<b$ are two consecutive points of $\Z\cup S,$ $V_\varepsilon:=W^S=\operatorname{dist}^2(\cdot,\Z\cup S)$ (independent of $\varepsilon$), $\bar F_\varepsilon:=F_\varepsilon^{\Z\cup S}$ as in \eqref{F0} and $\Psi^{(1)}:=\lambda_{s,n}^{-1}\Psi_{\mathrm{sup}}$ the one-dimensional cost.
     \item In regime \hyperlink{R2}{(R2)}: $a<b$ are two consecutive points of $\Z,$ $V_\varepsilon:=W_\varepsilon $, $\bar F_\varepsilon:=F_\varepsilon$ as in \eqref{eq:main_functional}, and  $\Psi^{(1)}:=\lambda_{s,n}^{-1}\Psi_{\mathrm{sub}}$ the one-dimensional cost.
\end{itemize}
With this notation, we have
$$
\bar F_\varepsilon(v;Q):=\frac{1}{\theta_\varepsilon}\left(
\varepsilon^{2s}[v]^2_{H^s(Q)}+\int_Q V_\varepsilon(v(t))\,dt
\right),\quad\text{for all } v \in L^1(Q).
$$
On each slice, the constants in front of the seminorm and of the potential are no longer 1. We first need to establish how the multiplicative constants in front of the two terms affect the cost of a unit jump.
Let $I\subset \R$ be a bounded open interval, $V:\R\rightarrow [0,+\infty)$ a continuous potential and $c>0.$ We define 
\begin{equation}\label{eq:funzionale_slicing}
    \mathcal{G}_\varepsilon(v;I,V,c):=\frac{1}{\theta_\varepsilon}\left(
\varepsilon^{2s}c[v]^2_{H^s(I)}+\int_I V(v(t))\,dt
\right),\quad\text{for all } v \in L^1(I).
\end{equation}
The following lemma extends \Cref{lemma:lb_1d_critico} and \Cref{lemma:lb_1d_sub} to arbitrary constants in front of the two terms.
\begin{lemma}\label{lemma:slicing}  Let $n\ge2$ and $s\in\left[\frac12,1\right)$, and assume either \hyperlink{R1}{(R1)} or \hyperlink{R2}{(R2)}. Let $a,b,V_\varepsilon,\Psi^{(1)}$ be as above. Let $I \subset \R$ be a bounded open interval and $u \in BV(I;\{a,b\})$. Let $\{u_\varepsilon\}\subset L^1(I;[a,b])$ be a sequence that converges to $u$ in $L^1(I)$. Then for all $c,p>0 $ we have
    \begin{equation}
        \liminf_{\varepsilon\to 0}\mathcal{G}_\varepsilon(u_\varepsilon;I,pV_\varepsilon,c)\ge p^{1-\frac{1}{2s}}\,c^\frac{1}{2s}\,\Psi^{(1)}(a,b)\,\mathcal{H}^0(J_u).
    \end{equation}
    For $c=p=1$ this reduces to \Cref{lemma:lb_1d_critico,lemma:lb_1d_sub}.
\end{lemma}
\begin{proof}
Since $\mathcal G_\varepsilon(v;I,pV_\varepsilon,c)=c\,\mathcal G_\varepsilon(v;I,qV_\varepsilon,1)$ with $q:=p/c$,
and $c\,q^{1-\frac{1}{2s}}=p^{1-\frac{1}{2s}}c^{\frac1{2s}}$, we may take $c=1$.

Setting 
$\bar \varepsilon:=q^{-\frac{1}{2s}}\varepsilon$, we have $\varepsilon^{2s}=q\,\bar\varepsilon^{2s}$, hence
\[
\mathcal G_\varepsilon(v;I,qV_\varepsilon,1)=\frac{q\,\theta_{\bar\varepsilon}}{\theta_\varepsilon}\,
\frac1{\theta_{\bar\varepsilon}}\Big(\bar\varepsilon^{2s}[v]^2_{H^s(I)}+\int_IV_\varepsilon(v) \Big),
\]
$$\text{where}\quad 
V_{\varepsilon}=\begin{cases}
    W^S,\quad &\text{in regime \hyperlink{R1}{(R1)}},\\
    W_{ \varepsilon},\quad &\text{in regime \hyperlink{R2}{(R2)}},
\end{cases}
$$
$$\text{with}
\qquad \frac{q\,\theta_{\bar \varepsilon}}{\theta_\varepsilon}= q^{1-\frac{1}{2s}}\quad \text{for }s>\frac12 \quad \text{and}\quad 
\frac{q\,\theta_{\bar \varepsilon}}{\theta_\varepsilon}\to q^{1-\frac{1}{2s}}\quad \text{for }s=\frac12.
$$ 
Since, in regime \hyperlink{R2}{(R2)},
$$
\lim_{\varepsilon\to 0}\frac{\delta_{\varepsilon}}{\theta_{\bar \varepsilon}}=\lim_{\varepsilon\to 0}\frac{\delta_{\varepsilon}}{\theta_{ \varepsilon}}=+\infty,
$$
we can apply \Cref{lemma:lb_1d_critico,lemma:lb_1d_sub} along $\bar\varepsilon\to0$ and obtain 
$$
\liminf_{\varepsilon\to 0}\mathcal G_\varepsilon(u_\varepsilon;I,pV_\varepsilon,c)=c\liminf_{\varepsilon\to 0}\mathcal G_\varepsilon(u_\varepsilon;I,qV_\varepsilon,1)\ge cq^{1-\frac{1}{2s}}\Psi^{(1)}(a,b)\mathcal{H}^0(J_u).
$$

\end{proof}
In the following, we slice $\bar F_\varepsilon$ in all directions.
\begin{lemma}\label{lemma:slicing_F}
    Let $n\ge2$, $s\in\left[\frac12,1\right)$ and $Q\subset \R^n$ an open cube, and assume either \hyperlink{R1}{(R1)} or \hyperlink{R2}{(R2)}.
    Let $a,b,V_\varepsilon,\Psi^{(1)}$ be as above. Let $p $ be a probability density on $\mathbb{S}^{n-1}$ with $p>0 $ $\mathcal{H}^{n-1}$-a.e. If $\{u_\varepsilon\}\subset L^1(Q;[a,b])$ converges in $L^1(Q)$ to $u \in BV(Q;\{a,b\})$, then
    \begin{equation}\label{eq:slicing_1}
        \liminf_{\varepsilon\to 0}\bar F_\varepsilon(u_\varepsilon;Q)\ge 2^{-\frac{1}{2s}}\,\Psi^{(1)}(a,b)\,\int_{J_u\cap Q}d\mathcal{H}^{n-1}\int_{\mathbb{S}^{n-1}}d\mathcal{H}^{n-1}(\xi )\,p(\xi)^{1-\frac{1}{2s}}\,|\nu_u\cdot \xi |,
    \end{equation}
    where $\nu_u$ is the unit normal to the jump set $J_u$ of $u.$
\end{lemma}

\begin{proof}
    By the slicing formula for the $H^s$-seminorm (see \cite[Theorem 6.47]{Leoni}), we have
    \begin{equation}\label{1}
        [u]^2_{H^s(Q)}=\frac{1}{2}\int_{\mathbb{S}^{n-1}}d\mathcal{H}^{n-1}(\xi)\int_{\xi^\perp}d\mathcal{H}^{n-1}(y)\int_{Q_{\xi,y}}\int_{Q_{\xi,y}} \frac{|u(y+t\xi)-u(y+\tau \xi )|^2}{|t-\tau|^{1+2s}}\, dt d\tau 
    \end{equation}
    where $Q_{\xi,y}:=\{t \in \R\,:\,y+t\xi \in Q\}$. For the potential term, by Tonelli's theorem for every $\xi \in \mathbb{S}^{n-1}$ and every measurable $V\ge0$ we have
    $$
    \int_Q V(u(x)) \, dx=\int_{\xi^\perp}d\mathcal{H}^{n-1}(y)\int_{Q_{\xi,y}} V(u(y+t\xi))\,dt.
    $$
    Now, since $p$ is a probability density on $\mathbb{S}^{n-1},$ we have
    \begin{equation}\label{2}
    \int_Q V(u(x)) \, dx=\int_{\mathbb{S}^{n-1}}d\mathcal{H}^{n-1}(\xi)\int_{\xi^\perp}d\mathcal{H}^{n-1}(y)\,\left(p(\xi)\int_{Q_{\xi,y}} \, V(u(y+t\xi))\,dt\right).
    \end{equation}
    Using \eqref{1} and \eqref{2}, we can rewrite the functional as follows
    \begin{equation}
        \bar F_\varepsilon(u_\varepsilon)=\int_{\mathbb{S}^{n-1}}d\mathcal{H}^{n-1}(\xi)\int_{\xi^\perp}d\mathcal{H}^{n-1}(y)\,\mathcal{G}_\varepsilon(u_\varepsilon(y +\cdot\xi );Q_{\xi,y},p(\xi)V_\varepsilon,1/2),
    \end{equation}
    where the integrand vanishes if $Q_{\xi,y}=\emptyset.$ 

    Now, observe that for a.e. $(\xi, y) $, we have
    $$
    (u_\varepsilon)_{\xi,y}:=u_\varepsilon(y+ \cdot\xi)\rightarrow u(y+\cdot\xi )=:u_{\xi,y}\in BV(Q_{\xi,y};\{a,b\})\quad \text{in }L^1(Q_{\xi,y}),
    $$
    as a consequence of Tonelli's theorem. Hence \Cref{lemma:slicing} gives, for a.e. $(\xi,y)$,
    $$
    \liminf_{\varepsilon\to 0 }\mathcal{G}_\varepsilon((u_\varepsilon)_{\xi,y};\,Q_{\xi,y},\,p(\xi)V_\varepsilon,\,1/2)\ge 2^{-\frac{1}{2s}}p(\xi)^{1-\frac{1}{2s}}\Psi^{(1)}(a,b) \mathcal{H}^0(J_{u_{\xi,y}}).
    $$
    By \cite[Theorem 3.103]{Ambrosio}, we have
    $$
    \int_{\xi^\perp}\mathcal{H}^0(J_{u(y+\cdot\xi)})\,\,
    d\mathcal{H}^{n-1}(y)
    =(b-a)^{-1}\sup\left\{\int_Q u \frac{\partial \phi}{\partial \xi}\,:\, \phi\in C^\infty_c(Q),\,\,\,\|\phi\|_\infty\le 1 \right\}=\int_{J_u} |\nu_u\cdot \xi|\,d\mathcal{H}^{n-1}
    $$
    where $\nu_u(x)$ is the normal to the jump set $J_u$  at $\mathcal{H}^{n-1}$-a.e. point $x\in J_u$.
    Then, by Fatou's lemma, the thesis follows.
\end{proof}
The best choice of $p$ is given in the following result.
\begin{lemma}\label{probability_choice}
    Let $e\in \mathbb{S}^{n-1}$  and let $p $ be a probability density on $\mathbb{S}^{n-1}$ with $p>0 $ $\mathcal{H}^{n-1}$-a.e. Then 
    \begin{equation}\label{eq:sup_probabilita}
    \int_{\mathbb{S}^{n-1}}\,p(\xi)^{1-\frac{1}{2s}}\,|e\cdot \xi |\,d\mathcal{H}^{n-1}(\xi )\le (2C_{s,n})^\frac{1}{2s}.
     \end{equation}
    Moreover, if
  \begin{equation}\label{probabilita_giusta}
    p(\xi)=
       \frac{ |e\cdot \xi|^{2s}}{2\,C_{s,n}},
    \end{equation}
    then equality holds in \eqref{eq:sup_probabilita}.
\end{lemma}
\begin{proof}
    We first prove that, for all $e \in \mathbb{S}^{n-1}$,
    \begin{equation}\label{dim:probabilita}
    \int_{\mathbb{S}^{n-1}}|e\cdot \xi|^{2s}d\mathcal{H}^{n-1}(\xi )=2\,C_{s,n}.
        \end{equation}
    Let us decompose  $\R^n=e^\perp\oplus \R e \simeq\R^{n-1}\times \R$;  then for each $\xi \in \mathbb{S}^{n-1}$ we can write
    $$
\xi=(\xi',\xi_n)\in\mathbb{R}^{n-1}\times\mathbb{R}\quad \text{ with }
\qquad
\xi_n=\xi\cdot e\text{ and } \xi'=\xi-\xi_ne.
$$
If $\xi_n>0$ (analogously if $\xi_n<0$), we use the change of variables
$$
\zeta:=\frac{\xi'}{\xi_n}\in \R^{n-1}.
$$
Then $$\xi =\frac{(\zeta,1)}{\sqrt{1+|\zeta|^2}},\qquad 
d\mathcal{H}^{n-1}=d\zeta \frac{1}{(1+|\zeta|^2)^\frac{n}{2}}\quad \text{and}\quad |\xi \cdot e|^{2s}=\frac{1}{(1+|\zeta|^2)^s}.
$$
In conclusion, rewriting the integral \eqref{dim:probabilita} and using the definition of $C_{s,n}$ in \eqref{dhef:C_ns},  we have
$$
 \int_{\mathbb{S}^{n-1}}|e\cdot \xi|^{2s}d\mathcal{H}^{n-1}(\xi )=2 \int_{\mathbb{S}^{n-1}\cap \{\xi_n>0\}}|e\cdot \xi|^{2s}d\mathcal{H}^{n-1}(\xi )=2\int_{\R^{n-1}} \frac{d\zeta}{(1+|\zeta|^2)^\frac{n+2s}{2}}=2C_{s,n}.
$$
Now, observe that \eqref{eq:sup_probabilita} follows by H\"older's inequality with exponents $2s$ and $2s/(2s-1)$, indeed
$$
\int_{\mathbb{S}^{n-1}}\,p(\xi)^{1-\frac{1}{2s}}\,|e\cdot \xi |\le \left(\int_{\mathbb{S}^{n-1}}\,p(\xi) \right)^{1-\frac{1}{2s}}
\left(\int_{\mathbb{S}^{n-1}}\,|e\cdot \xi |^{2s}\right)^\frac{1}{2s}=(2C_{s,n})^\frac{1}{2s}.
$$
\end{proof}

\begin{proof}[Proof of \Cref{lb_mu}]
    Using the unified notation defined above we can treat both cases at once. Then, let us assume $a,b,V_\varepsilon,\Psi^{(1)}$ as before.
    We want to prove that: if $u\in BV(Q;\{a,b\})$ and $\{u_\varepsilon\}\subset L^1(Q;[a,b])$ converges to $u$ in $L^1(Q)$, then
\[
\liminf_{\varepsilon\to0}\bar F_\varepsilon(u_\varepsilon)\ \ge\ \lambda_{s,n}\Psi^{(1)}(a,b)\,\mathcal H^{n-1}(J_u).
\]
We define the sequence of Radon measures $\{\mu_\varepsilon\}$ on the cube $Q$ as follows
$$
\mu_\varepsilon(A):=\frac{1}{\theta_\varepsilon}\left(\varepsilon^{2s}\int_Q\int_A \frac{|u_\varepsilon(x)-u_\varepsilon(y)|^2}{|x-y|^{n+2s}}\,dx\,dy+\int_A V_\varepsilon(u_\varepsilon)\right),\qquad \text{for all measurable } A\subset Q,
$$
with $\mu_\varepsilon(Q)=\bar F_\varepsilon(u_\varepsilon).$
By the Besicovitch differentiation theorem, to prove the thesis, it is enough to prove that for $\mathcal{H}^{n-1}$-a.e. $y \in J_u$, we have
\begin{equation}\label{blowup}
    \liminf_{r \to 0}\,\liminf_{\varepsilon\to 0}\frac{\mu_\varepsilon(Q_r)}{r^{n-1}}\ge\lambda_{s,n}\Psi^{(1)}(a,b) ,
\end{equation}
where, with a little abuse of notation, $Q_r$ denotes the cube centered at $y$, with side $r$, parallel to the normal $\nu_{u}(y)$ on $J_u$ at $y$. Observe that $\bar F_\varepsilon(u_\varepsilon;Q_r)\le \mu_\varepsilon(Q_r)\le \mu_\varepsilon(Q)\le F_\varepsilon(u_\varepsilon)\le \sup_\varepsilon F_\varepsilon(u_\varepsilon)$, where the last term can be assumed finite. Hence, up to a subsequence $\mu_\varepsilon\rightharpoonup^*\mu$ and 
$$
\mu(\overline {Q_r})\ge \liminf_{\varepsilon\to 0}\mu_\varepsilon(Q_r)\ge \liminf_{\varepsilon\to 0}\bar F_\varepsilon(u_\varepsilon;Q_r).
$$

Now, fixing $y \in J_u$, 
by \Cref{lemma:slicing_F}, for each $p\ge 0$ with $\int_{\mathbb{S}^{n-1}}p(\xi)\,d\xi=1$ and $p>0$ $\mathcal{H}^{n-1}$-a.e., we have
\begin{equation}
\liminf_{\varepsilon\to 0}\bar F_\varepsilon(u_\varepsilon;Q_r)\ge 2^{-\frac{1}{2s}}\,\Psi^{(1)}(a,b)\,\int_{J_u\cap Q_r}d\mathcal{H}^{n-1}(x)\int_{\mathbb{S}^{n-1}}d\mathcal{H}^{n-1}(\xi )\,p(\xi)^{1-\frac{1}{2s}}\,|\nu_u(x)\cdot \xi |;
\end{equation}
in particular, we can choose $p(\xi)=|\nu_u(y)\cdot \xi|^{2s}(2C_{s,n})^{-1}$. Hence, using Fatou's lemma
\begin{multline*}
\liminf_{r \to 0}\,\liminf_{\varepsilon \to 0}\,
\frac{\bar{F}_\varepsilon(u_\varepsilon; Q_r)}{r^{n-1}}
\;\ge\;
\Psi^{(1)}(a,b)
\int_{\mathbb{S}^{n-1}}
\frac{|\nu_u(y)\cdot\xi|^{2s-1}}{2C_{s,n}^{\,1-\frac{1}{2s}}} \\
\cdot
\left(
\liminf_{r \to 0}\,\frac{1}{r^{n-1}}
\int_{J_u \cap Q_r} |\nu_u(x)\cdot\xi|\,
d\mathcal{H}^{n-1}(x)
\right)
d\mathcal{H}^{n-1}(\xi).
\end{multline*}
By rectifiability of $J_u$, we have
$$
\lim_{r\to 0}\frac{1}{r^{n-1}}\int_{J_u\cap Q_r}d\mathcal{H}^{n-1}(x )\,|\nu_u(x)\cdot \xi |= |\nu_u(y)\cdot \xi |.
$$
Hence, \eqref{blowup} follows by \Cref{probability_choice}, recalling that $\lambda_{s,n}=(C_{s,n})^\frac{1}{2s} $ (see \eqref{dhef:C_ns}).
\end{proof}

\subsection{Proof of the lower bound}\label{sub:proof_lower}
\begin{proof}[Proof of \Cref{thm:main} (lower bound)]
The lower bound is an immediate consequence of \Cref{lb_mu}, the truncation inequality \eqref{eq:truncation} and the lower estimate in \eqref{eq:two_sided_estimate}. We give details for regime \hyperlink{R1}{(R1)}; \hyperlink{R2}{(R2)} is similar.

    \noindent\textbf{Regime \hyperlink{R1}{(R1)}}.
    Let $\mathcal{V}=\union=\{v_i\}_{i\in\Z}$ with $v_i<v_{i+1}$.
    Let $\{u_\varepsilon\}\subset L^1(Q)$ be a sequence that converges to $u$ in $L^1(Q)$ such that $F_\varepsilon(u_\varepsilon)\leq M$. Then, by \Cref{thm:compactness}, $u\in BV(Q;\mathcal{V}).$ 
    Let $T_i$ be the truncations between consecutive points of $\mathcal{V}$ in \eqref{eq:truncation}.
By \eqref{eq:two_sided_estimate} and the truncation inequality \eqref{eq:truncation}, for any $\delta \in (0,1)$ and small $\varepsilon$, we have that
\begin{equation}\label{0}
    F_\varepsilon(u_\varepsilon)\ge F^{\Z\cup S}_\varepsilon(u_\varepsilon) + \frac{\delta_\varepsilon}{\theta_\varepsilon} |A_{\delta}(u_\varepsilon)| \ge \sum_{i\in \Z}
    F^{\Z\cup S}_\varepsilon(T_{i}u_\varepsilon) + \frac{\delta_\varepsilon}{\theta_\varepsilon} |A_{\delta}(u_\varepsilon)|.
\end{equation}
Now, for each $i \in \Z$, $\{T_{i}u_\varepsilon\}$ converges in $L^1$ to $T_{i}u$ and satisfies the hypotheses of \Cref{lb_mu} (if $n\geq 2$) and \Cref{lemma:lb_1d_critico} (if $n=1$); hence, we have
\begin{equation}
\liminf_{\varepsilon \to 0} F^{\Z\cup S}_\varepsilon(T_{i} u_\varepsilon) 
\geq \Psi_{\mathrm{sup}}(v_i,v_{i+1})\, \mathcal{H}^{n-1}(J_{T_{i} u})=\int_{J_{T_{i} u}}\Psi_{\mathrm{sup}}(v_i,v_{i+1})\, d\mathcal{H}^{n-1}(x)
\end{equation}
Now observe that, up to an $ \mathcal{H}^{n-1}\text{-null } $ set,
\[
J_u = \bigcup_{j \in \mathbb{Z}} J_{T_{j} u}.
\]
We can use additivity of $\Psi_{\mathrm{sup}} $ to conclude that 
\begin{equation}
\begin{split}
\sum_{i \in \mathbb{Z}} \Psi_{\mathrm{sup}}(v_i, v_{i+1})\, \chi_{J_{T_{i} u}}(x)
&= \sum_{i \,:\, [v_i, v_{i+1}] \subset [u^-(x), u^+(x)]} \Psi_{\mathrm{sup}}(v_i, v_{i+1})\, \chi_{J_u}(x) \\
&= \Psi_{\mathrm{sup}}(u^-(x), u^+(x))\, \chi_{J_u}(x),
\end{split}
\end{equation}
hence
   \begin{equation}\label{1_1}
\begin{aligned}
\liminf_{\varepsilon \to 0} F^{\Z\cup S}_\varepsilon(u_\varepsilon) 
\ge \int_{J_u} \Psi_{\mathrm{sup}}(u^+(x), u^-(x)) \, d\mathcal{H}^{n-1}(x).
\end{aligned}
\end{equation}
    Moreover, by the convergence in measure, $\liminf_\varepsilon |A_{\delta}(u_\varepsilon)|\ge |\{u\in S\}|$, so that 
    \begin{equation}\label{2_1}
    \liminf_{\varepsilon\to 0 }\frac{\delta_\varepsilon}{\theta_\varepsilon} |A_{\delta}(u_\varepsilon)|\ge K|\{u\in S\}|.
    \end{equation}
    By \eqref{0}, using \eqref{1_1} and \eqref{2_1}, we conclude that
    $$
    \liminf_{\varepsilon \to 0}F_\varepsilon(u_\varepsilon)\ge
    F_{\mathrm{sup}}(u).
    $$
\\\noindent
\textbf{Regime \hyperlink{R2}{(R2)}}.
We repeat the truncation argument with $\mathcal{V}=\Z$,
  for a sequence $\{u_\varepsilon\}\subset L^1(Q)$ converging to $u$ in $L^1(Q)$, and such that $F_\varepsilon(u_\varepsilon)\leq M$. By \Cref{thm:compactness} we have that $u\in BV(Q;\mathcal{V})$. 
By \Cref{lb_mu} if $n\ge 2$ and \Cref{lemma:lb_1d_sub} if $n=1$, we obtain
  $$
    \liminf_{\varepsilon \to 0}F_\varepsilon(u_\varepsilon)\ge 
    \int_{J_u} \Psi_{\mathrm{sub}}(u^+(x), u^-(x)) \, d\mathcal{H}^{n-1}(x)=F_{\mathrm{sub}}(u).
    $$
\end{proof}

\section{Upper bound}\label{sec:upper_bound}
The proof of the upper bound is organized according to the fractional parameter $s$ and the regime under consideration: $s=\frac12$ in \Cref{up_1/2} and $s>\frac12$ in \Cref{up_1}.
In each case, we first construct a one-dimensional recovery sequence,
and then extend the construction to higher dimensions
via a flat-interface approximation argument.

We denote by $\mathcal{V}$ the set of wells of the limit, i.e. $\mathcal{V}=\Z \cup S$ in regime \hyperlink{R1}{(R1)} and $\mathcal{V}=\Z $ in regime \hyperlink{R2}{(R2)}, ordered as $\{v_i\}_{i \in \Z}$ with $v_i<v_{i+1},$ and by $F$ the candidate $\Gamma$-limit, i.e. $F_{\rm sup}$ in regime \hyperlink{R1}{(R1)} and $F_{\rm sub}$ in regime \hyperlink{R2}{(R2)}.
By standard properties of $\Gamma$-convergence, it suffices to construct, for every $u $ in a class $ \mathcal{D} \subseteq BV(Q;\mathcal{V})$, dense in energy, a recovery sequence $\{u_\varepsilon\}$ satisfying
\[
u_\varepsilon \to u \quad \text{in } L^1(Q) 
\quad \text{and} \quad
\limsup_{\varepsilon \to 0} F_\varepsilon(u_\varepsilon) \le F(u).
\]

Let $\mathcal{D}$ be the set of functions $u \in L^\infty(Q; \mathcal{V}) \cap BV(Q; \mathcal{V})$ whose jump set $J_u$ is a finite union of $(n-1)$-dimensional simplices, 
\(
J_u = \bigcup_{i=1}^m J^i
\), with $\operatorname{dim}_{\mathcal{H}}(J^i\cap J^j)\le n-2$
 for $i \neq j$, and such that, for $\mathcal{H}^{n-1}$-almost every $x \in J_u$, the trace values $u^+(x)$ and $u^-(x)$ are \textit{consecutive} elements of $\mathcal{V}$: 
\begin{equation}
    \label{eq-split}
    \{u^+(x), u^-(x)\} = \{v_j, v_{j+1}\} \quad \text{for some } j \in \mathbb{Z}.
\end{equation}
In regime \hyperlink{R2}{(R2)} this means $|u^+-u^-|(x)= 1$; in regime \hyperlink{R1}{(R1)}, the jump amplitude is either $s_k - k$ or $(k+1) - s_k$.
The density of the class \(\mathcal{D}\) is a consequence of \cite{Bra-Con-Gar}, which provides an approximation in energy of any $u\in BV(Q; \mathcal{V})$ by functions whose jump sets are finite unions of \((n-1)\)-dimensional simplices. The property \eqref{eq-split} can then be obtained by splitting each flat interface into parallel ones which decompose any large jump into jumps between consecutive wells.

In what follows $u\in \mathcal{D}$. To avoid boundary effects, we  extend $u$ to a larger cube $ \tilde Q$ containing $Q$,  reflecting $u$ across the faces of $Q$; the extension, still denoted by $u$, has no jumps across $\partial Q$. Let $Q\subset \subset Q'\subset\subset \tilde Q.$ 
In the following, $\varphi_\varepsilon$ is the standard mollifier supported in $B_\varepsilon(0).$
\paragraph{\emph{Higher dimension notation}.}
Since the jump set of $u$ is a union of simplices, namely $J_u=\bigcup_{i=1}^m J^i$, for each $i$ we denote by $u^i_-<u^i_+$ the traces of $u $ on $J^i$ and by $\nu^i$ the unit normal to $J^i$ pointing towards $\{u=u^i_+\};$ in regime \hyperlink{R2}{(R2)} we define also the intermediate well $s^i:=s_{u_-^i}$ and $w^i:=s^i-u^i_-.$
Let $$
\Sigma:=\bigcup_{i\neq j}(J^i \cap J^j),\quad  \mathcal{N}_\sigma := \{ x \in \tilde Q : \, \operatorname{dist}(x, J_u) < \sigma \}.
$$
Since $\Sigma$ is contained in a finite union of $(n-2)$-dimensional simplices, for $\sigma$ small
\begin{equation}\label{stima_N}
|\mathcal{N}_\sigma|\le C\sigma,\quad |\{x\in \tilde Q\,:\,\operatorname{dist}(x,\Sigma)<\sigma\}|\le C\sigma^2,
\end{equation}
where $C$ depends on $u.$ 
\subsection{\texorpdfstring{Upper bound for $s=\frac12$}{Upper bound for s=1/2}}\label{up_1/2}
For $s = \frac12$ the energy is dominated by the logarithmic divergence of the $H^{1/2}$ seminorm. In regime \hyperlink{R1}{(R1)} the recovery sequence is a transition between consecutive wells at scale $\varepsilon$.
In regime \hyperlink{R2}{(R2)}, in addition, it rests at the intermediate well on a layer of width about $\varepsilon/\delta_\varepsilon$.

\subsubsection{Dimension one}\label{subsub:ub_1_1}
\begin{proof}[Proof of \Cref{thm:main}]\textit{(Upper bound, \( n=1 ,\) $s=\frac12$, \hyperlink{R2}{(R2)})} We extend $u$ to $\R$ with a constant without creating a jump point at $\partial Q$ (in the two extremal points of the interval).

 Let $\rho:=\min\{|x-y|\,:\,x,y\in J_u, \,x\neq y\}$ ($\rho:=1$ if $J_u=\{x\}$).
Let \( \varepsilon_0 \) be such that \( \delta_\varepsilon < 3/4 \) and $4\varepsilon/\delta_\varepsilon<\rho$ for  \( \varepsilon < \varepsilon_0\) (this is possible by \eqref{eq-sub-regime}). For $x_0\in J_u$, let $u^-(x_0)<u^+(x_0)=u^-(x_0)+1$ be the two traces of $u$ at $x_0$, let $s(x_0):=s_{u^-(x_0)}$ be the intermediate well  and $w(x_0)=s(x_0)-u^-(x_0)$.
We define
\begin{equation}\label{recovery_1_R2}
\tilde u_\varepsilon:=\displaystyle\begin{cases} 
\displaystyle s(x_0), & \text{in } (x_0-\tfrac{\varepsilon}{\delta_\varepsilon},x_0+\tfrac{\varepsilon}{\delta_\varepsilon}), \,x_0\in J_u,\\
u, & \text{otherwise },
\end{cases}
\quad  \text{and}\quad u_\varepsilon:= \tilde u_\varepsilon \,*\, \varphi_\varepsilon.
\end{equation}

By definition,   
\(
\norm{u - u_\varepsilon}_\infty \leq 1, \,\text{and} \,
|\{u \neq u_\varepsilon\}| \leq 4 \mathcal{H}^0(J_u) \varepsilon/\delta_\varepsilon, 
\)
so \( u_\varepsilon \to u \) in \( L^1(Q) \).

Now we estimate the \( H^{1/2} \)-seminorm. First of all, we notice that we can ignore the interaction between distant points, i.e., 
$$
\begin{aligned}
[u_\varepsilon]^2_{H^{1/2}(Q)} &= O(1) + \sum_{x_0 \in J_u} [u_\varepsilon]^2_{H^{1/2}(B_{\rho/2}(x_0))}.
\end{aligned}
$$
Fixing $x_0\in J_u$, in $B_{\rho/2}(x_0)$ let
\(
A_\varepsilon := \{u_\varepsilon = u^-(x_0)\}, \, P_\varepsilon := \{u_\varepsilon = s(x_0)\} \text{ and } C_\varepsilon := \{u_\varepsilon = u^+(x_0)\}.
\)
Then, using that $u_\varepsilon$ is $C/\varepsilon$-Lipschitz, the contribution of the remaining interactions is $O(1).$ Hence, we have
$$
\begin{aligned}
[u_\varepsilon]^2_{H^{\frac{1}{2}}(B_{\rho/2}(x_0))} &= O(1) + 2 \int_{A_\varepsilon} \int_{P_\varepsilon} \frac{|u_\varepsilon(x) - u_\varepsilon(y)|^2}{|x - y|^2} \, dx \, dy \\
& + 2 \int_{A_\varepsilon} \int_{C_\varepsilon} \frac{|u_\varepsilon(x) - u_\varepsilon(y)|^2}{|x - y|^2} \, dx \, dy + 2 \int_{C_\varepsilon} \int_{P_\varepsilon} \frac{|u_\varepsilon(x) - u_\varepsilon(y)|^2}{|x - y|^2} \, dx \, dy \\
&= O(1) + 2w(x_0)^2 \int_{A_\varepsilon} \int_{P_\varepsilon} \frac{ dx \, dy}{|x - y|^2}  \\
&\quad + 2 \int_{A_\varepsilon} \int_{C_\varepsilon} \frac{1}{|x - y|^2} \, dx \, dy  + 2(1-w(x_0))^2 \int_{C_\varepsilon} \int_{P_\varepsilon} \frac{ dx \, dy}{|x - y|^2}.
\end{aligned}
$$
By \eqref{conto_espl}, since $A_\varepsilon,$ $P_\varepsilon,$ $C_\varepsilon$ are intervals separated by gaps of lengths at least $\varepsilon$ and $P_\varepsilon$ has length proportional to $\varepsilon/\delta_\varepsilon-\varepsilon,$ we have
$$
\int_{C_\varepsilon} \int_{P_\varepsilon} \frac{ dx \, dy}{|x - y|^2}  = \int_{A_\varepsilon} \int_{P_\varepsilon} \frac{ dx \, dy}{|x - y|^2}  = |\log  \delta_\varepsilon|+O(1)\quad 
\text{and} 
\quad
\int_{C_\varepsilon} \int_{A_\varepsilon} \frac{ dx \, dy}{|x - y|^2}  = \log \frac{\delta_\varepsilon}{\varepsilon}+O(1).
$$
For the potential, we have
$$
\begin{aligned}
\int_Q W_\varepsilon(u_\varepsilon(x)) \, dx &= \sum_{y \in J_u} \int_{B_{\rho/2}(y)} W_\varepsilon(u_\varepsilon(x)) \, dx+C\varepsilon \le \mathcal{H}^0(J_u) \int_{B_{\rho/2}(x_0)} W_\varepsilon(u_\varepsilon(x)) \, dx+C\varepsilon,
\end{aligned}
$$
with $x_0\in J_u$. Now, since $W_\varepsilon(s(x_0))=\delta_\varepsilon$ and it is bounded, we have
$$
\int_{B_{\rho/2}(x_0)} W_\varepsilon(u_\varepsilon(x)) \, dx \leq |P_\varepsilon| \delta_\varepsilon + \varepsilon\sup_\varepsilon W_\varepsilon \leq C\varepsilon.
$$
In conclusion, since $\theta_\varepsilon=\varepsilon|\log \varepsilon|$ and $|\log \delta_\varepsilon|/|\log \varepsilon|\to \ell,$ we have
$$
\begin{aligned}
\limsup_{\varepsilon \to 0} F_\varepsilon(u_\varepsilon) 
&\leq  \limsup_{\varepsilon \to 0} 
\frac{1}{|\log \varepsilon|} \sum_{x_0 \in J_u}
\Bigl[ 
    2 |\log \varepsilon| - \bigl( 2 - 2 w(x_0)^2 - 2 (1-w(x_0))^2 \bigr) |\log \delta_\varepsilon| 
\Bigr] \\&=
\sum_{x_0 \in J_u} 
\Bigl( 2 - 4w(x_0)(1-w(x_0))\ell \Bigr)=F_{\mathrm{sub}}(u).
\end{aligned}
$$

\end{proof}

\begin{ohss}\label{shape}
The specific shape of \( u_\varepsilon \) near the jump is not relevant for the upper bound, as already observed in  
\cite{ABS1} and \cite{Kurzke}. 
Indeed, in the previous computation we only used that, near each jump point $x_0$, $u_\varepsilon$ is $C/\varepsilon$-Lipschitz, that $|\{u_\varepsilon=s(x_0)\}|\sim 
\varepsilon/\delta_\varepsilon$ and that $u_\varepsilon$ takes values outside $\{u^-(x_0),s(x_0),u^+(x_0)\}$ on a set of measure $O(\varepsilon).$
\end{ohss}
\begin{proof}[Proof of \Cref{thm:main}]\textit{(Upper bound, \( n=1 ,\) $s=\frac12$, \hyperlink{R1}{(R1)})}
We extend $u$ to $\R$ with a constant without creating a jump point at $\partial Q$ (in the two extremal points of the interval).
Let 
$ u_\varepsilon:= u \,*\, \varphi_\varepsilon$. Then $\{u_\varepsilon\neq u\}\subset \{x\in Q:\operatorname{dist}(x,J_u)<\varepsilon\}$, so $u_\varepsilon\to u$ in $L^1(Q)$ and $|A_1(u_\varepsilon)|\to |\{u\in S\}|.$ Hence, by the upper estimate \eqref{eq:two_sided_estimate},
$$
\limsup_{\varepsilon\to 0}F_\varepsilon(u_\varepsilon)\le \limsup_{\varepsilon\to 0}F^{\Z\cup S}_\varepsilon(u_\varepsilon)+K|\{u\in S\}|.
$$
Near each $x_0\in J_u$, the sequence $u_\varepsilon$  performs a transition between two consecutive wells $u^-(x_0),u^+(x_0)\in \Z\cup S$ within a layer of thickness $2\varepsilon$. The computation in the previous proof (with $P_\varepsilon=\emptyset$) gives $$[u_\varepsilon]^2_{H^{1/2}(Q)}=2\sum_{x_0}(u^+(x_0)-u^-(x_0))^2|\log\varepsilon|+O(1)\quad \text{and}\quad \int_Q \operatorname{dist}^2(u_\varepsilon,\Z\cup S)\le C\varepsilon.$$ Hence
$$
\limsup_{\varepsilon\to 0}F_\varepsilon(u_\varepsilon)\le 2\sum_{x_0}(u^+(x_0)-u^-(x_0))^2+K|\{u\in S\}|=F_{\rm sup}(u).
$$
\end{proof}

\subsubsection{Higher dimension}
The extension to  dimension $n\ge 2$ is based on the concentration of the energy in the jump set $J_u$, that we capture by analyzing the asymptotic properties of the measures associated with the fractional seminorm. 

For $\{u_\varepsilon\}\subset L^1(Q')$, we define the measures $\mu_\varepsilon$ on $\overline{Q'} \times \overline{ Q'}$ by
\begin{equation}\label{mu_eps}
\mu_\varepsilon(A \times B) := \frac{1}{|\log \varepsilon|} \int_A \int_B \frac{|u_\varepsilon(x) - u_\varepsilon(y)|^2}{|x - y|^{n+1}} \, dx \, dy, \quad A, B \subseteq Q'\text{ measurable.}
\end{equation}

\begin{lemma}\label{19} Let \( \{u_\varepsilon \} \) be a sequence such that \(\| u_\varepsilon\|_{L^\infty(Q')}\le M \) and \(\sup_\varepsilon\mu_\varepsilon(Q'\times Q') < +\infty \) with $\mu_\varepsilon$ defined by
\eqref{mu_eps}, and let \( \mu \) be any weak* limit of \( \mu_\varepsilon \) (up to a subsequence).
Then \( \mu \) is concentrated on the diagonal 
\(
D := \{(x,x) : x \in \overline{Q'}\}
\)
and 
\(
\lambda(A) := \mu(A \times A)
\)
 defines a Radon measure on $\overline{Q'}.$

\end{lemma}
\begin{proof}
To prove that \( \mu \) is concentrated on the diagonal \( D \),  
it is enough to show that for any continuous non-negative function \( \varphi : \tilde Q \times\tilde  Q \to \mathbb{R} \) such that \( \operatorname{supp} \varphi \cap D = \emptyset \), we have
\[
\int \varphi \, d\mu = 0.
\]
Let \( \delta := \operatorname{dist}(\operatorname{supp} \varphi, D) > 0 \), so that
\(
|x - y| \geq \delta \)  for all \( (x, y) \in \operatorname{supp} \varphi.
\)
Since \(-M \leq u_\varepsilon \leq M \), we obtain
\[
\begin{aligned}
 \int \varphi(x, y) \, d\mu
&= \lim_{\varepsilon \to 0} \frac{1}{|\log \varepsilon|} 
\int_{Q'} \int_{Q'} \varphi(x, y) 
\frac{|u_\varepsilon(x) - u_\varepsilon(y)|^2}{|x - y|^{n+1}} \, dx \, dy \leq \lim_{\varepsilon \to 0} \frac{4M^2\|\varphi\|_\infty|Q'|^2}{|\log \varepsilon| \, \delta^{n+1}} = 0.
\end{aligned}
\]
Hence \( \mu \) is concentrated on \( D \).
\end{proof}

  \begin{proof} [Proof of \Cref{thm:main}]\textit{(Upper bound, \( n\geq 2 \), $s=\frac12$)}
  Let $\sigma_\varepsilon:=\varepsilon/\delta_\varepsilon$ in regime \hyperlink{R2}{(R2)} and $\sigma_\varepsilon:=\varepsilon$ in regime \hyperlink{R1}{(R1)}. 
  
\textit{Step 1. Definition of the recovery sequence.}
In regime \hyperlink{R1}{(R1)} we set $u_\varepsilon:=u \,*\,\varphi_\varepsilon$. In regime \hyperlink{R2}{(R2)}
we modify \( u \) near each unit jump on $J^i$ by inserting a layer where the function takes the intermediate value \( s^i \). In order to keep these layers pairwise disjoint and away from $\Sigma,$ we shrink each simplex.

For  \( i = 1, \dots, m \), let  \( d^i(\cdot):=\operatorname{dist}(\cdot, J^i) \) and
\[
I^{i} := \left\{ x \in \tilde Q \,:\, d^i(x) \leq \frac{\varepsilon}{\delta_\varepsilon}\text{ and } d^j(x) >\frac{4\varepsilon}{\delta_\varepsilon}, \, \forall j \neq i \right\},
\qquad
\tilde J^i_\varepsilon := \left\{ x \in J^i \,:\, x + t \nu^i  \in I^i, \forall |t| \leq \frac{\varepsilon}{\delta_\varepsilon}\right\},
\]
$$
P^i_\varepsilon = \left\{  y + t \nu^i\,:\, \, y \in \tilde J^i_\varepsilon, \, |t| \leq \frac{\varepsilon}{\delta_\varepsilon}\right\}.
$$
Then 
\(
\mathcal{H}^{n-1}(\tilde J^i_\varepsilon) \to \mathcal{H}^{n-1}(J^i)
\), when $\varepsilon\to 0$
and the sets  \( P^i_\varepsilon + B_{\varepsilon/\delta_\varepsilon}(0) \) are pairwise disjoint.
Now, we define the sequence \( \tilde{u}_\varepsilon \), as follows:
\[
\tilde{u}_\varepsilon = 
\begin{cases} 
 s^i,& \text{in }  P^i_\varepsilon, \\
u ,& \text{otherwise in }\tilde Q
\end{cases}\quad\text{and}\quad u_\varepsilon:=\tilde{u}_\varepsilon * \varphi_\varepsilon.
\]
\noindent
In both regimes, we have
 \( \|u_\varepsilon\|_\infty\le M:=\|u\|_\infty  \) and \( u =u_\varepsilon\) in the set \( \tilde Q \setminus \mathcal{N}_{4\sigma_\varepsilon} \). Then
\[
\| u_\varepsilon - u \|_{L^1(Q)} \leq 2 M \left| \mathcal{N}_{4\sigma_\varepsilon} \right| \to 0,\quad \text{as }\varepsilon\to 0.
\]
Moreover, $|\nabla u_\varepsilon|\le C/\varepsilon.$

\textit{Step 2. Potential term.}
In regime \hyperlink{R1}{(R1)}, using the upper estimate \eqref{eq:two_sided_estimate}, the fact that $u_\varepsilon=u$ outside $\mathcal{N}_{4\varepsilon}$, the boundedness of the potential, and that $|A_1(u_\varepsilon)|\to |\{u \in S\}|$ by convergence in measure (for the definition of $A_1$ see \eqref{eq-set-u-comparison}), we have
\begin{equation}\label{ub_1_R1}
\frac{1}{\varepsilon |\log \varepsilon|} \int_{Q} W_\varepsilon(u_\varepsilon(x)) \, dx \le \frac{C|\mathcal{N}_{4\varepsilon}|}{\varepsilon|\log\varepsilon|}+\frac{\delta_\varepsilon}{\varepsilon|\log\varepsilon|}|A_1(u_\varepsilon)|\rightarrow K|\{u\in S\}|.
\end{equation}
In regime \hyperlink{R2}{(R2)}, we have that $\{u_\varepsilon\notin \Z\cup S\}\subset Z_\varepsilon:=\mathcal{N}_\varepsilon\cup 
\bigcup_i \{\operatorname{dist}(x, \partial P^i_\varepsilon) \leq \varepsilon \}, 
$ with $|Z_\varepsilon|\le C \varepsilon$. 
On the other hand, $\{u_\varepsilon \in S\}\subset P_\varepsilon:=\bigcup_i P^i_\varepsilon$ with $|P_\varepsilon|\le C\varepsilon/\delta_\varepsilon.$
Hence
\begin{equation}\label{ub_1_R2}
    \frac{1}{\varepsilon|\log\varepsilon|}\int_{Q} W_\varepsilon(u_\varepsilon(x)) \, dx 
\leq \frac{C |Z_\varepsilon| + \delta_\varepsilon |P_\varepsilon|}{\varepsilon|\log\varepsilon|}
\leq \frac{C}{|\log\varepsilon|}\rightarrow 0.
\end{equation}

\noindent
\textit{Step 3. Boundedness of the energy measures.}
Let \( \{\mu_\varepsilon\} \) be the sequence of measures on \( \overline{Q'} \times  \overline{Q'} \) defined in \eqref{mu_eps}. We show that
\[
\sup_{\varepsilon < \varepsilon_0} \mu_\varepsilon(Q' \times Q') < +\infty,
\]
and therefore \( \{\mu_\varepsilon\} \) is precompact in the weak* topology of measures on \( Q' \times Q' \).
We write 
\[
\int_{Q'}\int_{Q'} \frac{|u_\varepsilon(x) - u_\varepsilon(y)|^2}{|x - y|^{n+1}} \, dx \, dy =
\int_{\R^n} \frac{1}{|z|^{n+1}} \int_{Q'\cap (Q'-z)}  \, |u_\varepsilon(x+z) - u_\varepsilon(x)|^2\,dx\,dz
\]
and split the integral into $|z|\le \varepsilon,$ $\varepsilon<|z|<\frac34\sigma_\varepsilon$ and $|z|\ge \frac34 \sigma_\varepsilon;$ for regime \hyperlink{R1}{(R1)} the middle range is empty. For $r >0$, we define 
\begin{equation}\label{D_delta}
\mathcal{M}_r := \mathcal{N}_{\varepsilon+r} \cup \bigcup_i \left\{ x \in \tilde Q \,:\, \operatorname{dist}(x, \partial P^i_\varepsilon) \leq \varepsilon + r \right\},
\end{equation}
which contains all the points that are either close to the uncovered part of \( J_u \) or to the boundary of the transition regions (in regime \hyperlink{R1}{(R1)}, $\mathcal{M}_r=\mathcal{N}_{\varepsilon+r}$). Then 
\begin{equation}\label{stima_D_delta}
|\mathcal{M}_r| \leq C(\varepsilon + r),
\end{equation}
for some constant \( C > 0 \) independent of \( \varepsilon \) and $r$. Moreover, $u_\varepsilon(x+z)=u_\varepsilon(x)$ if \( x \notin \mathcal{M}_{|z|} \), since the segment $[x,x+z]$ is in the region where $u_\varepsilon$ is constant, and $u_\varepsilon(x+z)=u_\varepsilon(x)$ if $x \notin \mathcal{N}_{4\sigma_\varepsilon+|z|}.$
\begin{itemize}
    \item For $|z|\le \varepsilon$, since   
\(
|\nabla u_\varepsilon| \leq C/\varepsilon,
\)
we have
\[
\begin{aligned}
\int_{|z| \leq \varepsilon} \frac{dz}{|z|^{n+1}} \int_{Q'\cap{(Q'-z)}} \Bigl|u_\varepsilon(x+z)-u_\varepsilon(x)\Bigr|^{2} \,dx 
&\leq C \int_{|z| \leq \varepsilon} \frac{dz}{|z|^{n+1}}\,|\mathcal{M}_\varepsilon|\, \frac{|z|^{2}}{\varepsilon^{2}}\, \leq C \int_{|z| \leq \varepsilon} \frac{dz}{|z|^{n-1}\varepsilon} \leq C.
\end{aligned}
\]
\item For $\varepsilon<|z|<\tfrac34\sigma_\varepsilon$, by \eqref{stima_D_delta} and $\|u_\varepsilon\|_\infty\le M,$ we have
\[
\begin{aligned}
\int_{\varepsilon<|z|<\tfrac34\sigma_\varepsilon} \frac{1}{|z|^{n+1}} \int_{Q'\cap (Q'-z)} \Bigl|u_\varepsilon(x+z)-u_\varepsilon(x)\Bigr|^{2} \,dx\,dz 
&\leq 4M^2 \int_{\varepsilon<|z|<\tfrac34\sigma_\varepsilon} \frac{dz}{|z|^{n+1}}\,|\mathcal{M}_{|z|}|\\
&\leq 8CM^2 \int_{\varepsilon<|z|<\tfrac34\sigma_\varepsilon} \frac{dz}{|z|^{n}} \leq C\,(|\log\delta_\varepsilon|+1).
\end{aligned}
\]
\item For $|z|\ge \tfrac34\sigma_\varepsilon$, we have
$\mathcal{N}_{4\sigma_\varepsilon+|z|}\subseteq \mathcal{N}_{18|z|}$. Hence, using again that \(u\) is bounded by \(M\), we obtain
\[
\begin{aligned}
\int_{\frac34\sigma_\varepsilon\le |z|\le \operatorname{diam}(Q')} \frac{1}{|z|^{n+1}}
\int_{Q'\cap(Q'-z)} \bigl|u_\varepsilon(x+z)-u_\varepsilon(x)\bigr|^{2}\,dx\,dz
&\le 4M^{2}\int_{\frac34\sigma_\varepsilon\le |z|\le \operatorname{diam}(Q')}
\frac{|\mathcal{N}_{18|z|}|}{|z|^{n+1}}\,dz \\
&\le C\bigl(|\log\sigma_\varepsilon|+1\bigr).
\end{aligned}
\]\noindent
\end{itemize}
Dividing everything by $|\log \varepsilon|$, we obtain
\[
\mu_{\varepsilon}(Q' \times Q') \leq C + o(1)
\]
as \( \varepsilon  \) goes to $0$. In particular, up to a subsequence  \( \mu_\varepsilon\overset{*}{\rightharpoonup}\mu\) and, by \Cref{19}, \( \mu \) is concentrated on the diagonal and we can define $\lambda$. Moreover, $\lambda$ is concentrated on $J_u \cap \overline{Q'}.$ Indeed, considering $U_\rho:=\{\operatorname{dist}(x,\Sigma)\le\rho\},$
for $r<\rho$ we have $|\mathcal{N}_r\cap U_\rho|+|\mathcal{M}_r\cap U_\rho|\le C(\varepsilon+r)\rho$ and $|\mathcal{N}_r\cap U_\rho|\le C\rho^2$ for every $r$; hence repeating the same computation as above, we have
$$\mu_\varepsilon(U_\rho \times Q')\le C\rho +o(1)$$ and
then $\lambda(\Sigma)=0$. Moreover, $u_\varepsilon$ is constant away from $J_u$, so $\lambda(\overline{Q'}\setminus J_u)=0.$

\textit{Step 4. Density of $\lambda.$}
To conclude the proof, it suffices to show that for \( \mathcal{H}^{n-1} \)-a.e. \( x_{0} \in J^i \),  
\begin{equation}\label{tesi}
\begin{aligned}
\limsup_{r\to 0}
\frac{\lambda(Q_{2r}^{\nu^i}(x_0))}{(2r)^{n-1}}
&\leq \omega_{n-1}\Psi^{(1)}(u^i_-,u^i_+),\\
\text{where}\qquad
\Psi^{(1)}(u^i_-,u^i_+)
&=
\begin{cases}
2(u^i_+-u^i_-)^2,
& \text{in \hyperlink{R1}{(R1)},}\\[2mm]
2(1-2\ell w^i(1-w^i)),
& \text{in \hyperlink{R2}{(R2)},}
\end{cases}
\end{aligned}
\end{equation} 
\noindent
where \( Q_{2r}^{\nu^i}\left(x_{0}\right) \) is the cube  of side \( 2r \) centered at \( x_0 \) oriented in the direction \( \nu^i \).  

 Up to a rigid motion, we can assume that \( x_0 = 0 \), $\nu^i=e_1$ and \( J^i =\{x_1 = 0\}\) in a small cube \( Q_{2r}=(-r,r)^n \).
By construction, for small \( \varepsilon \),
\(
u_\varepsilon(x) = g_\varepsilon(x_1)
\) in $Q_{2r}$,
where \( g_\varepsilon \) is the one-dimensional recovery sequence of the previous subsection for the jump from $u^i_-$ to $u^i_+$ (see \Cref{subsub:ub_1_1} and \Cref{shape}). From this, \eqref{tesi} follows easily by reducing the integral to one dimension and applying the one-dimensional \(\Gamma\)-\(\limsup\) inequality for \( g_\varepsilon \). For clarity, let us show the steps.
\begin{align*}
&\frac{1}{(2r)^{n-1}|\log\varepsilon|}
\iint_{Q_{2r}\times Q_{2r}}
\frac{|u_\varepsilon(x)-u_\varepsilon(y)|^2}
{|x-y|^{n+1}}\,dx\,dy\\
&=\frac{1}{(2r)^{n-1}|\log\varepsilon|}
\iint_{Q_{2r}\times Q_{2r}}
\frac{|g_\varepsilon(x_1)-g_\varepsilon(y_1)|^2}
{|x-y|^{n+1}}\,dx\,dy\\
&=\frac{1}{(2r)^{n-1}|\log\varepsilon|}
\int_{-r}^{r}\int_{-r}^{r}
\int_{(-r,r)^{n-1}}\int_{(-r,r)^{n-1}}
\frac{|g_\varepsilon(x_1)-g_\varepsilon(y_1)|^2}
{\left((x_1-y_1)^2+|x'-y'|^2\right)^{\frac{n+1}{2}}}
\,dx'\,dy'\,dx_1\,dy_1.
\end{align*}
\\\noindent
Now, note that, for each $a>0,$
\[
\begin{aligned}
\frac{1}{2 r} \int_{-r}^{r} \int_{-r}^{r} \left(a+|x_i - y_i|^2\right)^{-\frac{n+1}{2}}\, dx_i \, dy_i & =
\frac{1}{2 r} \int_{-2 r}^{2 r} 
\left(a+z_i^2\right)^{-\frac{n+1}{2}}\left(2 r - \left|z_i\right|\right) \, dz_i \\
& \leq \int_{-\infty}^{\infty} \left(a+z_i^2\right)^{-\frac{n+1}{2}} \, dz_i.
\end{aligned}
\]
Using this fact, we can estimate from above:
\begin{equation}\label{slicing}
\begin{aligned}
\frac{\mu_\varepsilon(Q_{2r}\times Q_{2r})}{(2 r)^{n-1}}=&\frac{1}{(2r)^{n-1}|\log\varepsilon|} \iint_{Q_{2r} \times Q_{2r}} \frac{ \left| u_{\varepsilon}(x) - u_{\varepsilon}(y) \right|^2}{|x - y|^{n+1}} \, dx \, dy \\
&\leq \frac{1}{|\log\varepsilon|}\int_{-r}^r \int_{-r}^r |g_\varepsilon(x_1) - g_\varepsilon(y_1)|^2 \left( \int_{\mathbb{R}^{n-1}} \frac{1}{\left( |x_1 - y_1|^2 + |z|^2 \right)^{\frac{n+1}{2}}} \, dz \right) dx_1 \, dy_1 \\
&= \frac{\omega_{n-1}}{|\log\varepsilon|} \int_{-r}^r \int_{-r}^r \frac{|g_\varepsilon(x_1) - g_\varepsilon(y_1)|^2}{|x_1 - y_1|^2} \, dx_1 \, dy_1.
\end{aligned}
\end{equation}
Thus, we finally get
\[
\frac{\lambda\left(Q_{2r}\right)}{(2r)^{n-1}} \leq \omega_{n-1} \limsup_{\varepsilon \to 0} \frac{1}{|\log \varepsilon|} \int_{-r}^r \int_{-r}^r
\frac{|g_\varepsilon(x_1) - g_\varepsilon(y_1)|^2}{|x_1 - y_1|^2} \, dx_1 \, dy_1\le\omega_{n-1} \Psi^{(1)}(u^i_-,u^i_+),
\]
where we used the one-dimensional upper bound inequality for $\{g_\varepsilon\}$ (see \Cref{subsub:ub_1_1}).
In conclusion, by the Besicovitch differentiation theorem,
$$
\lambda \le \omega_{n-1}\Psi^{(1)}(u^-,u^+)\mathcal{H}^{n-1}\mres J_u. 
$$
Now, since $\overline Q \times \overline Q\subset\subset Q'\times Q'$, by \textit{Step 2} we have
\begin{align*}
\limsup_{\varepsilon\to 0} F_\varepsilon(u_\varepsilon)
&\le \mu(\overline{Q}\times\overline{Q})
   + \limsup_{\varepsilon\to 0}\frac{1}{\varepsilon|\log\varepsilon|}\int_Q W_\varepsilon(u_\varepsilon)\,dx\\
&\le \int_{J_u\cap Q}\omega_{n-1}\,\Psi^{(1)}(u^-,u^+)\,d\mathcal{H}^{n-1}
   + \begin{cases}
       K\,\bigl|\{u\in S\}\bigr|, & \text{in \hyperlink{R1}{(R1)}},\\[2pt]
       0, & \text{in \hyperlink{R2}{(R2)}}.
     \end{cases}
\end{align*}
where we used that $\mathcal{H}^{n-1}(J_u \cap \partial Q)=0$. This is exactly $F(u)$, since $\omega_{n-1}\Psi^{(1)}$ is the jump cost for $s=\frac12.$
\end{proof}

\subsection{Upper bound for \texorpdfstring{$s>\frac12$}{s > 1/2}}\label{up_1}

\begin{proof}[Proof of \Cref{thm:main} ($s>\frac12$)] Here $\theta_\varepsilon=\varepsilon.$
We set 
\begin{equation}\label{T_eps}
T_\varepsilon := \min \left\{\varepsilon^{-1/4}, \delta_\varepsilon^{-1/2} \right\} \quad\text{and}\quad  
\eta_\varepsilon := \varepsilon T_\varepsilon,
\end{equation}
so that $T_\varepsilon \to +\infty$ and $\eta_\varepsilon \to 0$ for $\varepsilon \to 0$. Let $v^{T_\varepsilon}_i$ be a (non-decreasing) minimizer of the optimal profile  problem $\Theta^{T_\varepsilon}_{s,n}(u_-^i, u_+^i)$ as defined in \eqref{theta_s_T}, with $C_{s,n}$ in front of the seminorm. We recall that $v^{T_\varepsilon}_i$ is equal to $u_-^i$ in $(-\infty, -T_\varepsilon)$ and to $u_+^i$ in $(T_\varepsilon,+\infty)$ and, by \Cref{costo_theta}, $\Theta_{s,n}^{T_\varepsilon}(u_-^i, u_+^i)\rightarrow \Theta_{s,n}(u_-^i, u_+^i)$ for $\varepsilon\to 0.$

\medskip
\noindent\textit{Case $n \ge 2$}

\textit{Step 1: Construction of the recovery sequence.}
 For each simplex $J^i$, let $d^i(x)$ be the signed distance from the hyperplane containing  $J^i$, positive in the direction $\nu^i$. We define disjoint tubular neighborhoods of width $\eta_\varepsilon$:
$$
I^i_\varepsilon := \left\{ x \in    \tilde Q : |d^i(x)| \le 2 \eta_\varepsilon \text{ and } \operatorname{dist}(x,J^j) > 4 \eta_\varepsilon, \ \forall j \neq i \right\}, \quad
\tilde J^i_\varepsilon := \left\{ y \in J^i : y + t \nu^i \in I^i_\varepsilon, \ \forall |t|\le  2 \eta_\varepsilon \right\}, $$
$$
P^i_\varepsilon := \left\{ y + t \nu^i: \ y \in \tilde J^i_\varepsilon, \ |t| \le \eta_\varepsilon \right\},
$$
so that $\mathcal{H}^{n-1}(\tilde J^i_\varepsilon) \to \mathcal{H}^{n-1}(J^i)$ and the sets $P^i_\varepsilon + B_{\eta_\varepsilon}(0)$ are pairwise disjoint.
 We set 
\begin{equation}
    \tilde{u}_\varepsilon(x) := 
    \begin{cases} 
        v^{T_\varepsilon}_i \left( d^i(x)/\varepsilon \right) & \text{if } x \in P^i_\varepsilon, \\
        u(x) & \text{otherwise in } \tilde Q.
    \end{cases}
\end{equation}
The function $\tilde u_\varepsilon$ still has some jumps on $J_u\setminus \bigcup \tilde J^i_\varepsilon$ and across the lateral boundary of $P^i_\varepsilon.$ Now, notice that there exists $L\ge 3$, only depending on $u$, such that the jump set of $\tilde u_\varepsilon$ in $Q'$ is contained in $\{\operatorname{dist}(\cdot,\Sigma)<(L-2)\eta_\varepsilon\}.$
Then, to smooth the function near $\Sigma $, we introduce a cut-off function $\psi_\varepsilon \in C^\infty(\tilde Q;[0,1])$ such that $\psi_\varepsilon = 1$ in $\{\operatorname{dist}(\cdot,\Sigma) \le (L-1) \eta_\varepsilon\}$ and $\psi_\varepsilon = 0$ in $\{\operatorname{dist}(\cdot,\Sigma) \ge L \eta_\varepsilon\}$, with $|\nabla \psi_\varepsilon| \le C / \eta_\varepsilon$. The recovery sequence is
\begin{equation}
u_\varepsilon := \psi_\varepsilon (\tilde{u}_\varepsilon * \varphi_{\eta_\varepsilon}) + (1 - \psi_\varepsilon) \tilde{u}_\varepsilon,
\end{equation}
where $\varphi_{\eta_\varepsilon}$ is a standard mollifier supported in $B_{\eta_\varepsilon}(0)$. 
With this choice, $u_\varepsilon$ has no jumps and $\{u_\varepsilon\neq u\}\subset \mathcal{N}_{L\eta_\varepsilon}$, so $u_\varepsilon \to u$ in $L^1(Q).$
This construction provides a smooth transition in a neighborhood of $\Sigma$, while preserving the optimal profile elsewhere.

In the following, let $M:=\|u\|_\infty$ and $U :=\{\operatorname{dist}(\cdot,\Sigma)<(L+1)\eta_\varepsilon\}$, so that $|U|\le C\eta_\varepsilon^2$ by \eqref{stima_N}.

\textit{Step 2: Estimate of the seminorm.}
Observe that if $|z|\ge \eta_\varepsilon $ and $x\notin \mathcal{N}_{(L+1)|z|}$, then
\(
u_\varepsilon(x+z) = u_\varepsilon(x),
\)
since $u_\varepsilon=u$ outside $\mathcal{N}_{L\eta_\varepsilon}.$
Hence, by \eqref{stima_N}, the long-range interaction vanishes: we have
\begin{equation}
\begin{aligned}
\varepsilon^{2s-1}\int_{\eta_\varepsilon\leq |z| \leq \operatorname{diam}(Q')} \frac{dz}{|z|^{n+2s}} &\int_{Q\cap (Q-z)} \Bigl|u_\varepsilon(x+z)-u_\varepsilon(x)\Bigr|^{2} \,dx 
\\&\leq\varepsilon^{2s-1} C \int_{\eta_\varepsilon\leq |z| \leq \operatorname{diam}(Q')} \frac{dz}{|z|^{n+2s-1}}=CT_\varepsilon^{1-2s}\to 0.
\end{aligned}
\end{equation}
For the short-range part ($|x-y|\le \eta_\varepsilon$), we distinguish two cases: when $x,y \notin U$ and when $x\in U.$
Let $x, y \notin U$ with $|x-y|\le \eta_\varepsilon.$ If $x \notin \mathcal{N}_{2\eta_\varepsilon}$, the segment $[x,y]$ does not meet the layers, so $u_\varepsilon(x)=u_\varepsilon(y).$
Otherwise, if $x\in \mathcal{N}_{2\eta_\varepsilon}$, there exists $i$ such that $x,y\in C_i:=\left\{ p + t \nu^i\,:\, \ p \in \tilde J^i_\varepsilon, \ |t| \le 3\eta_\varepsilon \right\}.$ In this set, $u_\varepsilon=v_i^{T_\varepsilon}(d^i/\varepsilon)$. Hence, we have 
\begin{equation}\label{5}
\begin{aligned}
\varepsilon ^{2s-1}&\iint_{(Q\setminus U)^2\cap \{|x-y|\le \eta_\varepsilon\}}
\frac{|u_\varepsilon(x)-u_\varepsilon(y)|^2}{|x-y|^{n+2s}}
\,dx\,dy
\\&\le\varepsilon ^{2s-1}
\sum_j
\int_{J^j}
\int_{-3\eta_\varepsilon}^{3\eta_\varepsilon}
\int_{J^j}
\int_{-3\eta_\varepsilon}^{3\eta_\varepsilon}
\frac{
\left|
v_j^{T_\varepsilon}\!\left(\frac{z}{\varepsilon}\right)
-
v_j^{T_\varepsilon}\!\left(\frac{w}{\varepsilon}\right)
\right|^2
}{
\bigl(
|\sigma-\tau|^2 + (z-w)^2
\bigr)^{\frac{n+2s}{2}}
}
\, d\sigma\,dz\,d\tau\,dw 
\\& \leq
 \sum_j \mathcal{H}^{n-1}(J^j) \int_{\mathbb{R}} \int_{\mathbb{R}} \left| v^{T_\varepsilon}_j(\xi) - v^{T_\varepsilon}_j(\zeta) \right|^2 \left( \int_{\mathbb{R}^{n-1}} \frac{d\theta}{\big( |\theta|^2 + (\xi - \zeta)^2 \big)^{\frac{n+2s}{2}}} \right) d\xi \, d\zeta,
\end{aligned}
\end{equation}
where we used the change of variables \[
    \xi = \frac{z}{\varepsilon }, \quad \zeta = \frac{w}{\varepsilon }, \quad \theta = \frac{\sigma - \tau}{\varepsilon }.
    \]
    Since $x,y\in Q$, in \eqref{5} the variables $\sigma,\tau$ range only in $J^j\cap(Q+B_{3\eta_\varepsilon})$; therefore, in the following, $\mathcal{H}^{n-1}(J^j)$ stands for $\mathcal{H}^{n-1}(J^j\cap(Q+B_{3\eta_\varepsilon}))$, which tends to $\mathcal{H}^{n-1}(J^j\cap\overline Q)=\mathcal{H}^{n-1}(J^j\cap Q)$ as $\varepsilon\to0$.
Now, notice that 
\begin{equation}\label{6}
\int_{\mathbb{R}^{n-1}}
\frac{1}{\big( |\theta|^2 + A^2 \big)^{\frac{n+2s}{2}}}
\, d\theta
=
\frac{1}{A^{n+2s}}
\int_{\mathbb{R}^{n-1}}
\frac{A^{n-1}}{( |\omega|^2 + 1 )^{\frac{n+2s}{2}}}
\, d\omega
=
\frac{C_{s,n}}{A^{1+2s}}.
\end{equation}
 Then, 
 $$
 \varepsilon ^{2s-1}\iint_{(Q\setminus U)^2\cap \{|x-y|\le \eta_\varepsilon\}}
\frac{|u_\varepsilon(x)-u_\varepsilon(y)|^2}{|x-y|^{n+2s}}
\,dx\,dy\le \sum_j \mathcal{H}^{n-1}(J^j)C_{s,n}[v_j^{T_\varepsilon}]^2_{H^s(\R)}.
 $$
Now, let $x \in U$ with $|x-y|\le \eta_\varepsilon.$ Let $h_\varepsilon:=\tilde u_\varepsilon\,*\,\varphi_{\eta_\varepsilon}$, which is $C/\eta_\varepsilon$-Lipschitz. By definition, we have
\begin{align}
u_\varepsilon(x)-u_\varepsilon(y)&=\psi_\varepsilon(x)(h_\varepsilon(x)-h_\varepsilon(y))+\\&+(\psi_\varepsilon(x)-
\psi_\varepsilon(y))
(h_\varepsilon(y)-\tilde u_\varepsilon(y))+(1-\psi_\varepsilon(x))(\tilde u_\varepsilon(x)-\tilde u_\varepsilon(y)).
\end{align}
The first two terms are bounded by $C|x-y|/\eta_\varepsilon$, hence their contribution vanishes; indeed
$$
C\frac{\varepsilon^{2s-1}|U|}{\eta^2_\varepsilon}\int_{|z|\le \eta_\varepsilon}\frac{dz}{|z|^{n+2s-2}}\le C \varepsilon T_\varepsilon^{2-2s}\rightarrow 0.
$$
The third term vanishes in $\{\operatorname{dist}(x, \Sigma)\le (L-1)\eta_\varepsilon\}$; otherwise the segment $[x,y]$ does not meet the jumps of $\tilde u_\varepsilon$. Hence we can repeat the same computation as in \eqref{5} and conclude that the energy contribution is bounded from above by
$$
C\mathcal{H}^{n-1}(J_u \cap \{\operatorname{dist}(x, \Sigma)< (L+2)\eta_\varepsilon\})\le C \eta_\varepsilon\rightarrow 0.
$$
In conclusion,
$$
\varepsilon^{2s-1}[u_\varepsilon]^2_{H^s(Q)}\le  \sum_{j=1}^m \mathcal{H}^{n-1}(J^j)C_{s,n}[v_j^{T_\varepsilon}]^2_{H^s(\R)}+o(1).
$$

\medskip

\textit{Step 3: Estimate of the potential term.}
By the upper estimate \eqref{eq:two_sided_estimate}, $W_\varepsilon\le \operatorname{dist}^2(\cdot,\Z\cup S)+\delta_\varepsilon\chi_{A_1}$. Moreover, outside $U$
we use the change of variables $t = d^i(x)/\varepsilon$ and we obtain
\begin{equation}
    \frac{1}{\varepsilon} \int_{Q} W_\varepsilon(u_\varepsilon) \, dx \leq \sum_{i=1}^m \mathcal{H}^{n-1}(J^i) \int_{-T_\varepsilon}^{T_\varepsilon} W^S(v_i^{T_\varepsilon}(t)) \, dt + \frac{\delta_\varepsilon}{\varepsilon} |A_{1}(u_\varepsilon)| + \frac{C|U|}{\varepsilon},
\end{equation}
where $C|U|/\varepsilon\le C\eta^2_\varepsilon /\varepsilon\leq C\varepsilon^{1/2}\to 0.$ In regime \hyperlink{R1}{(R1)}, the second term converges to $K|\{u \in S\}|.$ In regime \hyperlink{R2}{(R2)}, $u$ takes values in $\Z$, so $|A_1(u_\varepsilon)|\le |\mathcal{N}_{L\eta_\varepsilon}|\le C\eta_\varepsilon$, and the second term is bounded by $C\delta_\varepsilon^{1/2}\to 0.$

\textit{Step 4: Conclusion.}
Combining  \textit{Step 2} and \textit{Step 3}, and since 
$$
C_{s,n}[v_i^{T_\varepsilon}]^2_{H^s(\R)}+\int_\R\operatorname{dist}^2(v_i^{T_\varepsilon},\Z\cup S)=\Theta^{T_\varepsilon}_{s,n}(u_-^i,u_+^i), 
$$
 we have
\begin{equation}
\begin{aligned}
\limsup_{\varepsilon\to 0}F_\varepsilon(u_\varepsilon)\le
\sum_{i=1}^m
\mathcal{H}^{n-1}(J^i)\,
\Theta_{s,n}(u_-^i,u_+^i)
+
 \begin{cases} 
        K |\{u \in S\}|  & \text{in \hyperlink{R1}{(R1)}} , \\
       0 & \text{in \hyperlink{R2}{(R2)} },
    \end{cases}
\end{aligned}
\end{equation}
which is exactly $F(u)$, since, by \Cref{costo_theta}, $\Theta_{s,n}(u_-^i,u_+^i)=\lambda_{s,n}\kappa_s\Phi_{\Z\cup S}(u_-^i,u_+^i)$.

\medskip\noindent
\textit{Case $n=1$.} 
\\
Now $J_u=\{x_1, \dots, x_m\}$ and $\nu^i\in \{\pm1\}$ is the orientation of the jump. We define the recovery sequence as
\begin{equation}
    u_\varepsilon(x) := 
    \begin{cases} 
        v^{T_\varepsilon}_i \left( \frac{\nu^i(x - x_i)}{\varepsilon} \right), & \text{if } x \in [x_i - \eta_\varepsilon, x_i + \eta_\varepsilon]=:P^i_\varepsilon, \\
        u(x), & \text{otherwise in } Q,
    \end{cases}
\end{equation}
where $\eta_\varepsilon, T_\varepsilon$ are defined in \eqref{T_eps}. For $\varepsilon$ sufficiently small, the intervals $[x_i - \eta_\varepsilon, x_i + \eta_\varepsilon]$ are mutually disjoint. By the change of variables $\xi = (x-x_i)/\varepsilon$ and $\zeta = (y-x_i)/\varepsilon$, we have:
\begin{equation}
    \varepsilon^{2s-1} \int_{P^i_\varepsilon} \int_{P^i_\varepsilon} \frac{|u_\varepsilon(x) - u_\varepsilon(y)|^2}{|x - y|^{1+2s}} \, dx \, dy = \int_{-T_\varepsilon}^{T_\varepsilon} \int_{-T_\varepsilon}^{T_\varepsilon} \frac{|v_i^{T_\varepsilon}(\xi) - v_i^{T_\varepsilon}(\zeta)|^2}{|\xi - \zeta|^{1+2s}} \, d\xi \, d\zeta.
\end{equation}
\noindent
 Similarly, the potential term scales as
\begin{equation}
    \frac{1}{\varepsilon} \int_{P^i_\varepsilon} W^S(u_\varepsilon(x)) \, dx = \int_{-T_\varepsilon}^{T_\varepsilon} W^S(v_i^{T_\varepsilon}(\xi)) \, d\xi.
\end{equation}
Since the interactions between different jumps contribute to the energy as $o(1)$, we have
\begin{equation}
    \limsup_{\varepsilon \to 0} F_\varepsilon(u_\varepsilon) \leq \sum_{i=1}^m \Theta_{s,1}(u_-^i, u_+^i) + 
    \begin{cases} 
        K |\{u \in S\}| &  \text{in \hyperlink{R1}{(R1)}}, \\
        0 &  \text{in \hyperlink{R2}{(R2)}}.
    \end{cases}
\end{equation}
This concludes the proof.

\end{proof}

\section{Limiting behavior for \texorpdfstring{$s < \frac12$}{s < 1/2}}\label{s<1/2}

In the low-order regime $s < \frac12$, the characteristic feature is that bounded $BV$ functions belong to $H^s(Q)$ with $[u]^2_{H^s(Q)}\le C\|u\|_\infty|D u|(Q)$. 
\\\noindent
Assuming the scaling parameter $\theta_\varepsilon = \varepsilon^{2s}$ as in \cite{SavinValdinoci2010}, the functional \eqref{eq:main_functional} simplifies to
\begin{equation}\label{functiona_piccolo_s}
F_\varepsilon(u; Q) = \int_{Q} \int_{Q} \frac{|u(x) - u(y)|^2}{|x - y|^{n + 2s}} \, dx \, dy + \frac{1}{\varepsilon^{2s}} \int_{Q} W_\varepsilon(u(x)) \, dx.
\end{equation}
Also in this case, the limiting behavior is determined by the value of 
$
K:=\lim_{\varepsilon\to 0}\delta_\varepsilon/\theta_\varepsilon\in[0,+\infty].
$
Unlike the cases $s \ge \frac12$, the limit functional for $s < \frac12$ is non-local and does not concentrate on the jump set $J_u$. The scaling $\theta_\varepsilon = \varepsilon^{2s}$ is the natural one for this regime since it preserves the fractional seminorm, while the potential term acts as a penalization for intermediate values, either allowing them (if $K < \infty$) or forbidding them (if $K = \infty$).

For a closed set $\mathcal{V}\subset \R,$ we set
\begin{equation}\label{Hs}
    H^s(Q;\mathcal{V}):=\{u \in L^1(Q):[u]_{H^s(Q)}<+\infty, u(x)\in \mathcal{V}\text{ for a.e. }x\in Q\}.
\end{equation}
\begin{theo}[Compactness and $\Gamma$-convergence for $s < \frac12$]\label{thm:s_less_half}
Let $s \in \left(0,\frac12\right)$ and consider the sequence of functionals $\{F_\varepsilon\}$  defined in \eqref{functiona_piccolo_s}.
\begin{enumerate}[label=(\roman*)]
    \item \textbf{Compactness:} Let $\{u_\varepsilon\} \subset L^1(Q)$ be a sequence such that $\sup_\varepsilon F_\varepsilon(u_\varepsilon) \leq M < +\infty$. Then there exists a sequence of integers $\{a_\varepsilon\} \subset \mathbb{Z}$ such that the translated sequence
$\{u_\varepsilon - a_\varepsilon\}$ is relatively compact in $L^1(Q)$. If, in addition, $u_\varepsilon\to u$ in $L^1(Q)$ then $u\in H^s(Q; \union)$. Moreover, if $\delta_\varepsilon / \varepsilon^{2s} \to +\infty$, then $u\in H^s(Q; \Z).$
    \item \textbf{$\Gamma$-convergence:} 
    \begin{itemize}
        \item If $\frac{\delta_\varepsilon}{\varepsilon^{2s}} \to K \in [0, +\infty)$, then $F_\varepsilon$ $\Gamma$-converges in $L^1(Q)$ to:
        \begin{equation}
            G_{\mathrm{sup}}(u) = \begin{cases} 
            [u]^2_{H^s(Q)} + K |\{x \in Q : u(x) \in S\}| & \text{if } u \in H^s(Q;\union), \\
            +\infty & \text{otherwise.}
            \end{cases}
        \end{equation}
        \item If $\frac{\delta_\varepsilon}{\varepsilon^{2s}} \to +\infty$, then $F_\varepsilon$ $\Gamma$-converges in $L^1(Q)$ to:
        \begin{equation}
            G_{\mathrm{sub}}(u) = \begin{cases} 
            [u]^2_{H^s(Q)} & \text{if } u \in H^s(Q;\Z), \\
            +\infty & \text{otherwise.}
            \end{cases}
        \end{equation}
    \end{itemize}
\end{enumerate}
\end{theo}

\begin{proof}\textit{(i)}
To prove compactness, we observe that the uniform bound on the $H^s(Q)$ seminorm leads to a uniform bound in $H^s(Q)$ for $u_\varepsilon
$ minus its mean, due to the fractional Poincaré inequality. 
Then by the compact embedding of $H^s(Q)$ into $L^1(Q)$, the thesis follows. More precisely,
from the uniform energy bound
\(
\sup_{\varepsilon>0} F_\varepsilon(u_\varepsilon)\le M,
\)
we immediately obtain the estimates
\begin{equation}\label{comp_s_piccolo}
[u_\varepsilon]_{H^s(Q)}^2\le M\quad \hbox{and}
\quad
\int_Q \operatorname{dist}^2(u_\varepsilon,\mathbb{Z}\cup S)\,dx
\le
\int_Q W_\varepsilon(u_\varepsilon)\,dx
\le
M\varepsilon^{2s}\longrightarrow0.
\end{equation}
By the fractional Poincaré--Wirtinger inequality on the bounded connected domain \(Q\) (see \cite[Theorem 6.33]{Leoni}), we have 
\[
\|u_\varepsilon-\bar{u}_\varepsilon\|_{L^2(Q)}
\le
C [u_\varepsilon]_{H^s(Q)}
\le
C\sqrt{M}, \quad \hbox{where}\quad 
\bar{u}_\varepsilon:=\fint_Q u_\varepsilon.
\]
Since
$[u_\varepsilon-\bar{u}_\varepsilon]_{H^s(Q)}
=
[u_\varepsilon]_{H^s(Q)},$
the sequence \(\{u_\varepsilon-\bar{u}_\varepsilon\}\) is uniformly bounded in \(H^s(Q)\). Since the embedding
\(
H^s(Q)\hookrightarrow L^2(Q)
\)
is compact (see \cite[Theorem 7.1]{Frazionari}), up to a subsequence,
\[
u_\varepsilon-\bar{u}_\varepsilon\rightarrow v
\quad\text{strongly in }L^2(Q).
\]
Now, let \(a_\varepsilon\in\mathbb Z\) be an
integer closest to \(\bar u_\varepsilon\), so that
\(
\bar u_\varepsilon-a_\varepsilon
\in
\left[-\frac12,\frac12\right].
\)
Then, up to a subsequence, 
\(
\bar u_\varepsilon-a_\varepsilon\longrightarrow c,
\)
for some \(c\in\left[-\frac12,\frac12\right]\). Consequently,
\[
u_\varepsilon-a_\varepsilon
\longrightarrow
v+c
\quad\text{strongly in }L^2(Q).
\]
Therefore, the sequence \(\{u_\varepsilon-a_\varepsilon\}\) is relatively compact in
\(L^1(Q)\).
If $u_\varepsilon\to u$, then by \eqref{comp_s_piccolo} for $\theta_\varepsilon=\varepsilon^{2s}$, we have
\[
u(x)\in\mathbb Z\cup S
\quad\text{for a.e. }x\in Q.
\]
Finally, by the lower semicontinuity of the fractional seminorm,
\[
[u]_{H^s(Q)}
\le
\liminf_{\varepsilon\to0}
[u_\varepsilon]_{H^s(Q)}
<
+\infty.
\]
 Moreover, by \eqref{eq:two_sided_estimate}, we have that if $\delta_\varepsilon/\varepsilon^{2s}\rightarrow +\infty$ then 
$|\{u\in S\}|=0$ by the uniform energy bound. This proves the desired compactness result.

Regarding the $\Gamma$-convergence, the $\Gamma$-liminf is a direct consequence of \eqref{eq:two_sided_estimate} and the lower semicontinuity of the $H^s$-seminorm with respect to $L^1$-convergence; if $\delta_\varepsilon/\varepsilon^{2s}\to K <+\infty$, the bulk term follows by convergence in measure, i.e. $\liminf_\varepsilon|A_\delta(u_\varepsilon)|\ge |\{u \in S\}|$ with $A_\delta$ the one of \eqref{eq:two_sided_estimate}. For the $\Gamma$-limsup, given $u \in H^s(Q; \union)$, the constant sequence $u_\varepsilon = u$ is a recovery sequence.
\end{proof}
\addtocontents{toc}{\protect\setcounter{tocdepth}{0}}
\section*{Use of AI}
NO AI tools were used to produce mathematical content; the strategy and
the proofs were developed and written exclusively by the authors.
Claude Opus 5.5 (Anthropic) was used only for minor language editing and to
detect typographical and LaTeX errors on the very last version of the manuscript.

\section*{Acknowledgments}
A.G. and M.A.P. are members of the Gruppo Nazionale per l'Analisi
Matematica, la Probabilità e le loro Applicazioni (GNAMPA) of the
Istituto Nazionale di Alta Matematica (INdAM).

\addtocontents{toc}{\protect\setcounter{tocdepth}{1}}

\printbibliography[title={References}]

\end{document}